\documentclass{amsart}
\usepackage{amscd,amssymb}
\newtheorem{theorem}{Theorem}[section]
\newtheorem{corollary}[theorem]{Corollary}
\newtheorem{example}[theorem]{Example}
\newtheorem{proposition}[theorem]{Proposition}
\newtheorem{lemma}[theorem]{Lemma}
\newtheorem{criterion}[theorem]{Criterion}
\newtheorem{remark}[theorem]{Remark}
\theoremstyle{definition}
\newtheorem*{definition}{Definition}
\newtheorem*{example*}{Example}
\newtheorem*{remark*}{Remark}
\newtheorem*{problem*}{Problem}
\newtheorem*{notation*}{Notation}
\def\mymid{{\restriction}}
\def\bplus{\mathbin{[+]}}
\def\btimes{\mathbin{[\times]}}
\def\Ext{\operatorname{Ext}}
\def\Tor{\operatorname{Tor}}
\def\Hom{\operatorname{Hom}}
\def\Star{\operatorname{Star}}
\def\dim{\operatorname{dim}}
\def\support{\operatorname{support}}
\def\Im{\operatorname{Im}}
\def\Cl{\operatorname{Cl}}
\def\pr{\operatorname{pr}}
\def\ab{\operatorname{ab}}
\def\id{\operatorname{id}}
\def\Int{\operatorname{Int}}
\def\mesh{\operatorname{mesh}}
\def\diam{\operatorname{diam}}
\def\cone{\operatorname{cone}}
\def\R{\mathbb{R}}
\def\Q{\mathbb{Q}}
\def\Z{\mathbb{Z}}
\def\N{\mathbb{N}}
\def\A{\mathbb{A}}
\def\C{\mathbb{C}}
\def\E{\mathbb{E}}
\def\sD{\mathcal{D}}
\def\sF{\mathcal{F}}

\def\sL{\mathcal{L}}

\def\sM{\mathcal{M}}
\def\sP{\mathcal{P}}

\def\sS{\mathcal{S}}
\def\sU{\mathcal{U}}
\def\sV{\mathcal{V}}

\def\sR{\mathcal{R}}
\def\sB{\mathcal{B}}
\def\sA{\mathcal{A}}
\def\sC{\mathcal{C}}
\def\varinjLim{\mathop{\underset{\textstyle\rightarrow}{\text{Lim}}}}
\def\varprojLim{\mathop{\underset{\textstyle\leftarrow}{\text{Lim}}}}
\def\limto{\varinjlim}
\def\limleftarrow{\varprojlim}
\def\Limto{\varinjLim}
\def\Limleftarrow{\varprojLim}
\def\Limone{\operatorname{Lim}^1}
\begin{document}
{\noindent\small
Topology Atlas Invited Contributions {\bf 6} no.~3 (2001) 61 pp.
\vskip \baselineskip}
\author{A.N. Dranishnikov}
\title{Cohomological dimension theory of compact metric spaces}
\address{Department of Mathematics,
University of Florida,
444 Little Hall,
Gaines\-ville, FL 32611-8105}
\email{dranish@math.ufl.edu}
\urladdr{http://www.math.ufl.edu/$\tilde{~}$dranish/}
\thanks{Beverley L.~Brechner asked Alexander Dranishnikov to write this
survey for Topology Atlas. It was written in 1998.
This survey has not been refereed.}
\maketitle

\tableofcontents
\setcounter{section}{-1}
\section{Introduction}

The cohomological dimension theory has connections with many different 
areas of mathematics: dimension theory, topology of manifolds, group 
theory, functional rings and others. It was founded by P.S.~Alexandroff in 
late 20's. Many famous topologists have contributed to the theory. Among 
them are Hopf, Pontryagin, Bockstein, Borsuk, Dyer, Boltyanskij, Kodama, 
Kuzminov, Sitnikov.

There are only few introductory and survey texts on the theory. The book 
by Alexandroff `Introduction to homological dimension theory and general 
combinatorial topology' \cite{A1} is written in old fashion language and 
hardly readable. There are surveys by Kodama (Appendix in \cite{Na}) and 
by Kuzminov \cite{Ku}. The first part of Kuzminov paper is devoted to 
compact metric spaces and is an excellent reading. We don't consider 
noncompact spaces in this paper, since the cohomological dimension of 
noncompact spaces behaves differently and is not completely developed. A 
very special case of the cohomological dimension theory is the case of 
integer coefficients. An excellent survey on this case was written by 
Walsh \cite{WA} where a detailed proof of the Edwards resolution theorem 
first was published. In 1988, twenty years after Kuzminov's survey I wrote 
a sequel to that \cite{Dr1}. Since then ten years passed, new results 
appeared and a new understanding of the old results ripened. So time came 
for an update survey. A new compressed survey was given by Dydak 
\cite{Dy3} where the main applications of the cohomological dimensions are 
discussed. Here we present a detailed introductory survey of the theory.

This survey and Dydak's have the same origin. They appeared as the notes 
to our joint book that we planned to write \cite{D-D-W}. We still have a 
hope that someday we will accomplish that.

In this paper we assume that the reader is familiar with basic elements of 
the homotopy theory, with homology and cohomology theories, including the 
\v{C}ech cohomology, the Steenrod homology and extraordinary 
(co)homologies. Some knowledge in the dimension theory and the theory of 
absolute neighborhood retracts will be useful. Also we don't discuss here 
any applications of the cohomological dimension theory even to the 
dimension theory. Interested reader can find a discussion of some 
applications in \cite{Dy3}.

I am thankful to Topology Atlas for inviting me to write this survey.
I am also thankful to NSF, DMS-9971709, for the support.

\section{General properties of the cohomological dimension}

We define {\it the cohomological dimension} with respect to an abelian 
group $G$ of a topological space $X$ as the largest number $n$ such that 
there exists a closed subset $A\subset X$ with $\check H^n(X,A;G)\neq 0$. 
We denote it by $\dim_G X=n$. If there is no such number we set 
$\dim_G X=\infty$. This definition is good for any space. We restrict 
ourselves by compact metric spaces (we call them compacta). Actually 
everywhere in this paper one can replace compact spaces by 
$\sigma$-compact, i.e.\ countable unions of compacta.

\begin{theorem}
For any compactum $X$ and an abelian group $G$ the following conditions 
are equivalent
\begin{enumerate}
\item $\dim_GX\leq n$,
\item $\check H^{n+1}(X,A;G)=0$ for all closed $A\subset X$,
\item $H_c^{n+1}(U;G)=0$ for all open $U\subset X$,
\item for every closed subset $A\subset X$ the inclusion homomorphism 
$\check H^n(X;G)\to\check H^n(A;G)$ is an epimorphism,
\item $K(G,n)$ is an absolute extensor for $X$, $K(G,n)\in AE(X)$, i.e.\
every continuous map $f\colon A\to K(G,n)$ of a closed subset 
$A\subset X$ has a continuous extension over $X$.
\end{enumerate}
\end{theorem}

\begin{proof}
The implication (1) $\Rightarrow$ (2) follows from the definition.

The condition (3) equals (2) by virtue of the equality
$\check H^k(X,A;G)=H_c^k(X\setminus A;G)$.
The implication (2) $\Rightarrow$ (4) follows from the long exact sequence
of the pair $(X,A)$.

The conditions (4) and (5) are equivalent since (5) is (4) 
formulated in homotopy language.

To show (5) $\Rightarrow$ (1) first we prove that (5)=$(5)_n$ implies
$(5)_k$ for all $k\geq n$.
We consider a Serre fibration $p\colon E\to K(G,n+1)$ where $E$ is 
contractible, $K(G,n+1)$ is a simplicial complex representing the 
Eilenberg-MacLane space.
Then the homotopy fiber of $p$ is $K(G,n)$. Then 
$p^{-1}(\Delta)\in AE(X)$ for any simplex $\Delta$ (see Chp~1).
Then $p^{-1}(\Delta)\in AE(A)$ for any closed $A\subset X$.
Then any map $g\colon A\to K(G,n)$ there is a homotopy lift $\bar 
g\colon A\to E$.
Since $E$ is contractible, the map $\bar g$ and, hence, $g$ is 
homotopically trivial. 
Therefore $g$ is extendible over $X$. 
Thus, we proved that $(5)_n$ implies $(5)_{n+1}$. 
By induction we can prove all $(5)_k$ for $k\geq n$. 
If $K(G,k)\in AE(X)$, then $K(G,k)\in 
AE(X/A)$ for any closed subset $A\subset X$. Let $k>n$.
Then any map $f\colon X/A\to K(G,k)$ 
can be lifted to a map $\bar f\colon X/A\to E$. Since $\bar f$ is 
null homotopic, the map $f$ is null homotopic. Since $f$ is
arbitrary, we have $\check H^{k}(X/A;G)=0$. Since we have that for 
all $k>n$, $\dim_GX<n+1$
and (1) is proven.
\end{proof}

The property (5) automatically implies:

\begin{corollary}
For every closed subset $A\subset X$ there is the inequality;
$\dim_GA\leq \dim_GX$ for any $G$.
\end{corollary}

\begin{example}
\mbox{}
\begin{enumerate}
\item $\dim_GX\leq \dim_{\Z}X\leq \dim X$ for any abelian group $G$ and
compact space $X$.
\item $\dim_GX=0$ if and only if $\dim X=0$ for any (nontrivial) $G$.
\item $\dim_{\Z}X=1$ if and only if $\dim X=1$.
\item $\dim_GK=n$ for every $n$-dimensional polyhedron $K$ and any $G\neq 0$.
\end{enumerate}
\end{example}

\begin{proof}
(1). 
The first inequality follows from the Universal Coefficient Formula. 
The second inequality can be rewritten as an implication 
$S^n\in AE(X) \Rightarrow K(\Z,n)\in AE(X)$ which follows from the
fact that $S^n$ is an $n$-skeleton of $K(\Z,n)$ and the standard
homotopy theory.

(2). The space $K(G,0)$ contains $S^0$ as a retract.

(3). $S^1\in K(\Z,1)$.

(4). 
By (1) $\dim_GK\leq n$. Since $K$ contains an open set $U$ homeomorphic 
to $\R^n$, $H^n_c(U;G)=G\neq 0$ and Theorem 1.1 implies the inequality 
$\dim_GK\geq n$.
\end{proof}

\begin{theorem}[Alexandroff Theorem]
For finite dimensional compacta there is the equality
$\dim_{\Z}X=\dim X$. 
\end{theorem}

\begin{proof}
In view of 1.3 (1) it suffices to show that $\dim_{\Z}X\geq \dim X$. 
Assume the contrary: $\dim X=n$ and $\dim_{\Z}X\leq n-1$. Take an 
Eilenberg-MacLane complex $K=K(\Z,n-1)$ such that its $n$-dimensional 
skeleton $K^{(n)}$ is an $n{-}1$-sphere $S^{n-1}$. Show that $S^{n-1}\in 
AE(X)$. Take a continuous map $f\colon A\to S^{n-1}=K^{(n)}$ of a closed
subset $A\subset X$. By Theorem 1.1 there is a continuous extension 
$\bar f\colon X\to K$. Since the dimension of $X$ is $\leq n$, by the 
Cellular Approximation theorem there is a homotopy $H_t\colon X\to K$ such 
that 
\begin{itemize}
\item $H_0=\bar f$, 
\item $H_1(X)\subset K^{(n)}$ and
\item $H_t\mymid_A=f$ for all $t\in[0,1]$. 
\end{itemize}

Hence, $H_1\colon X\to S^{n-1}$ is an extension of $f$.
Thus, $S^{n-1}\in AE(X)$ and hence, $\dim X\leq n-1$. Contradiction.
\end{proof}

\begin{theorem}[Countable Union Theorem]
Suppose $X = \bigcup X_i$ and each $X_i$ is a compactum. Then 
$\dim_GX = \sup\{\dim_GX_i\}$.
\end{theorem}

\begin{proof}
If the family $\{\dim_GX_i\}$ is bounded, then the formula holds by the 
trivial reason. Now, we show that if all $\dim_GX_i\leq n$, then 
$\dim_GX\leq n$. We show that $K(G,n)\in AE(X)$. Although $X$ is not 
compact, this condition implies the inequality $\dim_GX\leq n$. Let 
$f\colon A\to K(G,n)$ be a continuous map of a closed subset $A\subset X$. 
We define a nested increasing sequence of open in $X$ sets $U_1\subset 
\Cl(U_1)\subset U_2\subset \Cl(U_2)\cdots$ and a sequence of maps 
$f_i\colon \Cl(U_i)\to K(G,n)$ such that 
\begin{itemize}
\item $X = \bigcup_{i=1}^{\infty}U_i$, 
\item $A\subset U_1$ and $f_1\mymid_A=f$,
\item $f_{i+1}\mymid_{U_i}=f_i$ for all $i$. 
\end{itemize}
Then such a sequence defines a continuous map 
$\bigcup_{i=1}^{\infty}\colon X\to K(G,n)$ which is an extension of $f$.

Do it by induction on $i$. Extend the map $f$ over an open neighborhood 
$V\supset A$ to a map $f_1'\colon V\to K(G,n)$. Take an open set $U_1$ 
such that $A\subset U_1\subset \Cl(U_1)\subset V$ and define 
$f_1 = f_1'\mymid_{\Cl(U_1)}$. To define $U_{k+1}$ and $f_{k+1}$ we 
extend a map $f_k$ restricted on $\Cl(U_k)\cap X_k$ over a space $X_k$ to 
a map $g_k\colon X_k\to K(G,n)$. Then the union of $f_k$ and $g_k$ 
defines a continuous map $q_k\colon \Cl(U_k) \cup X_k\to K(G,n)$. Extend 
that map over a neighborhood $V_{k+1}$ to a map $f_{k+1}'$ and define 
$U_{k+1}\supset \Cl(U_k) \cup X_k$ such that its closure lies in 
$V_{k+1}$. 
Define $f_{k+1}$ as the restriction of $f_{k+1}'$ onto $U_{k+1}$.
\end{proof}

\begin{theorem}
Let $G=\Limto G_i$ and $\dim_{G_i}X\leq n$.
Then $\dim_GX\leq n$.
\end{theorem}

\begin{proof}
The formula $\Limto H^n_c(U;G)=H^n_c(U;\Limto G_i)$ implies the proof.
\end{proof}

\begin{corollary}
If $G = \bigoplus G_s$, then for every compactum $X$ the following
formula holds $\dim_GX = \sup\{\dim_{G_s}X\}$.
\end{corollary}

\begin{proof}
Since $H^n_c(U;G_s \oplus G') = H^n_c(U;G_s) \oplus H^n_c(U;G')$, the 
inequality $\dim_G X \geq \dim_{G_s} X$ holds. 
Hence, $\dim_G X \geq \sup\{\dim_{G_s} X\}$. 
The opposite inequality follows from Theorem~1.6 applied to 
$G=\Limto \bigoplus_{s=1}^i G_s$ and the fact that
$\sup_i\{\dim_{\bigoplus_{s=1}^i\! G_s} X\} = 
\sup_s\{\dim_{G_s} X\}$
imply the proof.
\end{proof}

\begin{definition}
A compactum $X$ has an {\it $r$-dimensional obstruction} at its point $x$ 
with respect to a coefficient group $G$ if there is a neighborhood $U$ of 
$x$ such that for every smaller neighborhood $V$ of $x$ the image of the 
inclusion homomorphism $i_{V,U}\colon H^r_c(V;G)\to H^r_c(U;G)$ is 
nonzero.
\end{definition}

\begin{theorem}
Let $X$ be a compact with $\dim_GX=r$ then $X$ contains a compact subset 
$Y$ of $\dim_GY=r$ such that at every point $x\in Y$ the compact $X$ has 
an $r$-dimensional obstruction with respect to $G$.
\end{theorem}

\begin{proof}
Let $W$ be an open subset of $X$ with $H_c^r(W;G)\neq 0$. Because 
of the continuity of cohomology there is a closed in $U$
set $Z$ minimal with respect the property: the inclusion 
homomorphism $H^r_c(W;G)\to H^r_c(Z;G)$ is nonzero. Then
$\dim_GZ=r$ and by the Countable Union Theorem there exists a 
compact subset $Y\subset Z$ with $\dim_GY=r$. For every $x\in Y$ we
take $U=W$. Let $V\subset U$ be a neighborhood of $x$.
Consider the diagram generated by 
exact sequence of pairs $(U,U\setminus V$ and $(Y,Y\setminus V)$.
$$
\begin{CD}
H^r_c(V;G)	@>>i_{V,U}>		H_c^r(U;G)	@>>>			\cdots\\
@VVj_{V,V\cap Y}V			@VVj_{U,Y}V				@.\\
H^r_c(V\cap Y;G)@>>i_{V\cap Y,Y}>	H^r_c(Y;G)	@>>j_{Y,Y\setminus V}>	H^r_c(Y\setminus V;G)\\
\end{CD}
$$
Let $\alpha\in H^r_c(U;G)$ such that $j_{U,Y}(\alpha)\neq 0$.
Since $Y$ is minimal, $j_{Y,Y\setminus V}(j_{U,Y}(\alpha))=0$.
The exactness of the bottom row implies that there is
$\beta\in H^r_c(Y\cap V;G)$ such that $i_{Y\cap V,Y}(\beta)=
j_{U,Y}(\alpha)$. Since $\dim_GV\leq r$, the homomorphism
$j_{V,Y\cap V}$ is an epimorphism and hence there is
$\gamma\in H^r_c(V;G)$ with $j_{V,Y\cap V}(\gamma)=\beta$.
Therefore $j_{U,Y}i_{V,U}(\gamma)\neq 0$ and hence 
$i_{V,U}(\gamma)\neq 0$.
\end{proof}

\begin{definition}
A compactum $X$ is called {\it dimensionally full-valued} if 
$\dim_GX=\dim_{\Z}X$ for all abelian groups $G$. It is clear that every 
$n$-dimensional manifold or $n$-dimensional polyhedron is dimensionally 
full-valued. The following are examples of dimensionally nonfull-valued
compacta.
\end{definition}

\begin{example}[Pontryagin surfaces]
There are 2-dimensional compacta $\Pi_p$ indexed by prime numbers
having the following cohomological dimensions:
$\dim_{\Q}\Pi_p=\dim_{\Z_q}\Pi_p=1$ for prime $q\neq p$ and
$\dim_{\Z_p}\Pi_p=2$.
\end{example}

\begin{proof}
Denote by $M_p$ the mapping cylinder of $p$-to-one covering map of the 
circle to itself $f_p\colon S^1\to S^1$. Denote by $\partial M_p$ the 
domain of the map $f_p$. We construct $\Pi_p$ as the limit space of an 
inverse sequence of polyhedra $\{L_k;q^{k+1}_k\}$ where $L_1$ is a 
2-dimensional sphere and every $L_{k+1}$ is obtained from $L_k$ and a 
triangulation $\tau_k$ on $L_k$ by replacing all 2-simplexes $\Delta$ in 
$L_k$ by $M_p$ identifying the boundary of simplex $\partial\Delta$ with
$\partial M_p$. A bonding map $q^{k+1}_k$ is defined by collapsing the image 
$\Im(f_p)=S^1\subset M_p$ to a point for all $M_p$ participating in the 
construction of $L_{k+1}$. We note that $M_p$ with $\Im(f_p)$ collapsed 
to a point is homeomorphic to a 2-simplex $\Delta$. Denote by $\xi\colon 
M_p\to\Delta$ the corresponding quotient map. In the above construction we 
chose triangulations $\tau_k$ such that preimages
$(q^{\infty}_k)^{-1}(\Delta)$ of 2-dimensional simplexes form a basis of 
topology on $\Pi_p$.

We note that
\begin{enumerate}
\item $H^2(M_p,\partial M_p;\Q)=H^2(M_p,\partial M_p;\Z_q)=0$,
\item $\xi^*\colon H^2(\Delta,\partial\Delta;\Z_p)\to H^2(M_p,\partial 
M_p;\Z_p)$ is an
isomorphism.
\end{enumerate}
To observe (1), (2) we suggest to use the simplicial homology with 
coefficients $\Q$, $\Z_q$ and $\Z_p$. The cohomological results follow 
from the Universal Coefficient Theorem.

By the property (1), $H^2_c\left((q^{k+1}_k)^{-1}(\Int\Delta);F\right)=0$ 
for any 2-simplex $\Delta$ in $L_k$ and for $F=\Q,\Z_q$, $q\neq p$. By 
the Mayer-Vietoris sequence we can get the equality 
$H^2_c\left((q^{k+1}_k)^{-1}(\Int A);F\right)=0$
for any subcomplex $A$ in $L_k$ for the same coefficients. 
Therefore 
\begin{align*}
H^2_c\left((q^{\infty}_k)^{-1}(\Int A);F\right)&
= \Limto H^2_c\left((q^{i+1}_i)^{-1}(q^i_k)^{-1}(\Int A);F\right)\\
&
= \Limto H^2_c\left((q^{i+1}_i)^{-1}(\Int(q^i_k(A)));F\right)\\
&
=0
\end{align*}
for any subcomplex $A\subset L_k$. Since every open set $U\subset\Pi_p$ 
can be presented as an increasing union of sets of the type 
$(q^{\infty}_k)^{-1}(\Int A)$, the formula 
$H^*_c(\limto U_j;F)=\Limto H^*_c(U_j;F)$ implies that 
$H^2_c(U;F)=0$ for every open set $U$ and $F=\Q,\Z_q$. Hence, 
$\dim_{\Q}\Pi_p\leq 1$ and $\dim_{\Z_q}\Pi_p\leq 1$.
The equality holds since $\Pi_p$ is not 0-dimensional.

Similarly the Mayer-Vietoris sequence implies that
$$(q^{k+1}_k)^*\colon H^2(L_k;\Z_p)\to H^2(L_{k+1};\Z_p)$$
is an isomorphism for all $k$.
Hence, $\check H^2(\Pi_p;\Z_p)\neq 0$ and, hence, $\dim_{\Z_p}\Pi_p=2$.
\end{proof}

According to the following theorem a Pontryagin compactum $\Pi_p$ cannot be
imbedded in $\R^3$.

\begin{theorem}
Every $n{-}1$-dimensional compact subset $X$ of the Euclidean space $\R^n$ 
is dimensionally full-valued.
\end{theorem}

\begin{proof}
By Alexandroff Theorem $\dim_{\Z}X=n-1$. According to Theorem 1.8
there is a point $x\in X$ having $n{-}1$-dimensional obstruction
in $X$ with respect to $\Z$. Consider a small ball $U$ in $\R^n$
centered at $x$. Then $H_c^{n-1}(X\cap U;\Z)\neq 0$. By the Alexander
duality $H_0(U\setminus X;\Z)\neq 0$. Since the singular 
0-dimensional homology is always a free group, it follows that the
group $H_c^{n-1}(X\cap U;\Z)$ is free abelian and nontrivial.
The Universal coefficient formula completes the proof.
\end{proof}

A family of subsets $\sU$ of a given set $X$ we call {\it multiplicative}
If $U,V\in\sU$ implies $U\cap V\in\sU$.

\begin{proposition}
Suppose that a compactum $X$ has a multiplicative basis $\sU$ having the 
property $H^k_c(U;G)=0$ for all $k>n$ and for all $U\in\sU$. Then 
$\dim_GX\leq n$.
\end{proposition}

\begin{proof}
Consider a family of open sets 
$\sV = \{V\subset X \mid H^k_c(V;G)=0\ \text{for all}\ k>n\}$.
The Mayer-Vietoris exact sequence
$$
\cdots
\to H^k_c(U;G) \oplus H^k_c(V;G)
\to H^k_c(U \cup V;G)
\to H^{k+1}_c(U\cap V)
\to \cdots
$$
implies that $U \cup V \in \sV$ provided $U,V \in \sV$. Since $\sV$ 
contains a basis $\sU$, it follows that every open set in $X$ is an 
increasing union of sets from $\sV$. The continuity of the cohomology 
implies that every open set in $X$ lies in $\sV$.
\end{proof}

\begin{proposition}
If $\dim_GX<\infty$, then the multiplicativity of the basis $\sU$
in Proposition 1.11 can be omitted.
\end{proposition}

\begin{proof}
If $\dim_GX=r>n$, then according to Theorem 1.8 there is an r-dimensional 
obstruction at some point $x$, which contradicts with the property of the 
basis $\sU$.
\end{proof}

According to Theorem 1.1 for a compactum $X$ to be cohomologically at most 
$n$-dimensional with respect to a coefficient group $G$ it suffices to 
have the property that $H^k_c(U;G)=0$ not for all $k>n$ but just for 
$k=n+1$ and for all open sets $U\subset X$. If instead of all open sets we 
consider only a basis $\sU$, then that property is insufficient even if 
$\sU$ is multiplicative. For example, the unit cube $I^n$ has a 
multiplicative basis $\sU$ consisting of open `rectangles'
$U=I_1\times\dots\times I_n\subset I^n$ of diameter less than one. Since 
every $I_j$ is homeomorphic to an open interval or a half interval, every 
$U$ is homeomorphic to Euclidean space $\R^n$ or half space $\R^n_+$. In 
both cases $H^1_c(U;G)=0$. Thus, $H^1_c(U;G)=0$ for all $U\in\sU$ but 
$I^n$ is far from being 0-dimensional.

\section{Bockstein theory}

As we have seen in \S1 the cohomological dimension of a given compactum 
depends on coefficient group. Any abelian group can be the coefficient 
group of a cohomology theory and there are uncountably many of them. 
It turns out to be that in the case of compacta it suffices to consider 
only countably many groups. 
Solving Alexandroff's problem \cite{A2}, M.F. Bockstein found a countable 
family of abelian groups $\sigma$ and an algorithm for computation of the 
cohomological dimension with respect to a given abelian group by means of 
cohomological dimensions with coefficients taken from $\sigma$. 
The Bockstein basis $\sigma$ consists of the following groups: 
rationals $\Q$, $p$-cyclic groups $\Z_p=\Z/p\Z$, $p$-adic circles 
$\Z_{p^{\infty}}=\Q_p/\A_p$, $p$-adic field factored out by $p$-adic 
integers, and $p$-localizations of integers 
$\Z_{(p)}=\{\frac{m}{n}\in\Q \mid \text{$n$ is not divisible by $p$}\}$ 
where $p$ runs over all primes. 
The set of all $p$-related groups in $\sigma$ we denote by 
$\sigma_p = \{\Z_p,\Z_{p^{\infty}},\Z_{(p)}\}$. 
Thus, 
$\sigma = \bigcup_p\sigma_p \cup \Q$. 
We note that the $p$-adic circle $\Z_{p^{\infty}}$ is the direct limit of 
groups $\Z_{p^k}$.

\begin{definition}
Given an abelian group $G\neq 0$ its Bockstein family 
$\sigma(G)\subset\sigma$ is defined by the following rule:
\begin{enumerate}
\item $\Z_{(p)}\in\sigma(G)$ if and only if $G/\Tor G$ is not divisible 
by $p$,
\item $\Z_p\in\sigma(G)$ if and only if $p-\Tor G$ is not divisible by 
$p$,
\item $\Z_{p^{\infty}}\in\sigma(G)$ if and only if $p-\Tor G\neq 0$ is 
divisible by $p$,
\item $\Q\in\sigma(G)$ if and only if $G/\Tor G\neq 0$ is divisible by 
all $p$.
\end{enumerate}
\end{definition}

\begin{example*}
\mbox{}
\begin{enumerate}
\item $\sigma(\Z)=\{\Z_{(p)} \mid \text{$p$ is prime}\}$,
\item If $G\in\sigma$, then $\sigma(G)=\{G\}$,
\item $\sigma(G) = \sigma(\Tor G) \cup \sigma(G/\Tor G)$ for any abelian 
group $G$.
\end{enumerate}
\end{example*}

\begin{theorem}[Bockstein Theorem]
For any compactum $X$ and for any abelian group $G$,
$\dim_GX=\sup\{\dim_H X \mid H\in\sigma(G)\}$.
\end{theorem}

\begin{lemma}
For any short exact sequence of abelian groups
$0\to G\to E\to\Pi\to 0$ and for any compactum $X$ the following
inequalities hold:
\begin{itemize}
\item[(a)]
$\dim_EX\leq\max\{\dim_GX,\dim_{\Pi}X\}$,
\item[(b)]
$\dim_GX\leq\max\{\dim_EX,\dim_{\Pi}X+1\}$,
\item[(c)]
$\dim_{\Pi}X\leq\max\{\dim_EX,\dim_GX-1\}$.
\end{itemize}
\end{lemma}

\begin{proof}
(a).
Let $n=\max\{\dim_GX,\dim_{\Pi}X\}$. The epimorphism $E\to\Pi$
defines a map $K(E,n)\to K(\Pi,n)$. Turn this map into a Serre
fibration $p$, then the exact sequence of fibration implies that the
homotopy fiber of $p$ is $K(G,n)$. By Theorem 1.1 we have
$K(G,n)\in AE(X)$ and $K(\Pi,n)\in AE(X)$. Then the extension 
theory implies (\cite{Dr4}) that $K(E,n)\in AE(X)$, i.e.\ $\dim_EX\leq n$.

(b).
Let $m=\max\{\dim_EX,\dim_{\Pi}X+1\}$. Here we realize the
monomorphism $G\to E$ by fibration $K(G,m)\to K(E,m)$. The 
homotopy fiber of that is $K(\Pi,m-1)$. Then the result follows.

(c).
The fibration $p\colon K(E,n)\to K(\Pi,n)$ of a. for
$n=\max\{\dim_EX+1,\dim_GX\}$ as any other fibration defines a map
$f\colon \Omega K(\Pi,n)=K(\Pi,n-1)\to p^{-1}(x_0)=K(G,n)$. The Serre 
construction
turns $f$ into a fibration with a fiber $K(E,n-1)$. 
Note that $K(G,n)\in AE(X)$ and $K(E,n-1)\in AE(X)$. Then the 
extension theory implies that $K(\Pi,n-1)\in AE(X)$, i.e.\
$\dim_{\Pi}X\leq n-1=\max\{\dim_EX,\dim_GX-1\}$.
\end{proof}

\begin{proposition}
Every compactum $X$ satisfies the equality 
$\dim_{\Z_p}X=\dim_{\Z_{p^k}}X$ for any $k$ and any prime $p$.
\end{proposition}

\begin{proof}
Induction on $k$. Lemma 2.2 (a) applied to the sequence
$0\to\Z_p\to\Z_{p^{k+1}}\to\Z_{p^k}\to 0$ 
and the induction
assumption establish the inequality 
$\dim_{\Z_{p^{k+1}}}X\leq \dim_{\Z_p}X$. Lemma 2.2 (c) together
with the induction assumption give an opposite inequality.
\end{proof}

\enlargethispage{\baselineskip}
\begin{theorem}[Bockstein Inequalities]
For any compactum $X$ the following inequalities hold:
\begin{itemize}
\item[BI1]
$\dim_{\Z_{p^{\infty}}}X\leq \dim_{\Z_p}X$,
\item[BI2]
$\dim_{\Z_p}X\leq \dim_{\Z_{p^{\infty}}}X+1$,
\item[BI3]
$\dim_{\Z_p}X\leq \dim_{\Z_{(p)}}X$,
\item[BI4]
$\dim_{\Q}X\leq \dim_{\Z_{(p)}}X$,
\item[BI5]
$\dim_{\Z_{(p)}}X\leq\max\{\dim_{\Q}X,\dim_{\Z_{p^{\infty}}}X+1\}$,
\item[BI6]
$\dim_{\Z_{p^{\infty}}}X\leq\max\{\dim_{\Q}X,\dim_{\Z_{(p)}}X-1\}$.
\end{itemize}
\end{theorem}

\begin{proof}
Since the $p$-adic circle can be presented as the direct
limit of groups $\Z_{p^k}$, Lemma 2.2 and Theorem 1.6 imply BI1.

Lemma 2.2 (b) applied to the sequence
$0\to\Z_p\to\Z_{p^{\infty}}\to\Z_{p^{\infty}}\to 0$ implies BI2.

Lemma 2.2 (c) applied to the sequence
$0\to\Z_{(p)}\to\Z_{(p)}\to\Z_p\to 0$ implies BI3.

Lemma 2.2 (a) applied to the sequence
$0\to\Z_{(p)}\to\Q\to \Z_{p^{\infty}}\to 0$, BI1 and BI3 imply BI4.

Lemma 2.2 (b) applied to the above sequence gives BI5.

Lemma 2.2 (c) applied to the same sequence gives BI6.
\end{proof}

\begin{lemma}
Let $G$ be an abelian group, then
$$\dim_GX=\max\{\dim_{\Tor G}X,\dim_{G/\Tor G}X\}$$ for every compactum $X$.
\end{lemma}

\begin{proof}
Since $H^{k+1}(K(\Tor G,k);\Q)=0$, it follows that the Bockstein long 
exact sequence generated by $0\to \Tor G\to G\to G/\Tor G\to 0$ is split 
into short exact sequences
$$
0\to \check H^k(Y;\Tor G)\to \check H^k(Y;G)\to \check H^k(Y;G/\Tor G)\to 0.
$$
Then the result follows.
\end{proof}

\begin{proof}[Proof of Bockstein Theorem]
First we consider the case when $G$ is a torsion group.
Then $G = \Tor G = \bigoplus_p p-\Tor G$. By 1.7 it follows that 
$\dim_GX=\sup\{\dim_{p-\Tor G}X\}$. 
Since $\sigma(\Tor G) = \bigcup\sigma (p-\Tor G)$, it suffices to show 
that 
$$\dim_{p-\Tor G}X=\sup\{\dim_H X \mid H\in\sigma(p-\Tor)\}.$$
Indeed, then
\begin{align*}
\dim_G X&
= \sup_p\{\dim_{p-\Tor G}X\}\\
&
= \sup_p\sup\{\dim_H X \mid H\in\sigma(p-\Tor G)\}\\
&
= \sup\{\dim_H X \mid H\in\bigcup_p\sigma(p-\Tor)\}\\
&
= \sup\{\dim_H X \mid H\in\sigma(G)\}.
\end{align*} 
If the group $p-\Tor G$ is not divisible by $p$, then it contains 
$\Z_{p^k}$ as a direct summand of $G$ for some $k\geq 1$. In that case 
$\sigma(p-\Tor G) = \{\Z_p\}$. By 1.7 we have 
$$\dim_{p-\Tor G}X
\geq \dim_{\Z_{p^k}}X
= \dim_{\Z_p}X
= \sup\{\dim_H X \mid H\in\sigma(p-\Tor G)\}.
$$ 
Here we applied Proposition 2.3 to obtain the second equality.
On the other hand, $p-\Tor G$ is a direct limit of finite abelian 
$p$-groups which are direct sums of groups isomorphic to $\Z_{p^m}$ for 
some $m$. Thus, by Theorem 1.6, 1.7 and Proposition 2.3, 
$$\dim_{p-\Tor G}X
\leq \dim_{\Z_p}X
= \sup\{\dim_H X \mid H\in\sigma(p-\Tor G)\}.
$$

Now we consider the case when $G$ is a torsion free group.
By the Universal Coefficient Formula $\check H^{n+1}(X,A;G)\neq 0$ if and 
only if $\check H^{n+1}(X,A)\otimes G\neq 0$ which is equivalent to 
$\check H^{n+1}(X,A)\otimes \Z_{(p)}\neq 0$ for all $p$ such that 
$\Z_{(p)}\in\sigma(G)$. 
By the Universal Coefficient Formula the latter is equivalent to $\check 
H^{n+1}(X,A;\Z_{(p)})\neq 0$ for all $p$ such that
$\Z_{(p)}\in\sigma(G)$.
Now the result follows from Theorem 1.1.

If $G$ is an arbitrary abelian group, then by Lemma 2.5, 
\begin{align*}
\dim_GX&
= \max\{\dim_{\Tor G}X, \dim_{G/\Tor G}X\}\\
&
= \sup\{\dim_H X \mid H\in\sigma(\Tor G)\cup\sigma(G/\Tor G)\}\\
&
= \sup\{\dim_H X \mid H\in\sigma(G)\}.
\end{align*}
\end{proof}

\begin{definition}
A compactum $X$ is $p$-{\it regular} if all its $p$-dimensions agree and 
coincide with the rational dimension: 
$$\dim_{\Z_p}X 
= \dim_{\Z_{p^{\infty}}}X 
= \dim_{\Z_{(p)}}X 
= \dim_{\Q}X.$$
Otherwise we call a compactum $X$ $p$-{\it singular}.
\end{definition}

\begin{lemma}
A compact $X$ is $p$-regular if and only if 
$\dim_{\Z_{(p)}}X = \dim_{\Z_{p^{\infty}}}X$.
\end{lemma}

\begin{proof}
Bockstein inequalities BI1 and BI3 imply that 
$\dim_{\Z_p}X = \dim_{\Z_{p^{\infty}}}X = \dim_{\Z_{(p)}}X$,
The inequalities BI4 and BI6 imply that
$\dim_{\Z_{(p)}}X=\dim_{\Q}X$.
\end{proof}

The following theorem we call the Bockstein Alternative (BA).

\begin{theorem}[Bockstein Alternative (BA)]
For any compactum $X$ there is an alternative: either 
\begin{align*}
\dim_{\Z_{(p)}}X&
=\dim_{\Q}X \quad\text{or}\\
\dim_{\Z_{(p)}}X&
=\dim_{\Z_{p^{\infty}}}X+1.
\end{align*}
\end{theorem}

\begin{proof}
It is clear that BA holds when $X$ is $p$-regular. Consider $p$-singular 
$X$ and assume that $\dim_{\Z_{(p)}}X\neq \dim_{\Q}X$. Then by BI4, 
$\dim_{\Z_{(p)}}X\geq \dim_{\Q}X$. Then BI5 implies that 
$\dim_{\Z_{(p)}}X\leq \dim_{\Z_{p^{\infty}}}X+1$. Since $X$ is
$p$-singular and in the view of BI1, BI3, Lemma 2.6 implies that 
$\dim_{\Z_{(p)}}X=\dim_{\Z_{p^{\infty}}}X+1$.
\end{proof}

\begin{remark*}
In the case of $p$-singular $X$, 
$\dim_{\Z_{(p)}}X=\max\{\dim_{\Q}X,\dim_{\Z_{p^{\infty}}}X+1\}$.
\end{remark*} 

\begin{definition}
$p$-{\it deficiency} $\epsilon_X(p)$ of a compactum $X$ is the difference
$\dim_{\Z_p}X - \dim_{\Z_{p^{\infty}}}X$. The inequalities BI1, BI2 imply 
that $\epsilon_X(p)\in\{0,1\}$.
\end{definition}

Let $\sP$ be the set of all prime numbers. For every compactum $X$ by 
$\sS_X\subset\sP$ we denote the set of $p$ for which $X$ is $p$-singular 
and by $\sD_X\subset\sP$ the set of all $P$ for which $X$ is 
$p$-deficient. It is clear that $\sD_X\subset\sS_X$. Then the deficiency 
function $\epsilon_X(\ )$ is just the characteristic function of the set 
$\sD_X$. Additionally we introduce {\it the field dimensional function} 
$d_X\colon \sP\cup\{0\}\to\N\cup\{\infty\}$ by the formulas: 
$d_X(p)=\dim_{\Z_p}X$ and $d_X(0)=\dim_{\Q}X$.

\begin{lemma}
The family $(\sS_X,\sD_X;d_X)$ consisting of the pair of the singularity
set and the deficiency set $\sD_X\subset\sS_X\subset\sP$ together with 
the field dimensional function $d_X$ completely determine cohomological 
dimensions of a given compactum $X$. Moreover for the groups from the 
basis $\sigma$ there are formulas:
\begin{enumerate}
\item $\dim_{\Q} X = d_X(0)$,
\item $\dim_{\Z_p} X = d_X(p)$,
\item $\dim_{\Z_{p^{\infty}}}X = d_X(p) - \chi_{\sD_X}(p)$ and 
\item $\dim_{\Z_{(p)}}X 
= \left(\max\{d_X(0), d_X(p) - \chi_{\sD_X}(p) + 1\}\right) \chi_{\sS}(p)
+ d_X(0)\chi_{\sP\setminus\sS}(p)$ 
\end{enumerate}
where $\chi_A$ denotes the characteristic function of a set $A$.
\end{lemma}

\begin{proof}
In view of Bockstein Theorem it is sufficient to prove the formulas.
The first formula is obvious. If $p\in\sP\setminus\sS$, then $X$ is 
$p$-regular and the formula (2) holds. If $p\in\sS$, then by BI5
$$\dim_{\Z_{(p)}}X\leq\max\{\dim_{\Q}X,\dim_{\Z_{p^{\infty}}}X+1\}$$ 
and
$$\dim_{\Z_{(p)}}X\geq\max\{\dim_{\Q}X,\dim_{\Z_{p^{\infty}}}X+1\}$$ 
by BI4, BI1, BI3 and Lemma 2.6.
\end{proof}

\begin{lemma}
For every compactum $X$ there is an additive group of a field 
$F\in\sigma$ such that $\dim_{\Z}X\leq \dim_FX+1$.
\end{lemma}

\begin{proof}
By the Bockstein Theorem (2.1), $\dim_{\Z} X =\dim_{\Z_{(p)}} X$ for some 
$p$. By the Bockstein Alternative (Theorem 2.7), either 
$\dim_{\Z_{(p)}} X=\dim_{\Q}X$ or 
$\dim_{\Z{(p)}} X=\dim_{\Z_{p^{\infty}}}X+1$. In the first case we take 
$F=\Q$, in the second case $F=\Z_p$. The inequality BI1 completes the 
proof in the second case.
\end{proof}

\begin{example}
A Pontryagin surface $\Pi_p$ has the following cohomological dimensions 
with respect to Bockstein groups $G\in\sigma$:
$\dim_{\Z_{p^{\infty}}} \Pi_p
= \dim_{\Q}\Pi_p
= \dim_{\Z_q} \Pi_p
= \dim_{\Z_{q^{\infty}}} \Pi_p
= \dim_{\Z_{(q)}} \Pi_p
= 1$ for $q\neq p$ and
$\dim_{\Z_p} \Pi_p = \dim_{\Z_{(p)}} \Pi_p = 2$.
\end{example}

\begin{proof}
First we note that a compactum $\Pi_p$ is $q$-regular for $q\neq p$. 
Since it is 2-dimensional, by the Bockstein theorem 
$\dim_{\Z_{(p)}} \Pi_p=2$. 
By BA we have $\dim_{\Z_{p^{\infty}}} \Pi_p = 1$.
\end{proof}

\section{Cohomological dimension of Cartesian product}

Theorem 1.5 allows to compute easily the cohomological dimension of the 
union of two compacta: $\dim_GX\cup Y=\max\{ \dim_GX, \dim_GY\}$. 
Unfortunately there is no easy way to compute the cohomological dimension 
of the product of two compacta. The natural formula $\dim_G(X\times 
Y)=\dim_GX+\dim_GY$ can be violated in both directions.

\begin{proposition}
Let $X$ and $Y$ be compacta and $G$ an abelian group. Then the following 
conditions are equivalent:
\begin{enumerate}
\item $\dim_G(X\times Y)\leq n$,
\item $H^k_c(U\times V;G)=0$ for all $k>n$ and all open subsets $U$ of 
$X$ and $V$ of $Y$,
\item $\check H^k((X/A)\wedge (Y/B);G)=0$ for all $k>n$ and all closed 
subsets $A$ of $X$ and $B$ of $Y$.
\end{enumerate}
\end{proposition}

\begin{proof}
(1) $\Rightarrow$ (2). 
It follows from Theorem 1.1.

(2) $\Rightarrow$ (1). 
Note that the family 
$\sU = \{U\times V \mid \text{$U$ is open in $X$}, \text{$V$ is open in $Y$}\}$ 
forms a multiplicative basis in $X\times Y$. 
Now Proposition 1.11 implies the proof.

(2) $\Leftrightarrow$ (3). 
Denote $A=X\setminus U$ and $B=Y\setminus V$. 
Then
\begin{align*}
H^k_c(U\times V;G)&
=\check H^k(X\times Y,X\times Y\setminus U\times V;G)\\
&
=\check H^k(X\times Y,X\times B\cup A\times Y;G)\\
&
=\check H^k(X\times Y/X\times B\cup A\times Y;G)\\
&
=\check H^k((X/A)\wedge (Y/B);G)
\end{align*} 
and the result follows.
\end{proof}

\begin{proposition}
Let $X$ and $Y$ be compacta and let $G\neq 0$ be an abelian group.
\begin{enumerate}
\item If $k\geq dim_GY$ is a number such that $\dim_{H^{k-i}_c(V;G)}X\leq i$ 
for all $i\geq 0$ and
all open subsets $V$ of $Y$, then $\dim_G(X\times Y)\leq k$,
\item If $\dim_{H^n_c(V;G)}X\geq m$, then $\dim_G(X\times Y)\geq n+m$.
\end{enumerate}
\end{proposition}

\begin{proof}
(1). 
Since $\dim_{H^{k-i}_c(V;G)}X\leq i$, we have
$H^{i+l}_c(U;H^{k-i}_c(V;G))=0$ for any $l>0$ and any open subset 
$U\subset X$ for all $i\geq 0$. 
By the Kunneth formula we have
\begin{align*}
H^{k+l}_c(U\times V;G)&
= \bigoplus_{j=0}^{k+l}H^j_c(U;H^{k+l-j}_c(V;G))\\
&
= \bigoplus_{j=0}^{l-1}H^j_c(U;H^{k+l-j}_c(V;G))
\oplus \bigoplus_{i=0}^kH^{i+l}_c(U;H^{k-i}_c(V;G))\\
&
= 0
\end{align*}
The first sum is zero by the assumption $k\geq \dim_GY$ and the second 
part is zero by the above formula. Proposition 1.11 completes the proof. 

(2). 
Since $\dim_{H^n_c(V;G)}X\geq m$, by virtue of Theorem 1.1 there exists 
an open subset $U\subset X$ such that $H^m_c(U;H^n_c(V;G))\neq 0$. By the 
Kunneth formula we have $H^{n+m}_c(U\times V;G)\neq 0$. 
Hence, $\dim_G(X\times Y)\geq n+m$.
\end{proof}

\begin{proposition}
For an additive group of a field $F$ the formula 
$\dim_F(X\times Y)=\dim_FX+\dim_FY$
holds for all compacta.
\end{proposition}

\begin{proof}
Let $\dim_FX=m$ and $\dim_FY=n$. Note that $\dim_{H^{n+m-i}_c(V;F)}X=0$ 
if $i<m$ and $\dim_{H^{n+m-i}_c(V;F)}X=\dim_{\bigoplus F}X\leq \dim_FX =m$ 
if $i\geq m$. In both cases $\dim_{H^{n+m-i}_c(V;F)}X\leq i$ and by 
Proposition 3.2 (1) it follows $\dim_F(X\times Y)\leq n+m$. Let $V$ be an 
open subset of $Y$ with $H^n_c(V;F)\neq 0$. 
Then $H^n_c(V;F)=\bigoplus F\neq 0$. 
Then $\dim_{H^n_c(V;F)}X\geq m$ and by Proposition 3.2 (2) we have 
$\dim_F(X\times Y)\geq n+m$. 
Therefore $\dim_F(X\times Y)=\dim_FX+\dim_FY$.
\end{proof}

\begin{proposition}
Suppose $X$ and $Y$ are compacta and $G$ is an abelian group.
\begin{enumerate}
\item $\dim_G(X\times Y)\leq \dim_GX+\dim_GY$ if $G$ is torsion free,
\item $\dim_G(X\times Y)\leq \dim_GX+\dim_GY+1$ in general case.
\end{enumerate}
\end{proposition}

\begin{proof}
(1). 
Let $\dim_GX=m$ and $\dim_GY=n$. Since $G$ is torsion free, 
$H^l_c(V;G)=H^l_c(V)\otimes G$ by virtue of the Universal Coefficient 
Formula. We note that if $\Z_{(p)}\in\sigma(H\otimes G)$, then $H\otimes 
G$ is not divisible by $p$ and, hence, $G$ is not divisible by $p$, 
therefore, $\Z_{(p)}\in\sigma(G)$. Then 
\begin{align*}
\dim_{H^l_c(V;G)}X&
= \dim_{H^l_c(V)\otimes G}X\\
&
= \sup\{\dim_{\Z_{(p)}}X \mid \Z_{(p)}\in\sigma(H^l_c(V)\otimes G)\}\\
&
\leq \sup\{\dim_{\Z_{(p)}}X \mid \Z_{(p)}\in\sigma(G)\}\\
&
= \dim_GX\\
&
=m.
\end{align*} 
Therefore $\dim_{H^{n+m-i}_c(V;G)}X\leq i$ for all $i\geq 0$. 
Then by Proposition 2.2 (1), $\dim_G(X\times Y)\leq n+m$.

(2). 
First, we prove the inequality for the $p$-adic circle 
$G=\Z_{p^{\infty}}$. We note that 
$\sigma(H\otimes\Z_{p^{\infty}})=\{\Z_{p^{\infty}}\}$ or $\emptyset$ and 
$\sigma(H\ast \Z_{p^{\infty}})\subset\{\Z_p,\Z_{p^{\infty}}\}$. 
Then it follows that
$\dim_{H^l_c(V;\Z_{p^{\infty}})}X
\leq \dim_{\Z_p}X
\leq \dim_{\Z_{p^{\infty}}}X+1
= m+1$.
Therefore $\dim_{H^{n+m-i}_c(V;\Z_{p^{\infty}})}X\leq i+1$ for all 
$i\geq 0$ and hence, $\dim_{\Z_{p^{\infty}}}(X\times Y)\leq n+m+1$. 
Now we have proven the inequality (2) for all groups from Bockstein basis 
(additionally to an above see also Proposition 3.3, Proposition 3.4 (1)). 
If $G$ is an arbitrary group, then by the Bockstein Theorem, 
\begin{align*}
\dim_G(X\times Y)&
= \sup\{\dim_H(X\times Y) \mid H\in\sigma(G)\}\\
&
\leq \sup\{\dim_HX+\dim_H Y+1 \mid H\in\sigma(G)\}\\
&
\leq \sup\{\dim_H X \mid H\in\sigma(G)\} 
+ \sup\{\dim_H Y \mid H\in\sigma(G)\}
+ 1\\
&
= \dim_GX + \dim_GY + 1.
\end{align*}
\end{proof}

\begin{proposition}
Suppose that compactum $X$ is not $p$-deficient, then 
$\dim_{\Z_{p^{\infty}}}(X\times Y)
= \dim_{\Z_{p^{\infty}}}X+\dim_{\Z_{p^{\infty}}}Y$ 
for every compactum $Y$.
\end{proposition}

\begin{proof}
Denote $\dim_{\Z_p}X=m$ and $\dim_{\Z_p}Y=n$.
Since $X$ is not $p$-deficient, we have $\dim_{\Z_{p^{\infty}}}X=m$. 
According to Bockstein inequalities BI1, BI2 there are two possibilities 
for $\dim_{\Z_{p^{\infty}}}Y$:
(a) to be equal $n$ and 
(b) to be equal $n-1$. 
In the case of (a) we can find an open set $V\subset Y$ with 
$H^n_c(V;\Z_{p^{\infty}})\neq 0$. 
Since $H^n_c(V;\Z_{p^{\infty}})$ is $p$-torsion group, 
$\sigma(H^n_c(V;\Z_{p^{\infty}}))\subset\{\Z_p,\Z_{p^{\infty}}\}$.
By the Bockstein Theorem and BI1, $\dim_{H^n_c(V;\Z_{p^{\infty}})}X=m$. 
Proposition 3.2 (2) implies that 
$\dim_{\Z_{p^{\infty}}}(X\times Y)\geq n+m$. 
The inequality $\dim_{\Z_{p^{\infty}}}(X\times Y)\leq n+m$ follows from 
Proposition 3.3 and BI1.
In the case of (b) one can show that 
$\dim_{H^{m+n-i-1}_c(V;\Z_{p^{\infty}})}X\leq i$
for all $i\geq 0$ and every open subset $V\subset Y$. 
For $i<m$ it is due to an obvious reason: 
$H^{m+n-i-1}_c(V;\Z_{p^{\infty}})=0$. 
For $i\geq m$ the inequality holds because of the inclusion 
$\sigma(H^{m+n-i-1}_c(V;\Z_{p^{\infty}}))
\subset \{\Z_p,\Z_{p^{\infty}}\}$ 
and the equality $\dim_{\Z_{p^{\infty}}}X=m$. 
Then Proposition 3.2 (1) implies that 
$\dim_{\Z_{p^{\infty}}}(X\times Y)\leq n+m-1$. 
The opposite inequality $\dim_{\Z_{p^{\infty}}}(X\times Y)\geq n+m-1$ 
follows from Proposition 3.3 and BI2.
\end{proof}

\begin{theorem}
Suppose that a compactum $X$ is $p$-regular for some prime $p$, then 
$\dim_G (X\times Y) = \dim_G X + \dim_G Y$ 
for all $G\in\sigma_p=\{\Z_p,\Z_{p^{\infty}},\Z_{(p)}\}$ and any other 
compactum $Y$.
\end{theorem}

\begin{proof}
Obviously theorem is true for $G=\Z_p$. 

Since $p$-regularity does not admit $p$-deficiency, the case 
$G=\Z_{p^{\infty}}$ follows from Proposition 3.5.

In the case of $G=\Z_{(p)}$ we denote by $n=\dim_{\Z_{(p)}}Y$ and 
$m=\dim_{\Z_{(p)}}X$.
Let $V$ be an open subset of $Y$ such that $A=H^n_c(V;\Z_{(p)})\neq 0$. 
If $A$ is not a torsion group, then by the Bockstein Theorem 
$\dim_AX\geq \dim_{\Z_{(q)}}X$. 
By BI4 we have $\dim_AX\geq \dim_{\Q}X=m$. 
Proposition 3.2 (2) implies that 
$\dim_{\Z_{(p)}}(X\times Y)\geq \dim_{\Z_{(p)}}X+\dim_{\Z_{(p)}}Y$. 
In the other direction the inequality follows by Proposition 3.4 (1). 
If $A$ is a torsion group, then $A$ is a $p$-torsion group, since
$A = H^n_c(V)\otimes {\Z_{(p)}}$ by the Universal Coefficient Formula. 
Therefore $\dim_AX\geq \dim_{\Z_{p^{\infty}}}X=m$ and the result follows.
\end{proof}

\begin{corollary} 
Suppose $X$ is a dimensionally full-valued compactum. 
Then\linebreak[5]
$\dim_G (X\times Y) = \dim_G X + \dim_G Y$
for any group $G$.
\end{corollary}

\begin{proof}
A compactum $X$ is $p$-regular for all $p$.
Hence $\dim_GX=\dim_{\Z}X$ for any group $G$.
Theorem 3.6 and Proposition 3.3 imply that the above formula holds for all 
$G\in\sigma$. 
If $G$ is an arbitrary abelian group, the Bockstein Theorem states that 
$\dim_G(X\times Y)
= \sup\{\dim_H(X\times Y) \mid H \in\sigma(G)\}
= \sup\{\dim_HX+\dim_HY \mid H\in\sigma(G)\}
= \dim_GX+\sup\{\dim_HY \mid H\in\sigma(G)\}
= \dim_GX+\dim_GY$.
\end{proof}

\begin{corollary}
\mbox{}
\begin{enumerate}
\item The product of two $p$-regular compacta is $p$-regular.
\item The product of $p$-regular and $p$-singular compacta is 
$p$-singular.
\end{enumerate}
\end{corollary}

\begin{example*}
Let $p\neq q$, then $\dim(\Pi_p\times\Pi_q)=3$ for different Pontryagin 
surfaces. Indeed, by theorems of Alexandroff and Bockstein,
$$\dim(\Pi_p\times\Pi_q)
= \max\{\dim_{\Z_{(r)}}(\Pi_p\times\Pi_q) \mid r\in\sP\}.$$
Since for every $r\in\sP$ one of the factors $\Pi_p$ or $\Pi_q$ is 
$r$-regular, by Theorem 3.6, 
$$\dim_{\Z_{(r)}}(\Pi_p\times\Pi_q)
= \dim_{\Z_{(r)}}\Pi_p+\dim_{\Z_{(r)}}\Pi_q.$$
Then $\dim_{\Z_{(r)}}(\Pi_p\times\Pi_q)=3$ if $r=p$ or $r=q$ and it 
equals $2$ if $r\neq p$ and $r\neq q$.
\end{example*}

\begin{lemma}
The deficiency set of the product is the union of deficiency sets of 
factors: $\sD_{X\times Y}=\sD_X\cup\sD_Y$.
\end{lemma}

\begin{proof}
By Propositions 3.3 and 3.5 the product $X\times Y$ cannot be 
$p$-deficient if both factors are not $p$-deficient. 
This implies the inclusion $\sD_{X\times Y}\subset\sD_X\cup\sD_Y$.

If $p\in\sD_X\setminus\sD_Y$, the $p$-deficiency of the product 
$X\times Y$ equals one by Propositions 3.3 and 3.5, and hence 
$p\in\sD_{X\times Y}$. 
Similarly if $p\in\sD_Y\setminus\sD_X$. 
If $p\in\sD_X\cap\sD_Y$, then by Proposition 3.4, 
$\dim_{\Z_{p^{\infty}}}(X\times Y)
\leq \dim_{\Z_{p^{\infty}}}X
+ \dim_{\Z_{p^{\infty}}}Y
+ 1
= \dim_{\Z_p}X
- 1
+ \dim_{\Z_p}Y
- 1
+ 1
= \dim_{\Z_p}(X\times Y)
- 1$. 
Then BI2 implies that
$\dim_{\Z_{p^{\infty}}}(X\times Y) = \dim_{\Z_p}(X\times Y) - 1$. 
It means that $p\in\sD_{X\times Y}$ in that case too. 
Thus, $\sD_{X\times Y}\supset\sD_X\cup\sD_Y$.
\end{proof}

\begin{corollary}
The $p$-deficiency of the product of two compacta can be computed by the 
following formula: 
$\epsilon_{X\times Y}(p)
= \epsilon_X(p) + \epsilon_Y(p) - \epsilon_X(p)\epsilon_Y(p)$.
\end{corollary}

\begin{proof}
The formula follows from the union formula for characteristic functions
$\chi_{A\cup B}
= 1 - (1 - \chi_{A})(1 - \chi_{B})
= \chi_{A} + \chi_{B} - \chi_{A}\chi_{B}$, 
Lemma 3.9 and the equality $\epsilon_X = \chi_{\sD_X}$.
\end{proof}

\begin{lemma}
The inequality 
$\dim_{\Z_{(p)}}(X\times Y)
\geq \dim_{\Z_{p^{\infty}}}X + \dim_{\Z_{p^{\infty}}}Y + 1$
holds for all $p$ and all $p$-singular compacta $X$ and $Y$.
\end{lemma}

\begin{proof}
Let $k = \dim_{\Z_{p^{\infty}}}X<\dim_{\Z_{(p)}}X$ and
$l = \dim_{\Z_{p^{\infty}}}Y < \dim_{\Z{(p)}}Y$. 
Consider a group $G=H^{l+1}_c(V;\Z_{(p)}) = H^{l+1}_c(V)\otimes\Z_{(p)}$ 
for an open subset $V\subset Y$ such that $H^{l+1}_c(V)\neq 0$.
Such a set $V$ exists because of Theorem 1.1 and the inequality 
$\dim_{\Z_{(p)}}Y\geq l+1$.
If the group $G$ has $p$-torsion, then $\Z_p$ or $\Z_{p^{\infty}}$ 
belongs to $\sigma(G)$.
In both cases $\dim_GX\geq \dim_{\Z_{p^{\infty}}}X=k$. 
By Proposition 3.2 (2), $\dim_{\Z_{(p)}}(X\times Y)\geq k+l+1$. 
If the group $G$ has no $p$-torsion, then $H^{l+1}_c(V)\otimes\Q\neq 0$ 
and hence, $\dim_{\Q}Y\geq l+1$.

Similarly, consider a group $G'=H^{k+1}_c(U;\Z_{(p)})$ and derive 
$\dim_{\Q}X\geq k+1$ or the required inequality 
$\dim_{\Z_{(p)}}(X\times Y)\geq k+l+1$. 
In the first case according to BI4 we have
$\dim_{\Z_{(p)}}(X\times Y)
\geq \dim_{\Q}(X\times Y)
\geq k+l+1$.
\end{proof}

\begin{corollary}
The product $X\times Y$ of two $p$-singular compacta is $p$-singular.
\end{corollary}

\begin{proof}
If one of the compacta is $p$-deficient, then by Lemma 3.8 the product is 
also $p$-deficient and, hence, $p$-singular. 
If both compacta are not $p$-deficient, then Lemma 3.11 implies 
$\dim_{\Z_{(p)}}(X\times Y)
\geq \dim_{\Z_{p^{\infty}}}X + \dim_{\Z_{p^{\infty}}}Y + 1
= \dim_{\Z_p}X + \dim_{\Z_p}Y + 1
= \dim_{\Z_p}(X\times Y) + 1$ 
and, hence, $X\times Y$ is $p$-singular.
\end{proof}

\begin{lemma}
The singularity set of the product of two compacta is the union of their 
singularity sets: $\sS_{X\times Y}=\sS_X\cup\sS_Y$.
\end{lemma}

\begin{proof}
Corollaries 3.8 and 3.12 imply the proof.
\end{proof}

The results of Lemmas 3.9, 3.13 and Proposition 3.3 can be summarize into 
the following.

\begin{theorem}
For any two compacta $X$ and $Y$ and their product $X\times Y$ there is 
the formula:
$(\sS_{X\times Y},\sD_{X\times Y}, d_{X\times Y})
= (\sS_X\cup\sS_Y,\sD_X\cup\sD_Y,d_X+d_Y)$.
\end{theorem}

If one of the factors is $p$-regular, then according to Theorem 3.6, the 
logarithmic law for the dimension of the product holds. If both factors 
are $p$-singular then the following deviation takes place for coefficient 
groups from Bockstein basis $\sigma$.

\begin{lemma}
Suppose $X$ and $Y$ are $p$-singular compacta. Then
\begin{enumerate}
\item 
$\dim_{\Z_{p^{\infty}}}(X\times Y)
= \dim_{\Z_{p^{\infty}}}X
+ \dim_{\Z_{p^{\infty}}}Y
+ \epsilon_X(p)\epsilon_Y(p)$
\item 
$\dim_{\Z_{(p)}}(X\times Y)
= \max\{\dim_{\Z_{p^{\infty}}}(X\times Y)+1, \dim_{\Q}(X\times Y)\}$.
\end{enumerate}
\end{lemma}

\begin{proof}
Proposition 3.3 and Corollary 3.10 imply that 
\begin{align*}
\dim_{\Z_{p^{\infty}}}(X\times Y)&
= \dim_{\Z_p}(X\times Y)-\epsilon_{X\times Y}(p)\\
&
= \dim_{\Z_p}X 
+ \dim_{\Z_p}Y 
- \epsilon_X(p) 
- \epsilon_Y(p) 
+ \epsilon_X(p)\epsilon_Y(p)\\
&
= \dim_{\Z_{p^{\infty}}}X 
+ \dim_{\Z_{p^{\infty}}}Y 
+ \epsilon_X(p)\epsilon_Y(p).
\end{align*}

Corollary 3.12 and Lemma 2.8(2) imply the second part of the theorem.
\end{proof}

\begin{theorem}
Let $X$ be a compactum, then
\begin{itemize}
\item[(a)]
$\dim_{\Z}(X\times X)=2\dim_{\Z}X$ or\/ $2\dim_{\Z}X-1$,
\item[(b)]
$\dim_{\Z}X^n=n\dim_{\Z}X$ or\/ $n\dim_{\Z}X-n+1$.
\end{itemize}
\end{theorem}

\begin{proof}
If there is a field $F$ such that $\dim_FX=\dim_{\Z}X$ then by 
Propositions 3.3 and 3.4 we have the first case. Now assume that there is 
no such a field. Then by Bockstein Theorem $\dim_{\Z}X=\dim_{\Z_{(p)}}X$ 
for some $p$. Our assumption implies that $X$ is $p$-singular and 
$\dim_{\Z_{(p)}}X>\dim_{\Q}X$. Lemma 3.15(1) states that 
$\dim_{\Z_{p^{\infty}}}(X\times X)
= \dim_{\Z_{p^{\infty}}}X+\epsilon^2_X(p)$.
By Lemma 3.15(2), we have
$\dim_{\Z_{(p)}}(X\times X)
= 2\dim_{\Z_{p^{\infty}}}X+\epsilon_X(p)+1$. 
Bockstein inequality BI1 and the assumption imply that 
$\dim_{\Z_{p^{\infty}}}X=\dim_{\Z_p}X$ and hence $\epsilon_X(p)=0$. By 
Lemma 2.9 there is a field $F'$ such that $\dim_{\Z}(X\times X)\leq 
\dim_{F'}(X\times X)+1$. Since $\dim_{F'}X\leq \dim_{\Z}X-1$, we have
\begin{align*}
2\dim_{\Z}X-1&
= \dim_{\Z_{(p)}}(X\times X)\\
&
\leq \dim_{\Z}(X\times X)\\
&
\leq \dim_{F'}(X\times X) + 1\\
&
= 2\dim_{F'}X + 1\\
&
\leq 2(\dim_{\Z}X-1) + 1\\
&
= 2\dim_{\Z}X - 1.
\end{align*}
Hence, $\dim_{\Z}(X\times X) = 2\dim_{\Z}X-1$.

Induction on $n$ implies part (b).
\end{proof}

\begin{definition}
A compactum $X$ is of the {\it basic type} if $\dim X^2=2\dim X$ and 
it is  called having the {\it exceptional type} if $\dim X^2=2\dim X-1$.
\end{definition}
 
This definition makes sense only for finite dimensional compacta. In that 
case $\dim X=\dim_{\Z}X$ by virtue of Alexandroff Theorem. Theorem 3.16 
proves that all compacta are split into these two classes.
Moreover, the dimension of the $n$-th power of $X$ equals 
$\dim X^n=n\dim X$ for compacta of the basic type and $\dim X^n=n\dim 
X-n+1$ for compacta of the exceptional type.

The proof of Theorem 3.16 suggests the following:

\begin{criterion}
A compactum $X$ is of the basic type if and only if there is a field 
$F\in\sigma$ such that $\dim_FX=\dim X$.
\end{criterion}

\section{Dimension type algebra}

Every compactum $X$ of positive dimension defines a function 
$\phi_X\colon \sigma\to\N\cup\{\infty\}$ by the formula 
$\phi_X(G)=\dim_GX$. 
This function $\phi_X$ satisfies the Bockstein Inequalities BI1--6. 
In Lemma 2.8 we defined a set $F=(\sS_X,\sD_X;d_X)$ where
$\sD_X\subset\sS_X\subset\sP$ is a pair of subsets of primes and
$d_X\colon \sP\cup\{0\}\to\N\cup\{\infty\}$
is the field dimensional function. 
The function $d_X$ has the property $d(\sP\setminus\sS) = d(0)$. 
The set $(\sS_X,\sD_X;d_X)$ completely defines $\phi_X$. 
Now if we forget that the function $\phi_X$ came from a compactum $X$, we 
can reformulate the results of \S2 in more abstract way. 
For every abstract function $\phi\colon \sigma\to\N\cup\{\infty\}$ one 
can define a regularity set 
$$\sR
= \{p\in\sP \mid \phi(\Z_{p^{\infty}})
= \phi(\Z_p)
= \phi(\Z_{(p)})
= \phi(\Q)\},$$
a singularity set
$\sS = \sP\setminus\sR$ and a deficiency set 
$\sD = \{p\in\sP \mid \phi(\Z_p) \neq \phi(\Z_{p^{\infty}})$.
The field dimensional function can be defined as $d(p)=\phi(\Z_p)$ and 
$d(0)=\phi(\Q)$. 
Thus the set $F_{\phi}=(\sS,\sD;d)$ is well defined.
On the other hand if we have a set $F=(\sS,\sD;d)$
where $\sD\subset\sS\subset\sP$ and 
$d\colon \sP\cup\{0\}\to\N\cup\{\infty\}$,
we can define a function
$\phi_F\colon \sigma\to\N\cup\{\infty\} $ by formulas: $\phi(\Z_p)=d(p)$, 
$\phi(\Q)=d(0)$,
$\phi(\Z_{p^{\infty}})=d(p)-\chi_{\sD}(p)$ and 
$\phi(\Z_{(p)})
= (\max\{d(0),d(p) - \chi_{\sD}(p)+1\})
\chi_{\sS}(p)+d(0)\chi_{\sP\setminus\sS}(p)$,
where $\chi_A$ denotes the characteristic function of a set $A\subset\sP$.

\begin{proposition}
The correspondence $\phi \to F_{\phi}$ defines a bijection between all 
functions $\phi\colon \sigma\to\N\cup\{\infty\}$ satisfying the Bockstein 
Inequalities BI1--BI6 and triples $F=(\sS,\sD;d)$ with 
$d(\sP\setminus\sS)=d(0)$. 
Its inverse is defined by the above correspondence $F \mapsto \phi_F$.
\end{proposition}

We denote the set of functions $\phi\colon \sigma\to\N\cup\{\infty\}$ 
satisfying the Bockstein inequalities by $\sB_+$ and the set of triples 
$(\sS,\sD;d)$ with the constrain $d(\sP\setminus\sS)=d(0)$ by $\sF_+$

On the class of all compacta we consider the following equivalence 
relation: $X\sim Y$ if and only if $\dim_GX=\dim_GY$ for all abelian 
groups $G$. An equivalence class under that relation is called a 
cohomological dimension type or briefly {\it cd-type}. We define zero 
cd-type as the type of 0-dimensional compacta. Every nonzero cd-type can 
be described by an element of $\sF_+$ as well as by an element of $\sB_+$.

\begin{definition}
We define two binary operations $\bplus$ and $\btimes$ on $\sF_+$ 
by the formulas:
$$(\sS_1,\sD_1;d_1) \bplus (\sS_2,\sD_2;d_2) 
= (\sS_1\cup\sS_2,\sD_1\cup\sD_2;d_1+d_2)$$
$$(\sS_1,\sD_1;d_1) \btimes (\sS_2,\sD_2;d_2) 
= (\sS_1\cap\sS_2,\sD_1\cap\sD_2;(d_1-d_1(0))(d_2-d_2(0))+d_1(0)d_2(0))$$
\end{definition}

\begin{proposition}
$F_1 \bplus F_2\in\sF_+$ and $F_1 \btimes F_2\in\sF_+$ for 
$F_1,F_2\in\sF_+$.
\end{proposition}

\begin{proof}
First
$(d_1+d_2)(\sP\setminus(\sS_1\cup\sS_2))
= (d_1+d_2)((\sP\setminus\sS_1)\cap(\sP\setminus\sS_2))
= d_1(0)+d_2(0)$.

Second, since $(d_1-d_1(0))(d_2-d_2(0))=0$ on
$(\sP\setminus \sS_1)\cup(\sP\setminus\sS_2)=\sP\setminus(\sS_1\cap\sS_2)$,
$d(\sP\setminus(\sS_1\cap\sS_2)) = d_1(0)d_2(0) = d(0)$.
\end{proof}

\begin{proposition}
The distributivity law holds for operations $\bplus$ and $\btimes$.
\end{proposition}

\begin{proof}
It is known that the distributivity law holds for $\cup$ and $\cap$. 
We omit an easy verification of the distributivity law for functions $d$.
\end{proof}

\begin{proposition}
The natural numbers $\N$ are imbedded into $\sF_+$ by homomorphism taking 
a number $n$ to $(\emptyset,\emptyset;n)$ where $n$ also denotes the 
corresponding constant function.
\end{proposition}

\begin{proof}
The proof is trivial.
\end{proof}

\begin{definition}
{\it The norm} of cd-type $F=(\sS,\sD;d)\in\sF_+$ is defined as
$$\| F\|=\sup_{\sP\cup\{0\}}\{d+\chi_{\sS\setminus\sD}\}.$$
\end{definition}

\begin{proposition}
Let $F\in\sF_+$ represent the cd-type of a compactum $X$,
then $\| X\|=\dim_{\Z}X$.
\end{proposition}

\begin{proof}
By Bockstein Theorem $\dim_{\Z}X=\sup\{\dim_{\Z_{(p)}}X \mid p\in\sP\}$. 
By Lemma 2.8
$$\dim_{\Z_{(p)}}X
= \begin{cases}
d(0)&			\text{if $p\in\sP\setminus\sS$,}\\
\max\{d(0),d(p)+1\}&	\text{if $p\in\sS\setminus\sD$,}\\
\max\{d(0),d(p)\}&	\text{if $p\in\sD$.}
\end{cases}
$$
Therefore, 
$\sup\{\dim_{\Z_{(p)}}X \mid p\in\sP\}
= \sup\{(d+\chi_{\sS\setminus\sD})(x) \mid x\in\sP\cup\{0\}\}$.
\end{proof}

On the set of functions $\sB_+$ there is the natural partial order $\leq$:
$$
\text{$\phi_1\leq\phi_2$ if and only if $\phi_1(G)\leq\phi_2(G)$ 
for all $G\in\sigma$.}
$$
Thus the bijection of Proposition 4.1 defines a partial order $\preceq$ 
on cd-types.

\begin{proposition}
Let $\phi_1,\phi_2\in\sB_+$, then $\phi$ defined as 
$\phi(G) = \max\{\phi_1(G), \phi_2(G)\}$ satisfies the Bockstein 
inequalities, i.e.\ $\phi\in\sB_+$.
\end{proposition}

\begin{proof}
Trivial.
\end{proof}

\begin{definition}
Let $F_1$ and $F_2$ be two cd-types, then we define the wedge 
$F_1\vee F_2$ as the cd-type corresponding to the function 
$\phi(G)=\max\{\phi_{F_1}(G),\phi_{F_2}(G)\}$.
\end{definition}

Proposition 4.6 is valid if one replaces the maximum by a supremum over 
an arbitrary index set. Thus an operation $\bigvee_{i\in J}F_i$ can be 
defined for any family $\{F_i \mid i\in J\}\subset\sF_+$.

\begin{proposition}
The distributivity law holds for $\vee$ and $\bplus$.
\end{proposition}

\begin{proposition}
For every family $\{F_i \mid i\in J\}\subset\sF_+$ there is a countable 
subset $J'\subset J$ such that 
$\bigvee_{i\in J}F_i = \bigvee_{i\in J'}F_i$.
\end{proposition}

\begin{proof}
Take an arbitrary group $G\in\sigma$. If the maximum of $\phi_{F_i}(G)$ 
is attained on some $i_G\in J$, we define $L_G=\{i_G\}$. If not, then 
there is a sequence $\{i^k_G\}$ such that 
$\lim_{k\to\infty}\phi_{F_{i^k_G}}(G)=\infty$. In that case we define
$L_G=\{i^k_G\}_{k\in\N}$. We do this for all groups $G\in\sigma$. Then we 
define $J'=\bigcup_{G\in\sigma}L_G$. Since $\sigma$ is a countable set and 
$L_G$ is countable for every $G\in\sigma$, the set $J'$ is countable.
\end{proof}
 
By $\delta_x$ we denote a characteristic function of one point set 
$\{x\}$, i.e. 
$$\delta_x(t)
= \begin{cases}
1&	\text{if $t=x$,}\\
0&	\text{if $t\neq x$.}
\end{cases}
$$
We define {\it Kuzminov's basis} as the set of the following cd-types:
\begin{align*}
\Phi(\Q,n)&
= (\sP,\emptyset; (n-1)\delta_0+1),\\
\Phi(\Z_{(p)},n)&
= (\sP\setminus\{p\},\emptyset;(n-1)(\delta_0+\delta_p)+1),\\
\Phi(\Z_p,n)&
=(\sP,\{p\};(n-1)\delta_p+1),\\
\Phi(\Z_{p^{\infty}},n)&
= (\sP,\emptyset; (n-2)\delta_p+1).
\end{align*}

Here we assume that $n>1$. Since all 1-dimensional compacta define the 
same cd-type, we let $\Phi(G,1)$ equal the cd-type of one-dimensional 
compacta for all $G\in\sigma$. For $G\neq\Z_p$ the singularity set 
consists of whole $\sP$ and hence all these $\Phi(G,n)$ belong to $\sF_+$. 
For $\Phi(\Z_p,n)$ the condition $d(\sP\setminus\sS)=d(0)$ turns into 
$d(p)=d(0)$ and it is easy to check that it holds. Hence, 
$\Phi(\Z_p,n)\in\sF_+$ too.

\begin{proposition}
For all $G\in\sigma$ and every $n$, $\| \Phi(G,n)\| = n$.
\end{proposition}

\begin{proof}
$\| \Phi(\Q,n)\| 
= \sup\{d(x) + \chi_{\sS\setminus\sD}(x) \mid x \in \sP \cup \{0\}\}
= \sup\{(n-1)\delta_0 + 1 + \chi_{\sP}\}
= \max\{n,2\}
= n$.

$\| \Phi(\Z_{(p)},n)\|
= \sup\{(n-1)(\delta_0+\delta_p) + 1 + \chi_{\sP\setminus\{p\}}\}
= \max\{n,2\}
= n$.

$\| \Phi(\Z_p,n)\| 
= \sup\{(n-1)\delta_p + 1 + \chi_{\sP\setminus\{p\}}\}
= \max\{n,2\}
= n$.

$\| \Phi(\Z_{p^{\infty}},n)\| 
= \sup\{(n-2)\delta_p + 1 + \chi_{\sP}\}
= \max\{n,2,1\}
= n$.
\end{proof}

\begin{proposition}
The field dimensional function $d$ has its maximum at $p$ for
cd-types $\Phi(\Z_{(p)},n)$, $\Phi(\Z_p,n)$ and 
$\Phi(\Z_{p^{\infty}},n)$. It has its maximum at 0 for $\Phi(\Q,n)$.
\end{proposition}

The proof is an easy observation.

\begin{theorem}
For any cd-type $F\in\sF_+$ there is a representation
$F = \bigvee\{\Phi(G,k_G) \mid G\in\sigma\}$.
\end{theorem}

\begin{proof}
Let $F=(\sS,\sD);d)$. If the norm of $F$ equals one, then we take
$k_G=1$ for all $G\in\sigma$. If the norm is greater than one, we
take $k_{\Q}=d(0)$,
$$ k_{\Z_{(p)}}
= \begin{cases}
d(p)&	\text{if $p\in\sP\setminus\sS$,}\\
1&	\text{otherwise}
\end{cases}
$$
$$k_{\Z_p}
= \begin{cases}
d(p)&	\text{if $p\in\sD,$}\\
1&	\text{otherwise}
\end{cases}
$$
and
$$
k_{\Z_{p^{\infty}}}
= \begin{cases}
d(p)+1&	\text{if $p\in\sS\setminus\sD$,}\\
1&	\text{otherwise.}
\end{cases}
$$
Then we consider a cd-type 
\begin{align*}
F'&
= \bigvee\{\Phi(G,k_G) \mid G\in\sigma\}\\
&
= \Phi(\Q,d(0)) 
\vee \bigvee_{p\in\sP\setminus\sS}\!\! \Phi(\Z_{(p)},d(p))
\vee \bigvee_{p\in\sD}\! \Phi(\Z_p,d(p))
\vee \bigvee_{p\in\sS\setminus\sD}\!\! \Phi(\Z_{p^{\infty}},d(p)+1).
\end{align*}
If $F'=(\sS',\sD';d')$, then in the view of Proposition 4.10,
\begin{align*}
d'(0)&
= \max\{d(0),d(0),1,1\} = d(0)
\quad\text{and}\\
d'(p)&
= \max\{1,d(p),d(p),d(p)+1-1\} = d(p).
\end{align*} 
Therefore $d'=d$. It is easy to verify that $\sD'=\sD$ and $\sS'=\sS$. 
Hence $F=F'$.
\end{proof}

\begin{definition}
{\it The inferior norm} $|F|$ of cd-type $F=(\sS,\sD;d)$ is 
defined as 
$$\min\{d(x) - \chi_{\sD}(x) \mid x\in\sP\cup\{0\}\}.$$
\end{definition}

\begin{proposition}
Let $F\in\sF_+$ and let $\phi_F\in\sB_+$ be its representative.
Then $\| F\| = \sup\{\phi_F(G) \mid G\in\sigma\}$ and
$|F| = \inf\{\phi_F(G) \mid G\in\sigma\}$.
\end{proposition}

\begin{proof}
Note that 
\begin{align*}
\sup\{\phi_F(G) \mid G\in\sigma\}&
= \sup\{\phi_F(\Z_{(p)}) \mid p\in\sP\}\\
&
= \max\left\{
	\sup\{\max\{d(0), d(p)+\chi_{\sS\setminus\sD}(p)\} \mid p\in\sS\},
	d(0)
\right\}\\
&
= \sup\{ d(x)+\chi_{\sS\setminus\sD}(x) \mid x\in\sP\cup\{0\} \}\\
&
= \|F\|.
\end{align*}

By virtue of Bockstein's inequalities,
$$
\inf\{\phi_F(G) \mid G\in\sigma\}
= \inf\{\phi_F(\Q),\phi_F(\Z_{p^{\infty}})\}
= \min\{d(0),d(x) - \chi_{\sD}(x)\}
= |F|.
$$
\end{proof}

\begin{lemma}
For any two cd-types $F_1$ and $F_2$ there are inequalities:
$| F_1 | + \| F_2 \| 
\leq \| F_1 \bplus F_2 \| 
\leq \| F_1 \| +\| F_2 \|$.
\end{lemma}

\begin{proof}
We may assume that all norms are finite. 
Suppose that the norm
$\| F_2 \| = \sup\{d_2(y)+\chi_{\sS_2\setminus\sD_2}(y)\}$ 
is achieved at $x$. Then it suffices to show that
$$
d_1(x)
- \chi_{\sD_1}(x)
+ d_2(x)+\chi_{\sS_2\setminus\sD_2}(x)
\leq d_1(x)
+ d_2(x)
+ \chi_{(\sS_1\cup\sS_2)\setminus(\sD_1\cup\sD_2)}(x).
$$ 
This is equal to the inequality
$\chi_{\sS_2\setminus\sD_2}
\leq \chi_{_{(\sS_1\cup\sS_2)\setminus(\sD_1\cup\sD_2)}}
+ \chi_{_{\sD_1}}$ 
which follows from the fact that 
$\sS_2 \setminus \sD_2
\subset \sS_2 \setminus (\sD_2\cup\sD_1) \cup \sD_1 
\subset \left( (\sS_2\cup\sS_1) \setminus (\sD_2\cup\sD_1) \right) 
\cup \sD_1$.

Assume that $\|F_1 \bplus F_2\| $ is achieved at $x$, then
$\|F_1 \bplus F_2\|
= d_1(x)
+ d_2(x)
+ \chi_{(\sS_1\cup\sS_2) \setminus (\sD_1\cup\sD_2)}(x)$.
Because of the inclusion 
$(\sS_1\cup\sS_2)\setminus(\sD_1\cup\sD_2)
\subset (\sS_1\setminus\sD_1)\cup(\sS_2\setminus\sD_2)$, 
we have the inequality
$\chi_{(\sS_1\cup\sS_2)\setminus(\sD_1\cup\sD_2)}
\leq \chi_{\sS_1\setminus\sD_1}
+ \chi_{\sS_2)\setminus\sD_2)}$.
Then 
$$\|F_1 \bplus F_2\|
\leq d_1(x)
+ \chi_{\sS_1\setminus\sD_1}(x)
+ d_2(x)
+ \chi_{\sS_2)\setminus\sD_2)}(x)
\leq \|F_1\| + \|F_2\|.$$
\end{proof}

By $\sF$ we denote the set of all triples $(\sS,\sD;d)$, where 
$\sD\subset\sS\subset\sP$ are subsets of primes and 
$d\colon \sP\cup\{0\}\to\Z\cup\{-\infty,+\infty\}$ has the
property $d(\sP\setminus\sS)=d(0)$. As one can see $\sF$ is the natural 
extension of $\sF_+$. All operations $\vee$, $\bplus$, $\btimes$ as well 
as the partial order $\preceq$ can be extended to $\sF$. The notions of 
the norm $\|\ \|$ and the inferior norm $|\ |$ can be defined on $\sF$ 
without changes. All propositions proven for $\sF_+$ can be repeated 
without changes for $\sF$. Similarly one can extend $\sB_+$ to $\sB$ 
together with the bijection $\sF \to \sB$.

\begin{definition}
A conjugation $\bar F$ of $F=(\sS,\sD;d)\in\sF$ is defined as
$$(\sS,\sS\setminus\sD;-d).$$
\end{definition}

It is clear that $\bar F\in\sF$.

\begin{proposition}
For every $F=(\sS,\sD;d)\in\sF$,
\begin{enumerate}
\item $\bar{\bar F}=F$,
\item $F \bplus \bar F=(\sS,\sS;0)$,
\item $\|F \bplus \bar F\|=0$.
\end{enumerate}
\end{proposition}

\begin{proof}
The statements (1),(2) are obvious; (3) follows from (2) and the 
definition of the norm.
\end{proof}

\begin{lemma}
The conjugate $\bar F$ of $F$ is the maximal element with respect to the
order $\preceq$ in the set $\{F'\in\sF \mid \| F \bplus F'\|\leq 0\}$.
\end{lemma}

\begin{proof}
First, from Proposition 4.14(3) we can see that 
$\bar F \preceq \bigvee\{F'\in\sF \mid \| F \bplus F' \| \leq 0\}$.
Then we show that $\bar F\succeq F'$ for every $F'$ having the property
$\|F \bplus F'\|\leq 0$. Let $F=(\sS,\sD;d)$ and $F'=(\sS',\sD';d')$.
The inequality $\|F \bplus F'\|\leq 0$ implies that
\begin{equation*}\tag{$\ast$}
d' + \chi_{(\sS\cup\sS')\setminus(\sD\cup\sD')}
\leq -d.
\end{equation*}
Then $d'\leq d$ and hence, $\phi_{F'}(\Z_p)\leq\phi_{\bar F}(\Z_p)$ and
$\phi_{F'}(\Q)\leq\phi_{\bar F}(\Q)$. 
We recall that $\phi_F$ stands for the function from $\sB$ corresponding 
to $F$ under extended bijection from Proposition 4.1.
We note that 
$\sS\setminus\sD
\subset \sS\setminus(\sD\cup\sD')\cup\sD'
\subset (\sS\cup\sS')\setminus(\sD\cup\sD')\cup\sD'$. 
Hence, 
$\chi_{\sS\setminus\sD}
\leq \chi_{(\sS\cup\sS')\setminus(\sD\cup\sD')} + \chi_{\sD'}$. 
Therefore by this and ($\ast$) we have the following inequality:
\begin{equation*}\tag{$\ast\ast$}
d'- \chi_{_{\sD'}}
\leq -d - \chi_{(\sS\cup\sS')\setminus(\sD\cup\sD')} - \chi_{\sD'} 
\leq -d-\chi_{\sS\setminus\sD}.
\end{equation*}
Hence, $\phi_{F'}(\Z_{p^{\infty}})\leq\phi_{\bar F}(\Z_{p^{\infty}})$.

To treat the group $\Z_{(p)}$ we consider three cases. 

(1) $p\in\sP\setminus\sS'$,
then $\phi_{F'}(\Z_{(p)})=d'(0)\leq -d(0)\leq\phi_{\bar F}(\Z_{(p)})$. 
Here we applied ($\ast$) and BI4. 

(2) $p\in\sS'\cap\sS$, then
\begin{align*}
\phi_{F'}(\Z_{(p)})&
= \max\{d'(0),d'(p)-\chi_{\sD'}(p)+1\}\\
&
\leq
\max\{-d(0),-d(p)-\chi_{\sS \setminus\sD}(p)+1\}\\
&
= \phi_{\bar F}(\Z_{(p)}). 
\end{align*}
Here we applied both 
($\ast$) and 
($\ast\ast$).

(3) Finally, if $p\in\sS'\setminus\sS$, then the inclusion 
$\sS'\setminus\sS
\subset (S\cup\sS')\setminus(\sD\cup\sD')\cup\sD'$ 
implies that
$\chi_{(\sS\cup\sS')\setminus(\sD\cup\sD')}(p) + \chi_{_{\sD'}}(p) \geq 1$.
Then, 
$d'(p) - \chi_{\sD'}(p)+1
\leq d'(p) - \chi_{(\sS\cup\sS')\setminus(\sD\cup\sD')}(p)
\leq -d(p)
= -d(0)
= \phi_{\bar F}(\Z_{(p)})$.
Because of this and ($\ast$) we have that
$\phi_{F'}(\Z_{(p)})
= \max\{d'(p),d'(p)-\chi_{\sD"}(p)+1\}
\leq \phi_{\bar F}(\Z_{(p)})$.

Thus, $\bar F = \bigvee\{F'\in\sF \mid \| F \bplus F' \| \leq 0\}$.
\end{proof}

\section{Realization theorem}

The main result of this Section is the following:

\begin{theorem}[Realization Theorem]
For every cd-type $F=(\sS,\sD;d)\in\sF_+$ there exists a compactum $X$
such that $F_X=(\sS_X,\sD_X;d_X)=(\sS,\sD;d)$. Moreover, 
$F_X \bplus F_Y=F_{X\times Y}$ and $F_X \vee F_Y = F_{X\vee Y}$.
\end{theorem}

Thus the name `cd-types' for elements of $\sF_+$ is justified.

A compactum representing a fundamental cd-type $\Phi(G,n)$ is called a 
{\it fundamental compactum} of type $(G,n)$. The notation for this is 
$X\in F(G,n)$. A fundamental compacta have the following cohomological 
dimension with respect to groups from $\sigma$:

\begin{figure}[h]
\begin{tabular}{r|c|c|c|c|c|c|c}
&
$\Z_{(p)}$& 
$\Z_p$&
$\Z_{p^\infty}$&
$\Q$&
$\Z_{(q)}$&
$\Z_q$&
$\Z_{q^\infty}$\\
\hline
$F(\Q,n)$&
$n$&
$1$&
$1$&
$n$&
$n$&
$1$&
$1$\\
\hline
$F(\Z_{(p)},n)$&
$n$&
$n$&
$n$&
$n$&
$n$&
$1$&
$1$\\
\hline
$F(\Z_p,n)$&
$n$&
$n$& 
$n{-}1$&
$1$& 
$1$&
$1$&
$1$\\
\hline
$F(\Z_{p^\infty},n)$&
$n$&
$n{-}1$&
$n{-}1$& 
$1$& 
$1$&
$1$&
$1$\\
\end{tabular}
\caption{Cohomological dimension w.r.t.\ groups from $\sigma$}
\end{figure}

Here $p,q$ are primes, $q$ runs over all primes $\ne p$.

Let $h$ be a reduced homology (or cohomology) theory. A map between two 
topological spaces is called $h_*$-essential (or $h^*$-essential) if it 
induces nonzero homomorphism in $h$-homologies (or $h$-cohomologies). If 
one of the spaces is not a CW-complex, then we consider the \v{C}ech 
extension $\check h$. We recall that a cohomology theory $h^*$ is called 
{\it continuous} if for every direct limit 
$L=\limto \{L_i;\lambda^i_{i+1}\}$ of finite CW-complexes the formula 
$h^*(L)=\limleftarrow h^*(L_i)$ holds. We note that a cohomology 
$h^*(\ ;F)$ with coefficients in a field $F$ is continuous.

In this section we give a proof of Realization Theorem based on the 
following general theorem.

\begin{theorem}
Let $P$ and $K$ be simplicial complexes and assume that $K$ is countable 
complex.
Let $h_*$ ($h^*$) be a reduced generalized homology (continuous 
cohomology) theory.
If $h_n(P)\neq 0$ ($h^n(P)\neq 0$) and $h_k(K)=0$ ($h^k(K)=0$) for
all $k<n$, then there exist a compactum $X$, having the property $K\in AE(X)$,
and an $h_n$-essential ($h^n$-essential) map $f\colon X\to P$. 
\end{theorem}

\begin{corollary}
For every $n\in\N$ there are $n$-dimensional fundamental compacta of types 
$(\Q,n)$, $(\Z_p,n)$, $(\Z_{(p)},n)$ and $(\Z_{p^{\infty}},n)$ for all 
primes $p$ and any $n$.
\end{corollary}

\begin{proof}
To realize the type $(\Q,n)$ we take $P=S^n$ and for $K$ we take the 
wedge of an Eilenberg-Maclane complex and $n$-sphere 
$K=K(\bigoplus_{p\in\sP}\Z_p,1) \vee S^n$ and we consider a continuous 
cohomology $h^*=H^*(\ ;\Q)$.
We note that $h^*(K)=0$ for $k<n$. By Theorem 5.2 there exists a 
compactum $X$ having nontrivial $n$-dimensional rational \v{C}ech 
cohomology.
The property $K\in AE(X)$ implies that 
$K(\bigoplus_{p\in\sP}\Z_p,1)\in AE(X)$ and $S^n\in AE(X)$. The second 
condition implies that $\dim X\leq n$. The first implies the inequality 
$\dim_{\bigoplus\Z_p}X\leq 1$ by virtue of Theorem 1.1.
Corollary 1.7 implies that $\dim_{\Z_p}X\leq 1$ for all $p$.
Since $X$ is not 0-dimensional (see 1.3), $\dim_{\Z_p}X=1$ for all $p$.
Therefore by BI1, $\dim_{\Z_{p^{\infty}}}X=1$. 
The equality $\dim_{\Q}X=n$ follows form the $n$-dimensionality of $X$ 
and the fact that $X$ has nontrivial rational $n$-dimensional cohomology.
The Bockstein inequality BI4 imply that $\dim_{\Z_{(p)}}X=n$ for all $p$.
Now according to the above table $X$ has cohomological dimensions as 
$F(\Q,n)$, hence $X$ is of type $(\Q,n)$.

For the type $(\Z_p,n)$ we take $P=S^n$, 
$K=K(\Z[\frac{1}{p}],1) \vee S^n$ and $h^*=H^*(\ ;\Z_p)$.
Then the compactum $X$ of Theorem 5.2 has a cohomological dimension
$\dim_{\Z[\frac{1}{p}]}X=1$ and the covering dimension $\dim X\leq n$. 
By virtue of Bockstein Theorem, $\dim_{\Z_{(q)}}X=1$ for prime $q\neq p$. 
By BI1,3,4 we have $\dim_{\Z_{q^{\infty}}}X= \dim_{\Z_q}X=\dim_{\Q}X=1$. 
Since $X$ has nontrivial $n$-dimensional cohomology with 
$\Z_p$-coefficients, we have $\dim_{\Z_p}X\geq n$. The equality holds, 
since $X$ is $n$-dimensional. Since $\dim X\leq n$, by BI3 
$\dim_{\Z_{(p)}}X=n$. We may assume that $n>1$. Then $X$ is $p$-singular. 
Since $\dim_{\Z_{(p)}}X= \max\{\dim_{\Q}X,\dim_{\Z_{p^{\infty}}}X+1\}$ for 
$p$-singular compacta, we have that $\dim_{\Z_{p^{\infty}}}X=n-1$. Thus, 
according to the above table $X\in F(\Z_p,n)$.

For the type $(\Z_{(p)},n)$ we take $P=K(\Z_{p^{\infty}},n)$, 
$K=K(\bigoplus_{q\neq p}\Z_q,1) \vee S^n$, and $h_*=H_*(\ ;\Z_{(p)})$. 
Since $\Z_{p^{\infty}}\otimes\Z_{(p)}\neq 0$, by Hurewicz theorem and the 
Universal Coefficient Formula, $H_n(P;\Z_{(p)})\neq 0$. 
Note that $h_k(K)=0$ for all $k<n$.
Apply Theorem 5.2 to obtain a compact $X$ and a map $f\colon X\to P$ with 
the certain properties. The property $K\in AE(X)$ implies, by virtue
Theorem 1.1, equalities $\dim_{\Z_q}X=\dim_{\Z_{q^{\infty}}}X=1$ and 
the inequality $\dim X\leq n$. The essentiality of the map $f\colon X\to 
K(\Z_{p^{\infty}},n)$ gives nontrivial element in cohomology 
$\check H^n(X;\Z_{p^{\infty}})$. Hence $\dim_{\Z_{p^{\infty}}}X=n$. The 
inequalities BI1 and BI3 imply that $\dim_{\Z_{(p)}}X=n$. Hence, by Lemma 
2.6, $X$ is $p$-regular. Therefore, $\dim_{\Q}X=n$. We assume that $n>1$, 
since any 1-dimensional compactum can serves as $F(\Z_{(p)},1)$. Then $X$ 
is $q$-singular for all prime $q\neq p$.
Then by Lemma 2.8 it follows that $\dim_{\Z_{(q)}}X=\dim_{\Q}X=n$. Thus, 
$X$ has cohomological dimensions with respect to groups from $\sigma$ as 
it prescribed for $F(\Z_{(p)},n)$ by the table in the beginning of this 
section.

For the type $(\Z_{p^{\infty}},n)$ we take $P=S^n$,
$K = K(\Z[\frac{1}{p}],1) \vee K(\Z_p,n-1) \vee S^n$ 
and $h_*=H_*(\ ;\Z_{p^{\infty}})$.
Note that $ H_k(K;\Z_{p^{\infty}})=0$ for $k<n$.
We apply Theorem 5.2 to obtain a compactum $X$ having the property $K\in 
AE(X)$ and an essential map onto $n$-dimensional sphere. This properties 
imply that $\dim X=n$, 
$\dim_{\Z_{(q)}}X
= \dim_{\Z_q}X
= \dim_{\Z_{q^{\infty}}}X
= \dim_{\Q}X=1$ 
and $\dim_{\Z_p}X\leq n-1$.
Since $\dim_{\Z_{(p)}}X=n$, by the Bockstein Alternative it follows that
$\dim_{\Z_{p^{\infty}}}X=n-1$. Then by BI1, $\dim_{\Z_p}X=n-1$. Then 
$X\in F(\Z_{p^{\infty}},n)$.
\end{proof}

\begin{definition}
An {\it extension problem} $(A,\alpha)$ on a topological space $X$ is a 
map $\alpha\colon A\to K$ defined on a closed subset $A\subset X$ with 
the range a CW-complex (or ANE). A {\it solution} of an extension problem 
$(A,\alpha)$ is a continuous extension $\bar\alpha\colon X\to K$ of a map 
$\alpha$. A {\it resolution} of an extension problem $(A,\alpha)$ is a map 
$f\colon Y\to X$ such that the induced extension problem 
$f^{-1}(A,\alpha) = (f^{-1}(A), \alpha\circ f\mymid_{\cdots})$ on $Y$ 
has a solution.
\end{definition}

Because of the Homotopy Extension Theorem the solvability of an extension 
problem $(A,\alpha)$ is an invariant of homotopy class of $\alpha$. A 
family of extension problems $\{(A_i,\alpha_i)\}_{i\in J}$ forms a {\it 
basis} if for every extension problem $(B,\beta)$ there is $i\in J$ such 
that $B\subset A_i$ and the restriction $\alpha_i\mymid_B$ is homotopic to 
$\beta$. In that case we say that $(A_i,\alpha_i)$ contains $(B,\beta)$.

In view of the Homotopy Extension Theorem the following Proposition is 
obvious.

\begin{proposition}
Suppose that a map $f\colon Y\to X$ resolves extension problems on $X$ 
from a given basis $\{(A_i,\alpha_i)\}_{i\in J}$. Then $f$ resolves all 
extension problems on $X$.
\end{proposition}

\begin{proposition}
Let $K$ be fixed. Let $X$ be the limit space of an inverse sequence of 
compacta $\{X_k,q_k^{k+1}\}$ and let $\{(A^k_i,\alpha^k_i)\}_{i\in J_k}$ 
be a basis of extension problems for every $k$. Then 
$\{(q_k^{\infty})^{-1}(A^k_i,\alpha^k_i)\mid k\in\N,i\in J_k\}$ is a 
basis of extension problems on $X$ where $q^{\infty}_k\colon X\to X_k$ 
denotes the infinite projection in the inverse sequence.
\end{proposition}

\begin{proof}
Since $K\in ANE$, for every extension problem $(A,\alpha)$ on $X$ there 
is a number $k$ and a map $\beta\colon q^{\infty}_k(A)\to K$ such that 
$\beta\circ q^{\infty}_k\mymid_A$ is homotopic to $\alpha$. Take a problem 
$(A^k_i,\alpha^k_i)$ containing $(q^{\infty}_k(A),\beta)$. Then 
$\alpha^k_i\mymid_{q^{\infty}_k(A)} \sim\beta$. The extension problem 
$(q^{\infty}_k)^{-1}(A^k_i,\alpha^k_i)$ contains the problem 
$(A,\alpha)$. 
\end{proof}

\begin{lemma}
For any extension problem $(A,\alpha\colon A\to K)$ on $X$ there is a 
resolution of it $g\colon Y\to X$ such that every preimage $g^{-1}(x)$ is 
a point or homeomorphic to $K$. If additionally $X$ and $K$ are 
simplicial complexes, $A$ is a subcomplex and $\alpha$ is a simplicial 
map, then the resolving map $g$ can be chosen simplicial.
\end{lemma}

\begin{proof}
Let $\pi\colon K\times I\to \cone(K)$ be the standard projection onto the 
cone. So, the preimage $\pi^{-1}(x)$ is one point set if $x$ is not the 
cone vertex, and it is homeomorphic to $K$ if $x$ is the cone vertex. We 
identify $K$ with the bottom of the cone $\cone(K)$. Since $\cone(K)\in 
ANE$, there is an extension $\bar\alpha\colon X\to \cone(K)$ of $\alpha$. 
We define $Y$ as a pullback of the diagram:
$$\begin{CD}
Y	@>>\gamma>	K\times I\\
@VVgV			@VV\pi V\\
X	@>>\bar\alpha>	\cone(K)\\
\end{CD}
$$
Then $\pr\circ\gamma\colon Y\to K$ is a solution of the extension problem 
$g^{-1}(A,\alpha)$ where $\pr\colon K\times I$ is the projection. Thus, 
the map $g\colon Y\to X$ resolves the problem $(A,\alpha)$. Since $g$ is 
parallel to $\pi$ in the pullback diagram, $g$ has the same set of 
topological types of point preimages, i.e.\ the set consisting of the one 
point space and $K$.

If $\alpha$ is simplicial and $A\subset X$ is a subcomplex, then we consider
the natural structure of a simplicial complex on the cone $\cone(K)$. Take 
all vertices of $X$ which do not belong to $A$ to the cone vertex and 
thus, define a simplicial extension $\bar\alpha$ of $\alpha$. Consider a 
product simplicial structures on $K\times I$ and $X\times(K\times I)$. 
Then the projection $\pi\colon K\times I\to \cone(K)$ is a simplicial 
map. Consider the induced triangulation on the pullback 
$L\subset X\times(K\times I)$. The map $g$ is simplicial with respect to 
that triangulation.
\end{proof}

\begin{proposition}
Let $X$ be the limit space of an inverse sequence $\{X_k;q^{k+1}_k\}$ and 
let $\{(A^k_i,\alpha^k_i)\}_{i\in J_k}$ be a basis of extension problems 
for each $k$.
Assume that $q^{\infty}_k$ resolves all problems $(A^k_i,\alpha^k_i)$ for 
all $k$.
Then $K\in AE(X)$.
\end{proposition}

\begin{proof}
According to Proposition 5.5 $X$ has a basis of solvable extension 
problems. Then by Proposition 4.4 all extension problems on $X$ have 
solutions. It means that $K\in AE(X)$.
\end{proof}

\begin{remark*} 
If a map $f\colon Y\to X$ resolves some extension problem $(A,\alpha)$ on 
$X$, then for any map $g\colon Z\to Y$ the composition $f\circ g$ resolves 
$(A,\alpha)$.
\end{remark*}

\begin{lemma}
Let $g\colon L\to M$ be a simplicial map onto a finite dimensional 
complex $M$ and let $h_*$ be a reduced homology theory such that 
$h_k(g^{-1}(x))=0$ for all $k<n$ ($k\in\Z$). Then $g$ induces an 
isomorphism $g_*\colon h_k(L)\to h_k(M)$ for $k<n$ and an epimorphism for 
$k=n$.
\end{lemma}

\begin{proof}
We prove it by induction on $m=\dim M$.

If $\dim M=0$, then Lemma holds.

Let $\dim M=m>0$. 
We denote by $A$ a regular neighborhood in $M$ of $(m-1)$-dimensional
skeleton $M^{(m-1)}$. Since the map $g\colon L\to M$ is simplicial, 
$g^{-1}(A)$ has a deformation retraction onto $g^{-1}(M^{(m-1)}$. By the 
induction assumption Lemma holds for 
$g\mymid_{\cdots}\colon g^{-1}(M^{(m-1)})\to M^{(m-1)}$. 
Hence, the conclusion of Lemma holds for 
$g\mymid_{\cdots}\colon g^{-1}(A)\to A$. 
We define $B=M\setminus \Int A$, i.e.\ $B$ is a union of disjoint 
$m$-dimensional PL-cells $B=\bigcup B_i$. 
Since $g$ is simplicial, $g^{-1}(B_i)\simeq g^{-1}(c_i)\times B_i$ where 
$c_i\in B_i$. 
Therefore the conclusion of Lemma holds for 
$g\mymid_{\cdots}\colon g^{-1}(B)\to B$. 
Note that $\dim(A\cap B)=m-1$ and, hence,Lemma holds for 
$g\mymid_{\cdots}\colon g^{-1}(A\cap B)\to A\cap B$. 
The Mayer-Vietoris sequence for the triad $(A,B,M)$ produces the following 
diagram:
$$
\begin{CD}
h_k(A'\cap B')	@>>>	h_k(A')\oplus h_k(B')	@>>>	h_k(L)	@>>>	h_{k-1}(A\cap B')	@>>>\\
@VVV			@VVV				@Vg_*VV		@VVV			\\
h_k(A\cap B)	@>>>	h_k(A)\oplus h_k(B)	@>>>	h_k(M)	@>>>	h_{k-1}(A\cap B)	@>>>\\
\end{CD}
$$
Here $A'=g^{-1}(A)$ and $B'=g^{-1}(B)$.
The Five Lemma implies that $g_*$ is an isomorphism for $k<n$. The 
epimorphism version of the Five Lemma implies that $g_*$ is an epimorphism 
for $k=n$.
\end{proof}

\begin{lemma}
Let $g\colon L\to M$ be a simplicial map onto a finite dimensional 
complex $M$ and let $h^*$ be a reduced cohomology theory such that 
$h^k(g^{-1}(x))=0$ for all $k<n$ ($k\in\Z$). Then $g$ induces an 
isomorphism $g^*\colon h^k(M)\to h^k(L)$ for $k<n$ and a monomorphism for 
$k=n$.
\end{lemma}

\begin{proof}
We can apply the argument of Lemma 5.8 with the only difference that at 
the very end we should apply the monomorphism version of the Five Lemma.
\end{proof}

\begin{proof}[Proof of Theorem 5.2]
Since $h_n(P)\neq 0$ ($h^n(P)\neq 0$), there exists a finite subcomplex 
$P_1\subset P$ such that the inclusion is $h_n$-essential 
($h^n$-essential). For cohomology this follows from the continuity of 
$h^*$, for homology it follows from the fact that every homology
has a compact support. 
We construct $X$ as the limit space of an inverse sequence of polyhedra
$\{P_k;q_k^{k+1}\}$ where $f\colon X\to S^n$ will be the composition of 
$q^{\infty}_1$ and the inclusion $P_1\subset P$. 
We construct this sequence by induction on $k$ such that
\begin{enumerate}
\item 
for every $k$ there is fixed some countable basis of extension problems 
$\sA^k=\{(A^k_i,\alpha^k_i)\}$ on $P_k$,
\item 
for every $k$ some nonzero element $a_k\in h_n(P_k)$ ($a_k\in h^n(P_k)$) 
is fixed such that $(q^{k+1}_k)_*(a_{k+1})=a_k$ 
($q^{k+1}_k)^*(a_k)=a_{k+1}$) 
for all $k$.
\item 
for every problem $(A^k_i,\alpha^k_i)\in\sA^k$ there is $j>k$ such that 
$q^j_k$ resolves it.
\end{enumerate}
If we manage to construct such a sequence, then by Proposition 5.7 
$K\in AE(X)$.
The property (2) would imply that $f$ is $h_n$-essential 
($h^n$-essential).
Thus, Theorem 5.2 would be proven. 

Enumerate all prime numbers $2=p_1<p_2<p_3<\dots<p_k<\dots$. 
We are going to work with homology first. 
We fix some element $a_1\in h_*(P_1)$ which goes to a nonzero element 
$a\in h_n(P)$. 
Denote by $\tau_1$ a triangulation on $P_1$ and by $\beta^k\tau$ $k$-th 
barycentric subdivision of $\tau$.
There are only countably many subpolyhedra in $P_1$ with respect to all 
subdivisions $\beta^k\tau$. 
Since the set of homotopy class $[L,K]$ is countable for every compact 
$L$, we have only countably many different extension problems 
$(A,\alpha)$ defined on those subpolyhedra. 
Denote the set of all these extension problems $(L,\alpha)$ on $P_1$ with 
simplicial maps $\alpha$ by $\sA^1$. 
Since $K\in ANE$, it easy to show that $\sA^1$ form a basis of extension
problems on $P_1$. 
We enumerate elements of $\sA^1$ by all powers of 2. 
Let $N\colon \sA^1\to\N$ be enumeration function. 
Then we consider an extension problem from $\sA^1$ having number one in 
our list and resolve it by a simplicial map $g\colon L\to P_1$ by means 
of Lemma 5.6.
By Lemma 5.8 $g_*\colon h_n(L)\to h_n(P_1)$ is an epimorphism. 
Take $a_2'\in h_n(L)$ such that $g_*(a_2')=a_1$. 
Since a homology $a_2'$ has a compact support, there is a finite 
subcomplex $P_2\subset L$ and an element $a_2\in h_n(P_2)$ which goes to
$a_2'$ under the inclusion homomorphism. 
We define a bonding map $q^2_1\colon P_2\to P_1$ as the restriction 
$f\mymid_{P_2}$ of $f$ onto $P_2$. 
Then the condition (2) holds:
$(q^2_1)_*(a_2)=a_1$. 
Then we define a countable basis $\sA^2=\{(A^2_i,\alpha^2_i)\}$ of 
extension problems such that every $A^2_i$ is a subcomplex of $P_2$ with 
respect to iterated barycentric subdivision of the triangulation on $P_2$. 
Enumerate elements of $\sA^2$ by all numbers of the form $2^k3^l$ with 
$k\geq 0$ and $l>0$. 
Lift all the problems from the list $\sA^1$ to a space $P_2$, i.e.\ 
consider $(q^2_1)^{-1}(\sA^1)$. 
Thus the family $(q^2_1)^{-1}(\sA^1)\cup\sA^2$ is enumerated by all 
numbers of the form $2^k3^l$, let 
$N\colon (q^2_1)^{-1}(\sA^1)\cup\sA^2\to\N$ be the enumeration function.
Now consider the extension problem having number 2 in updated list and 
apply the whole staff from the above to obtain $P_3$ and so on.

Thus, all problems in $\sA^k$ will be enumerated by numbers of the form 
$p_1^{l_1}p_2^{l_2} \cdots p_k^{l_k}$ with $l_k>0$. 
Since $k\leq p_k$, we have 
$k\in N((q^k_1)^{-1}(\sA^1)\cup (q^k_2)^{-1}(\sA^2)\cup \cdots \sA^k)$.
Hence we can keep going for any $k$. 
As the result of this construction we have that if a problem 
$(A^l_i,\alpha^l_i)$ has number $k$, then $l\leq k$ and the problem
is resolved by $q^{k+1}_l$. Thus, the conditions (1)--(3) hold.

If we consider a continuous cohomology $h^*$ instead of homology, we 
apply Lemma 5.9 instead of Lemma 5.8. Then we apply the continuity to get 
a finite subcomplex $P_2$.
The rest is the same.
\end{proof}

\begin{proof}[Proof of Theorem 5.1]
By Corollaries 5.3 and 5.9 we can realize by compacta all fundamental 
cd-types.
According to Theorem 4.11 an arbitrary cd-type $F\in\sF_+$ can be 
presented as $\bigvee\{\Phi(G,k_G)\mid G\in\sigma\}$. Then the one-point 
compactification of the disjoint union of fundamental compacta 
$\bigcup\{F(G,k_G)\mid G\in\sigma\}$ realizes the cd-type $F$.

The property $F_X \bplus F_Y=F_{X\times Y}$ follows from the definition 
of the operation $\bplus$ and Lemmas 3.3, 3.9 and 3.13. 
The equality $F_X\vee F_Y=F_{X\vee Y}$ follows from the formula 
$\dim_GX\vee Y=\max\{\dim_GX,\dim_GY\}$ which is the consequence of 
Theorem 1.5.
\end{proof}

\section{Test spaces}

Given an abelian group $G$, a compactum $X$ is said to be $G$-{\it 
testing space} for some class of compacta $\sC$ if for all spaces 
$Y\in\sC$ the following equality holds:
$$\dim_G=\dim(X\times Y)-\dim X.$$

\begin{theorem}
For any abelian group $G$ and any natural number $n$, there exists an 
$n$-dimensional compactum $T_n(G)$ which is a $G$-testing space for class 
of compacta $Y$ satisfying the inequality $\dim Y-\dim_GY<n$.
\end{theorem}

The following is the table of the dimension of the product of two 
fundamental compacta with $n\geq m$:

\begin{figure}[h]
\begin{tabular}{r|c|c|c|c|c|c|c}
&
$(\Z_{(p)},n)$&
$(\Z_p,n))$&
$(\Z_{p^\infty},n))$&
$(\Q,n))$&
$(\Z_{(q)},n))$&
$(\Z_q,n))$&
$(\Z_{q^\infty},n)$\\
\hline
$F(\Q,m)$&
$m{+}n$&
$n{+}1$&
$n{+}1$&
$m{+}n$&
$m{+}n$&
$n{+}1$& 
$n{+}1$\\
\hline
$F(\Z_{(p)},n)$&
$m{+}n$&
$n{+}1$&
$n{+}1$&
$m{+}n$&
$m{+}n$&
$n{+}1$&
$n{+}1$\\
\hline
$F(\Z_p,n)$&
$m{+}n$&
$m{+}n$&
$n{+}1$&
$n{+}1$&
$m{+}n$&
$n{+}1$&
$n{+}1$\\
\hline
$F(\Z_{p^\infty},n)$&
$m{+}n$&
$m{+}n{-}1$&
$m{+}n{-}1$&
$n{+}1$&
$m{+}n$&
$m{+}n{-}1$&
$n{+}1$\\
\end{tabular}
\caption{Dimension of the product of two fundamental compacta with
$n\geq m$}
\end{figure}
 
Here $q\neq p$. 
We leave to the reader the computations in this table. 
They are based on Proposition 4.5 and the formula 
$F_X \bplus F_Y = F_{X\times Y}$.
The result of calculations, presented in the table, can be summarized in 
the following formula ($n\geq m$):
$$\dim (F(G,n)\times F(G',m))=\dim_GF(G',m)+n.$$

\begin{proposition}
For any fundamental cd-type $\Phi(G,n)$ and any other cd-type $F$ there 
is the formula:
$$\|F \bplus \Phi(G,n)\|
= \begin{cases}
\max\{\|F\|+1, n+\phi_F(G)\}&	\text{if $\| F \| \geq n$,}\\
n+\phi_F(G)&			\text{if $\| F \| \leq n$.}\\
\end{cases}
$$
\end{proposition}

\begin{proof}
The function $\phi_F$ was defined in the beginning of \S4. The 
fundamental cd-types can be given via functions $\phi_F$ by means of the 
table of \S5.
If $F$ is a fundamental cd-type, then the result follows from the table.
In general case by Theorem 4.11 
$F = \bigvee\{\Phi(G',k_{G'}) \mid G'\in\sigma\}$.
Then 
\begin{align*}
\| F \bplus \Phi(G,n)\|&
= \sup\left\{ \| \Phi(G',k_{G'}) \bplus \Phi(G,n) \| \right\}\\
&
= \sup\left\{ \max\{k_{G'}\}+1, n+\phi_{\Phi(G',k_{G'})}(G)\} \right\}\\
&
= \max\left\{
\sup\{k_{G'}\}+1,n+\sup\{\phi_{\Phi(G',k_{G'})}(G)\}
\right\}\\
&
= \max\{\|F\|+1,n+\phi_F(G)\}.
\end{align*}
\end{proof}

\begin{proof}[Proof of Theorem 6.1]
We define $T_n(G)$ as a compactum representing the following cd-type 
$\bigvee\{\Phi(h,n) \mid H\in\sigma(G)\}$. 
Let us consider a compactum $X$ with $\dim X-\dim_GX<n$. 
If $\|F_X\|<n$, then by Proposition 6.2,
\begin{align*}
\dim(X\times T_n(G))&
= \| F_X \bplus \bigvee\{\Phi(H,n) \mid H\in\sigma(G)\} \|\\
&
= \sup\{\| F_X \bplus \Phi(H,n) \| \mid H\in\sigma(G)\}\\
&
= n + \phi_{F_X}(G)\\
&
= n + \dim_GX. 
\end{align*}
So, the testing formula holds.
If $\|F_X\|\geq n$, then by Proposition 6.2,
\begin{align*}
\dim(X\times T_n(G))&
= \|F_X \bplus \bigvee\{\Phi(H,n) \mid H\in\sigma(G)\}\|\\
&
= \sup\{\| F_X \bplus \Phi(H,n)\| \mid H\in\sigma(G)\}\\
&
= \max\{\| F_X \| + 1, n+\phi_{F_X}(G)\}\\
&
= \max\{\dim X+1, n+\dim_GX\}. 
\end{align*}
Since $\dim X-\dim_G\leq n-1$, we have that $\dim X+1\leq n+\dim_GX$ and 
hence, $\dim(X\times T_n(G))=n+\dim_GX$.
\end{proof}

\begin{theorem}
For two finite dimensional compacta $X$ and $Y$ the following conditions 
are equivalent:
\begin{enumerate}
\item $X$ and $Y$ have the same cd-type: $F_X=F_Y$,
\item for every compactum $Z$ there is the equality 
$\dim(X\times Z)=\dim(Y\times Z)$.
\end{enumerate}
\end{theorem}

\begin{proof}
The cd-type $F_{X\times Z}$ equals $F_X \bplus F_Z$ and hence depends 
only on cd-type of $X$. Therefore, 
$\dim(X\times Z)
= \dim_{\Z}(X\times Z)
= \|F_{X\times Z}\|
= \|F_X \bplus F_Z\|
= \|F_Y \bplus F_Z\|
= \dim(Y\times Z)$.

Given group $G$, we take $Z=T_n(G)$ with $n>\max\{\dim X,\dim Y\}$. Then 
by the testing equality we obtain: 
$\dim_GX+n
= \dim(X\times T_n(G))
= \dim(Y\times T_n(G))
= \dim_GY+n$.
Hence, $\dim_GX=\dim_GY$.
\end{proof}

\begin{corollary}
A finite dimensional compactum $X$ is dimensionally full-valued if and 
only if $\dim(X\times Z)=\dim X+\dim Z$ for all compacta $Z$.
\end{corollary}

\begin{proof}
Let $n=\dim X$, take $Y=I^n$. If $X$ is dimensionally full-valued, then 
it has the same cd-type as an $n$-cube $Y$. Since $\dim(Y\times Z)=n+\dim 
Z$, then by Theorem 6.3, $\dim(X\times Z)=n+\dim Z=\dim X+\dim Z$. If 
$\dim(X\times Z)=\dim X+\dim Z=n+\dim Z$, then 
$\dim(X\times Z)=\dim(Y\times Z)$ for all compacta $Z$. Hence 
$\dim_GX=\dim_GI^n=n$ for all $G$. Therefore, $X$ is dimensionally 
full-valued.
\end{proof}

The test spaces are very useful for extending some results of the 
Dimension Theory to a cohomological dimension.

\begin{theorem}
Let $f\colon X\to Y$ be a map between compacta and let $G$ be an abelian 
group.
\begin{enumerate}
\item If $f$ is $(k+1)$-to-1 map, i.e.\ the number of points in 
$f^{-1}(x)\leq k+1$, then $\dim_GX\geq \dim_GY-k$,
\item If $f$ is an open and all preimages of point are countable, then 
$\dim_GX=\dim_GY$.
\end{enumerate}
\end{theorem}

\begin{proof}
(1). Consider a map $f\times \id\colon X\times T_n(G)\to Y\times T_n(G)$ 
for large enough $n$ and apply the Hurewicz Theorem to obtain 
$\dim(X\times T_n(G))\geq \dim(Y\times T_n(G))-k$. Then the inequality 
$\dim_GX\geq \dim_GY-k$ follows from the $G$-testing formula.

(2). Consider the same map as in (1) and apply the Alexandroff Theorem to 
obtain the result.
\end{proof}

Let $F\in\sF$ be a cd-type, denote by $kF$ the sum $\bplus_{i=1}^kF$.
We recall that the integers $\Z$ are naturally imbedded in $\sF$.
For every $n\in\Z$ we denote by $\tilde n$ the image of $n$ in $\sF$
under that imbedding.

\begin{proposition}
Let $G\in\sigma$, then
\begin{enumerate}
\item $2\Phi(G,n)=\Phi(G,2n)\vee \tilde 2$ and in the general case
$k\Phi(G,n) = \Phi(G,kn) \vee \tilde k$ if 
$G\neq\Z_{p^{\infty}}$,
\item $2\Phi(G,n) = \Phi(G,2n-1) \vee \tilde 2$ and
$k\Phi(G,n) = \Phi(G,kn-k+1) \vee \tilde k$ if 
$G=\Z_{p^{\infty}}$.
\end{enumerate}
\end{proposition}

\begin{proof}
Let $\Phi(G,n)=(\sS_n,\sD_n;d_n)$, then 
$2\Phi(G,n)
= (\sS_n,\sD_n;2d_n)
= (\sS_{2n},\sD_{2n};2d_n)$.
Let 
$\Phi(G,2n)\vee\tilde 2=(\sS',\sD';d')$. Then the field function
$d'$ of $\Phi(G,2n)\vee\tilde 2$ is defined by the formula
$d'(x)=\max\{d_{2n}(x),2\}=2d_n(x)$. If $\Phi(G,2n)$ is 
$p$-regular, then $\Phi(G,2n)\vee\tilde 2$ is $p$-regular. If
$\Phi(G,2n)$ is $p$-singular, then $\Phi(G,2n)\vee\tilde 2$ is
$p$-singular provided $2n>2$. Hence, $\sS'=\sS_{2n}$. Similarly,
$\sD'=\sD_{2n}$.

The proof in the case of $k>2$ is not more difficult.

The difference in this case ($k=2$) is that the formula for 
$d'$ is the following $\max\{d_{2n-1}(x),2\}=2d_n(x)$. The rest of
the argument is the same.
\end{proof}

\begin{lemma}
Let $X$ be a fundamental compactum of the type $(G,n)$, 
$G\in\sigma$. Then for every $k$, the $k$-th power $X^k$ is a
$G$-testing space for the class of compacta $Y$ with 
$\dim Y-\dim_GY<n$.
\end{lemma}

\begin{proof}
By Proposition 6.6 the cd-type of a compactum $X^k$ is the same as
the cd-type of the union of $Z\coprod I^k$ where $Z$ is a 
fundamental compactum of the cd-type $(G,m)$ and $m=\dim X^k$.
Hence, $\dim(Y\times X^k)=\max\{\dim(Y\times Z),\dim Y+k\}$. By
Theorem 6.1 we can continue $=\max\{\dim_GY+m, \dim Y+k\}$. 
Since for $k>1$ the inequality $m-k\geq kn-k+1-k\geq 
n>\dim Y-\dim_GY$ holds, $\dim_GY+m\geq \dim Y+k$. Hence,
$\dim(Y\times X^k)=\dim_GY+m$.
\end{proof}

\begin{proposition}
Let $R$ be a principal ideal domain with unity $1\in R$, then
for no prime $p$, $\Z_{p^{\infty}}\in\sigma(R)$.
\end{proposition}

\begin{proof}
Assume that $\Z_{p^{\infty}}\in\sigma(R)$. 
Then it means that $p$-torsion subgroup $T=p-\Tor(R)$ is $p$-divisible.
Note that $T$ is an ideal in $R$. Therefore, $T=uR$ for some
$u\in R$. Since $1\in R$, it follows that $u\in T$. Let $p^k$ be
the order of $u$. Since $T$ is $p$-divisible, there is a quotient
$u/p^k\in T$. Then $u/p^k=uv$ for some $v\in R$. Hence,
$0=(p^ku)v=p^k(uv)=p^k(u/p^k)=u$. Contradiction.
\end{proof}

\begin{theorem}
Let $f\colon X\to Y$ be a continuous map between finite dimensional 
compacta. Then
\begin{enumerate}
\item $\dim_GX\leq \dim_GY+\max\{\dim f^{-1}(y) \mid y\in Y\}$ for 
any abelian group $G$,
\item $\dim_GX\leq \dim Y+\max\{\dim_Gf^{-1}(y) \mid y\in Y\}$ for 
any abelian group $G$,
\item $\dim_GX\leq \dim_GY+\max\{\dim_Gf^{-1}(y) \mid y\in Y\}$ if
$G$ is a principal ideal domain with the unity,
\item $\dim_GX\leq \dim_GY+\max\{\dim_Gf^{-1}(y) \mid y\in Y\}+1$ 
for any abelian group $G$.
\end{enumerate}
\end{theorem}

\begin{proof}
Let $n>\dim X,\dim Y$.

(1). 
We consider a map $f\times \id \colon X\times T_n(G)\to Y\times T_n(G)$. 
The Hurewicz Theorem from the Dimension Theory implies that 
\begin{align*}
\dim(X\times T_n(G))&
\leq \dim(Y\times T_n(G))\\
&
\quad
+ \max\{\dim(f\times \id)^{-1}(y,t) \mid (y,t)\in Y\times T_n(G)\}. 
\end{align*}
since 
$(f\times \id)^{-1}(y,t) = f^{-1}(y)$ 
for all $t$, we have the following 
$\dim_GX+n \leq \dim_GY+n+\max\{\dim f^{-1}(y) \mid y\in Y\}$.

(2). 
For that case we consider a map $f\circ\pi\colon X\times T_n(G)\to Y$, 
where $\pi\colon X\times T_n(G)\to X$ is the projection. 
By the Hurewicz theorem, we have 
$
\dim(X\times T_n(G))
\leq \dim Y
+ \max\{\dim(f\circ\pi)^{-1}(y) \mid y\in Y\}$. 
Then 
\begin{align*}
\dim_G X + n&
\leq \dim Y + \max\{\dim(f^{-1}(y)\times T_n(G)) \mid y\in Y\}\\
&
= \dim Y + \max\{\dim(f^{-1}(y))\} + n.
\end{align*}

(3). 
Let $G\in\sigma$. 
We consider a map 
$(f\circ\pi\times \id)\colon X\times T_n(G)\times T_n(G)\to Y\times T_n(G)$. 
Note that $(f\circ\pi\times \id)^{-1}(y,t)=f^{-1}(y)\times T_n(G)$. 
By Lemma 6.7 $T_n(G)\times T_n(G)$ is a $G$-testing space. 
This together with the Hurewicz theorem gives
\begin{equation*}\tag{$\ast$}
\dim_G X + \dim(T_n(G)\times T_n(G))
\leq \dim_G Y + n + \max\{\dim_Gf^{-1}(y)\} + n. 
\end{equation*}
If $G\neq\Z_{p^{\infty}}$, then $\dim(T_n(G)\times T_n(G))\leq 2n$ and 
hence, $\dim_GX\leq \dim_GY+\max\{\dim_Gf^{-1}(y)\}$. 
Let $G$ be a PID with the unity. 
Then by Proposition 6.8 no $\Z_{p^{\infty}}$ belongs to $\sigma(G)$. 
By the Bockstein theorem $\dim_GX=\dim_HX$ for some $H\in\sigma(G)$. 
Then 
\begin{align*}
\dim_G X&
= \dim_H X\\
&
\leq \dim_H Y + \max\{\dim_H f^{-1} \mid y\in Y\}\\
&
\leq \dim_G Y + \max\{\dim_G f^{-1}(y) \mid y\in Y\}.
\end{align*}

(4). Apply Proposition 3.4(2) to ($\ast$) to obtain the required
inequality.
\end{proof}

We recall that a map $f\colon X\to \R^n$ of a subset $X\subset\R^n$ is
called an $\epsilon$-move if $\|f(x)-x\|<\epsilon$ for all $x\in X$.

\begin{theorem}
For any compactum $X\subset\R^n$ and for any abelian group $G$ there is 
an $\epsilon>0$, such that the inequality $\dim_Gf(X)\geq \dim_GX$ holds
for every $\epsilon$-move $f\colon X\to\R^n$. 
\end{theorem}

\begin{proof}
Since a test space $T_n(G)$ is $n$-dimensional, we may assume that 
$T_n(G)\subset\R^{2n+1}$. 
The Alexandroff Theorem says that for a compactum 
$X\times T_n(G)\subset\R^n\times\R^{2n+1}$ there is a positive $\epsilon$ 
such that for every $\epsilon$-move $g\colon X\times T_n(G)\to\R^{3n+1}$ 
one has the inequality $\dim(g(X\times T_n(G)))\geq \dim(X\times T_n(G))$. 
Given an $\epsilon$-move $f\colon X\to\R^n$ we define another 
$\epsilon$-move $g\colon X\times T_n(G)\to\R^n\times\R^{2n+1}$ as 
$f\times \id$. 
We note that $g(X\times T_n(G))=f(X)\times T_n(G)$. 
Then 
\begin{align*}
\dim(g(X\times T_n(G)))&
= \dim(f(X)\times T_n(G))\\
&
= \dim_G X + n\\
&
\geq \dim(X\times T_n(G))\\
&
= \dim_G X + n.
\end{align*}
Hence, $\dim_Gf(X) \geq \dim_GX$.
\end{proof}

\section{Infinite-dimensional compacta of finite cohomological dimension}

According to the Realization Theorem (5.1) for any abelian group 
$G\in\sigma$ for any number $n\in\N$ there is an $n$-dimensional compactum 
$X_{n,G}$ with the cohomological dimension $\dim_GX_{n,G}=1$. Using this 
data it is easy to construct an infinite dimensional compactum $X$ with 
$\dim_GX=1$. It suffices to consider an one-point compactification 
$a(\bigcup_{n=1}^{\infty}X_{n,G})$ of a disjoint union of compacta 
$X_{n,G}$.

As it follows from 1.3(3) there is no such compactum for $G=\Z$. By the 
Alexandroff Theorem any $n$-dimensional compactum has $\dim_{\Z}X=n$. 
Nevertheless one can prove the following:

\begin{theorem}
There is an infinite-dimensional compactum $X$ having $\dim_{\Z}X\leq 3$.
\end{theorem}

\begin{proof}
The proof is based on the following result in K-theory: 
$$\tilde K_{\C}^*(K(\Z,n);\Z_p)=0
\quad \text{for $n\geq 3$ \cite{B-M},\cite{A-H}.}$$
Here $h^*=K_{\C}^*(\ ;\Z_p)$ is the reduced complex K-theory with $\Z_p$
coefficients, i.e.\ $h^*$ is generalized cohomology theory defined by the 
spectrum $E_{2n}=BU^{M(\Z_p,1)}$ and $E_{2n+1}=U^{M(\Z_p,1)}$. 
This cohomology theory is continuous, since $h^k(L)$ is a finite group for 
every compact polyhedron $L$. We apply Theorem 5.2 to $P=S^4$, $K=K(\Z,3)$ 
and $h^*$ for $n=0$ to obtain an essential map $f\colon X\to S^4$ of a
compactum $X$ having $\dim_{\Z}X\leq 3$. 
If we assume for a moment that the dimension of $X$ is finite, then by 
Alexandroff Theorem, $\dim X\leq 3$. 
But a map of 3-dimensional compactum to a 4-dimensional sphere cannot be 
essential. Hence $\dim X=\infty$.
\end{proof}

We note that one can use K-homology instead of K-cohomology here, since
$$\tilde K^{\C}_*(K(\Z,3);\Z_p)=0$$ 
as well. 
In that case a compactum $X$ has a $h_*$-essential map $f\colon X\to S^4$. 
Moreover by the proof of Theorem 5.2 one can assume that any given element 
$a\in h_*(P)$ lies in the image $\Im(f_*)$.
For applications we need a relative version of this.

\begin{theorem}
Let $h_*$ be a reduced generalized homology theory with 
$$h_*(K(G,n))=0,\quad n\in\N.$$ 
Then for every compact polyhedral pair $(P,L)$ and any element 
$a\in h_*(K,L)$ there is a compactum $X\supset L$ and a map $f\colon 
(X,L)\to (P,L)$ such that
\begin{enumerate}
\item $\dim_G(X\setminus L)\leq n$,
\item $a\in \Im(f_*)$ and
\item $f\mymid_L = \id_L$.
\end{enumerate}
\end{theorem}

This theorem is a relative version of Theorem 5.2 for $K=K(G,n)$.

If one applies this theorem to the pair $(B^4,\partial B^4)$ with $G=\Z$, 
$n=3$ and $h_*=\tilde K_*(\ ;\Z_p)$ for odd $p$ and some nontrivial 
element in $h_*(B^4,\partial B^4)$, he gets a compactum $X\supset S^3$ of 
$\dim_{\Z}X=3$ and an essential map onto $B^4$. Hence, $X$ is 
infinite-dimensional as in Theorem 7.1.

Using more advanced algebraic topology we are going to prove the following:

\begin{theorem}[\cite{D-W2}]
There is an infinite dimensional compactum $X$ with 
$\dim_{\Z}(X\times X)=3$.
\end{theorem}

We recall that a truncated spectrum is a sequence of pointed spaces 
$\E=\{E_i\}$, $i\leq 0$, such that $E_{i-1}=\Omega E_i$. Thus, any 
truncated spectrum is generated by one space $E_0$.
The lower half of every $\Omega$-spectrum is an example of a truncated 
spectrum.
The reduced truncated cohomology of a given space $X$ with coefficients 
in a given truncated spectrum $T^i(X;\E)$ is the set of pointed homotopy 
classes of mappings $X$ to $E_i$.
Note that $T^i(X)$ is a group for $i<0$ and it is an abelian group for 
$i<-1$. 
Truncated cohomologies possess many features of generalized cohomology. 
For every map $f\colon X\to Y$ there is the induced homomorphism ($i>0$) 
$f^*\colon T^i(Y)\to T^i(X)$. 
Homotopic maps induce the same homomorphism and a null-homotopic map 
induces zero homomorphism.
There is the natural Mayer-Vietoris exact sequence
$$
\cdots
\to T^r(A\cup B)
\to T^r(A)\times T^r(B)
\to T^r(A\cap B)
\to T^{r+1}(A\cup B)
\to \cdots
$$
of groups for $r\leq -1$ and abelian groups for $r\leq -2$. 
Therefore Lemma 5.9 holds for a truncated cohomology for $n\leq -2$. 
We call a truncated homology $T^*$ continuous if for every direct limit 
of finite CW-complexes $L=\limto \{L_i;\lambda^i_{i+1}\}$ the following 
formula holds $T^k(L)=\Limleftarrow T^k(L_i)$ for $k<0$. 
We note that the Milnor Theorem holds for truncated cohomologies: 
$$
0
\to \Limone \{T^{k-1}(L_i)\}
\to T^k(L)
\to \Limleftarrow \{T^k(L_i)\}
\to 0.
$$

Hence, if $T^k(M)$ is a finite group for every finite complex $M$ and 
every $k<-1$, By the Mittag-Lefler condition $T^*$ is continuous.

We consider a truncated cohomology $T^*$ generated by a mapping space 
$E_0=(S^7)^M$ where $M=M(\Z_2,1)=\R P^2$ is a Moore space of the type 
$(\Z_2,1)$ and $S^7$ is the 7-dimensional sphere.

\begin{lemma}
The truncated cohomology theory $T^*$ is continuous.
\end{lemma}

For the proof we need the following

\begin{proposition}
Let $\nu_2\colon S^1\to S^1$ be a map of the degree two.
Then the map $\nu_2\wedge \id\colon S^1\wedge\R P^2\to S^1\wedge\R P^2$ is 
null homotopic.
\end{proposition}

\begin{proof}
The space $S^1\wedge\R P^2$ is the suspension $\Sigma M$ over the 
projective space and it can be defined as a quotient map 
$p\colon B^3\to \Sigma M$.
Temporarily we denote by 2 a fixed map of degree 2 between 2-spheres and 
by 1, the identity map of the 2-sphere. Let $C_q$ denote the mapping cone 
of a map $q\colon X\to Y$ i.e.\ $C_q=\cone(X) \cup_q Y$. Consider the 
following commutative diagram:
$$
\begin{CD}
S^2	@>>1>	S^2	@>>>	C_1\\
@VV1V		@VV2V		@VVpV\\
S^2	@>>2>	S^2	@>>>	C_2\\
@AA2A 		@AA2A 		@AAgA\\
S^2	@>>2>	S^2	@>>>	C_2\\
\end{CD}
$$
Here the mapping cone $C_1$ is homeomorphic to a 3-ball $B^3$ and $C_2$ 
is homeomorphic to $\Sigma M$.
First we note that the map $g$ is homotopic to the map $\nu_2\wedge \id$.
Then we show that $g$ has a lift $g'\colon \Sigma M\to B^3$ with respect 
to $p$. 
In fact $g'$ is defined by the following diagram:
$$
\begin{CD}
S^2	@>>1>	S^2	@>>>	C_1\\
@AA2A		@AA1A		@AAg'A\\
S^2	@>>2>	S^2	@>>>	C_2\\
\end{CD}
$$
Since $B^3$ is contractible, $g'$ is null-homotopic and, hence $g$ is 
null-homotopic.
\end{proof}

\begin{proof}[Proof of Lemma 7.4]
Show that every element of a group $T^k(L)$ has an order 2 for $k<0$.
Indeed, 
$T^k(L)
= [L,\Omega^{-k}(S^7)^M]
= [\Sigma M,(S^7)^{\Sigma^{-k-1}L}]$.
For any space $N$ and for any element $a\in[\Sigma M,N]$ represented by a 
map $f\colon \Sigma M\to N$, the element $2a$ is represented by a map 
$f\circ(\nu_2\wedge \id)$ and it is homotopic to zero by virtue of 
Proposition 7.5. 
Note that $T^k(L)=[S^k\wedge L\wedge M,S^7]$. 
When a complex $L$ is finite this group is finitely generated. 
Hence in the case of $k<-1$, the group $T^k(L)$ of any finite complex $L$ 
is finite. 
As we know it suffices for the continuity.
\end{proof}

\begin{proposition}
For every $k<0$ we have $T^k(K(\Z[\frac{1}{2}],1))=0$.
\end{proposition}

\begin{proof}
We can present $K(\Z[\frac{1}{2}],1)$ as the direct limit of complexes 
$M_i$ where each $M_i$ is homotopy equivalent to the circle $S^1$ and 
every bonding map $\xi_i\colon M_i\to M_{i+1}$ is homotopy equivalent to a 
map of the degree two $S^1\to S^1$. Then
$T^k(K(\Z[\frac{1}{2}],1))
= [\limto \{M_i,\xi_i\},(S^7)^M]
= [(\limto \{M_i,\xi_i\})\wedge M,S^7]
= [\limto \{M_i\wedge M,\xi_i\wedge \id\},S^7]$. 
Consider a bonding map 
$\xi_i\wedge \id\colon M_i\wedge M\to M_{i+1}\wedge M$. 
This map is homotopy equivalent to the map $\nu_2\wedge \id$ and hence, it 
is homotopically trivial. 
Therefore the space 
$\limto \{M_i\wedge M,\xi_i\wedge \id\}$ is 
homotopically trivial. 
Hence, $T^k(K(\Z[\frac{1}{2}],1))=0$.
\end{proof}

We also need the following result.

\begin{theorem}[Miller Theorem (Sullivan Conjecture) \cite{Mi}]
Let $K$ be a CW-complex of finite dimension and $\pi$ be a finite group. 
Then the mapping space $K^{K(\pi,1)}$ is weakly homotopy equivalent to a 
point.
\end{theorem}

\begin{proposition}
For every $k$ we have $T^k(K(\Z_2,1))=0$.
\end{proposition}

\begin{proof}
We note that 
$T^k(K(\Z_2,1))
= [K(\Z_2,1),(S^7)^{\Sigma^kM}]
= [\Sigma^kM,(S^7)^{K(\Z_2,1)}]
= 0$
by Theorem 7.7.
\end{proof}

The following Proposition is a version of Theorem 5.2 for a truncated 
cohomology.

\begin{proposition}
Let $P$ and $K$ be simplicial complexes and assume that $K$ is countable 
complex. 
Let $T^*$ be a reduced truncated continuous cohomology theory.
If $T^n(P)\neq 0$ and $T^k(K)=0$ for some $n<-1$ and all $k<n$, then 
there exist a compactum $X$, having the property $K\in AE(X)$, and a 
$T^n$-essential map $f\colon X\to P$.
\end{proposition}

The proof is the same. 

\begin{proof}[Proof of Theorem 7.3]
We take $P=S^3$, 
$K = K([\Z[\frac{1}{2}],1)\vee K(\Z_2,1)$ and $K^*$ is as above.
Note that 
$T^{-2}(S^3)
= [S^3,\Omega^2(S^7)^M]=[S^3\wedge S^2\wedge M, S^7]
= [\Sigma^5M,S^7]
= [M(\Z_2,6),S^7]
= H^7(M(\Z_2,6))
= H_6(M(\Z_2,6))
= \Z_2
\neq 0$.
By Propositions 7.5 and 7.6 we have $T^k(K)=0$ for $k\leq-2$. 
Proposition 7.9 gives us a compactum $X$ with $\dim_{\Z_2}X\leq 1$ and 
$\dim_{\Z[\frac{1}{2}]}X\leq 1$ and an essential map $f\colon X\to S^3$. 
By the Bockstein Theorem $\dim_{\Z_{(q)}}X\leq 1$ for all prime $q\neq p$. 
Hence, the cohomological dimensions of $X$ with respect to all fields 
from the Bockstein basis $\sigma$ do not exceed one. 
Hence by Theorem 3.15 $\dim_{\Z}(X\times X)\leq 3$. 
Hence (see also Lemma 2.9), $\dim_{\Z}X\leq 2$.
Since $X$ admits an essential map onto $S^3$, the dimension of $X$ cannot 
be $\leq 2$.
Therefore by Alexandroff Theorem $\dim X=\infty$.
\end{proof}

We recall that a space $X$ is strongly infinite dimensional provided that 
there exists an essential map $f\colon X\to I^{\infty}$ of $X$ onto the 
Hilbert cube. 
A map $f\colon X\to I^{\infty}$ is essential provided that $p\circ f\colon 
X\to I^n$ is essential for each coordinate projection 
$p\colon I^{\infty}\to I^n$. 
It is known that this definition does not depend on the product structure 
on the Hilbert cube $I^{\infty}$. 
Finally we recall that a map $f\colon X\to I^n$ is essential provided the 
extension problem $(f^{-1}(\partial I^n),f\mymid_{\cdots})$ on $X$ for 
mappings to $\partial I^n$ has no solution.

\begin{theorem}
Let $h^*$ be a reduced continuous cohomology theory such that $h^*(K)=0$ 
for some countable simplicial complex $K$. Then there exists a strongly 
infinite dimensional compactum $X$ having the property $K\in AE(X)$.
\end{theorem}

\begin{corollary}
There exists a strongly infinite dimensional compactum $X$ with 
$\dim_{\Z}X\leq 3$.
\end{corollary}

\begin{proof}
Take $K=K(\Z,3)$ and $h^*=\tilde K_{\C}^*(\ ;\Z_p)$.
\end{proof}

\begin{corollary}
For every prime $p$ there is a strongly infinite dimensional compactum 
$X$ with $\dim_{\Z[\frac{1}{p}]}X=1$.
\end{corollary}

\begin{proof}
Take $K=K(\Z[\frac{1}{p}],1)$ and $h^*=\tilde H^*(\ ;\Z_p)$.
\end{proof}

\begin{proof}[Proof of Theorem 7.10]
By induction we construct two inverse sequences $\{P_k,q^{k+1}_k\}$ and
$\{I^k,\omega^{k+1}_k\}$ and a morphism between them, i.e.\ a sequence of 
maps $\{f_k\colon P_k\to I^k\}$ such that all squares are commutative. 
The first sequence consists of polyhedra and the second sequence consists 
of $k$-cubes, $k=1,2,\dots$ with bonding maps 
$\omega^{k+1}_k\colon I^{k+1}\to I^k$ defined as projections on factors.
For every $k$ we define by the same induction an element 
$\mu_k\in h^*(I^k\partial I^k)$ and a countable basis $\sA^k$ of extension 
problems on $P_k$ with respect to the complex $K$ and consisting of 
simplicial problems. 
We construct the sequences in such a way that
\begin{enumerate}
\item
$(q^n_k)*(f^*_k(\mu_k))\neq 0$ for any $k$ and every $n>k$,
\item
every extension problem $(A^k_i,\alpha^k_i)\in\sA^k$ is resolved by 
$q^j_k$ for some $j$.
\end{enumerate}

First, assume that we can construct such sequences. 
Then by Proposition 5.7 the limit space 
$X=\limleftarrow \{P_k,q^{k+1}_k\}$ has the property $K\in AE(X)$.
Since 
$(\omega^{\infty}_k\circ f)*(\mu_k)
= (q^{\infty}_k)^*(f_k)*(\mu_k)
\neq 0$,
the map $\omega^{\infty}_k\circ f\colon X\to I^k$ is essential for every 
$k$. 
Therefore the limit map $f\colon X\to I^{\infty}$ is essential and, hence 
$X$ is strongly infinite dimensional.

Now we present the induction. We define $P_1=I^1$ and $f_1=\id$. Take 
nonzero element $\mu_1\in h^*(I,\partial I)$ and fix a basis $\sA^1$.
Assume that the commutative diagram
$$
\begin{CD}
P_1	@<q^2_1<<	P_2	@<q^3_2<<	\cdots	@<q^k_{k-1}<<		P_k\\
@Vf_1VV			@Vf_2VV			@.				@Vf_kVV\\
I^1	@<\omega^2_1<<	I^2	@<\omega^3_2<<	\cdots	@<\omega^k_{k-1}<<	I^k\\
\end{CD}
$$
is already constructed, elements $\mu_i\in h^*(I^i,\partial I^i)$ are 
defined for $i\leq k$ and extension problem bases $\sA^i$, $i\leq k$, are 
fixed such that
\begin{enumerate} 
\item
$(q^k_i)*(f^*_i(\mu_i)\neq 0$ for all $i\leq k$,
\item
All problems $\bigcup_{i=1}^k\sA^i$ are enumerated by all numbers of 
the form $p_1^{l_1}p_2^{l_2}\dots p_k^{l_k}$ where $p_1,p_2,\dots, p_k$ 
are the first $k$ prime numbers,
\item
For every $i< k$ the problem having $i$-th number is resolved by some map 
$q^k_j$.
\end{enumerate}

To make an induction step we note that the number $k$ has the form 
$p_1^{l_1}\dots p_k^{l_k}$.
Hence there is an extension problem $(A^r_i,\alpha^r_i)\in\sA^r$ for 
$r\leq k$ having the number $k$ in our list. 
We lift that problem to the $k$-th level and apply Lemma 5.6 to resolve 
that lift by a simplicial (with respect to some subdivisions) map 
$g\colon L\to P_k$ having point preimages homeomorphic to $K$ or to 
one-point space. 
By virtue of Lemma 5.9, the induced homomorphism 
$g^*\colon H^*(P_k)\to h^*(L)$ is an isomorphism. 
Since $h^*$ is continuous, there exists a compact subcomplex $L'\subset L$ 
such that $g_1^*((f_i\circ q^k_i)^*(\mu_i))\neq 0$ for $i\leq k$ where 
$g_1$ is the restriction of $g$ onto $L'$. 
We define the complex $P_{k+1}=L'\times I$ and the bonding map 
$q^{k+1}_k\colon P_{k+1}\to P_k$ as the composition $g_1\circ\omega$ 
where $\omega\colon L'\times I\to L'$ is the projection. 
We define $f_{k+1}\colon P_{k+1}\to I^{k+1}$ as the product
$(f_k\circ g_1)\times \id\colon L'\times I\to I^k\times I$. 
We let $\mu_{k+1}$ to be the suspension $\Sigma\mu_k$. Then we define a 
countable basis $\sA^{k+1}$ of extension problems on $P_{k+1}$ consisting 
of simplicial problems. 
Enumerate all the problem in the list $\sA^{k+1}$ by all numbers of the 
form $p_1^{l_1}\dots p_k^{l_k}p_{k+1}^{l_{k+1}}$ with $l_{k+1}>0$.
 
Let us verify the properties (1)--(3) for $k+1$. 
It is clear that the conditions (2)--(3) hold.
By the construction the property (1) holds for $i<k$. 
Then all we need is to check that $f_{k+1}^*(\mu_{k+1})\neq 0$. 
We note that the homomorphism 
$f_{k+1}^*\colon h^*(I^{k+1},\partial I^{k+1})
\to h^*(P_{k+1},(f_{k+1})^{-1}(\partial I^k))$
is generated by the following map
$\id_{S^1} \wedge (f_k\circ g_1) \colon 
S^1 \wedge (L'/(f_k\circ g_1)^{-1}(\partial I^k)
\to S^1 \wedge (I^k/\partial I^k)$ 
which is the suspension $\Sigma(f_k\circ g_1)$. 
Since $(f_k\circ g_1)^*(\mu_k)\neq 0$, we have that 
$f^*_{k+1}(\mu_{k+1})\neq 0$.

Thus, the induction step is completed.
\end{proof}

\section{Resolution theorems}

In this section we are proving some resolution theorems for the 
cohomological dimension. We start from resolving of polyhedra. 
First we describe Williams' construction.

\begin{definition}
A simplicial complex over $n$-simplex $\Delta^n$ is a pair $(L,\xi)$ 
where $L$ is a simplicial complex and $\xi\colon L\to\Delta^n$ is 
nondegenerate simplicial map (no edge goes to a vertex).
\end{definition}

\begin{example*}
The first barycentric subdivision of any simplicial $n$-dimensional 
complex $K$ defines the natural complex over $\Delta^n$. The map
$\xi\colon \beta^1K\to\Delta^n=\{0,1,\dots,n\}$ assigns to every 
barycenter $c_{\sigma}\in\beta^1K$ the dimension of corresponding 
simplex $\sigma$. 
\end{example*}

Now for every resolution $f\colon X\to \Delta^n$ of a simplex $\Delta^n$ 
we can define a resolution of a simplicial complex $(L,\xi)$ over 
$\Delta^n$ by taking the pullback:
$$
\begin{CD}
X\Delta L	@>>\xi'>	X\\
@VVf'V				@VVfV\\
|L|		@>>\xi>		\Delta^n\\
\end{CD}
$$
For example Pontryagin surfaces (Example 1.9) were constructed by taking 
resolutions of some triangulations of 2-dimensional polyhedra which are 
induced by a resolution $\xi\colon M_p\to\Delta^2$. Recall that $\xi$ is 
a simplicial map of $M_p$ onto a 2-simplex $\Delta^2$.
Here $M_p$ is the mapping cylinder of a map of degree $p$ between two 
circles.

\begin{definition}
Let $G$ be an abelian group and $L$ be a simplicial complex.
An {\it Edwards-Walsh resolution} of $L$ in the dimension $n$ is a pair
$(EW(L,G,n),\omega)$ consisting of a CW-complex $EW(L,G,n)$ and a map
$\omega\colon EW(L,G,n)\to |L|$ onto a geometric realization of 
$L$ 
such that
\begin{enumerate}
\item 
$\omega$ is 1-to-1 over the $n$-skeleton $L^{(n)}$, hence it defines an 
inclusion $j\colon L^{(n)}\subset EW(L,G,n)$,
\item 
for every simplex $\Delta$ of $L$, $\omega^{-1}(\Delta)$ is a subcomplex
of $EW(L,G,n)$ having the type of Eilenberg-MacLane space $K(\bigoplus 
G,n)$,
\item 
for every simplex $\Delta$ of $L$ the inclusion 
$\omega^{-1}(\partial\Delta)\subset\omega^{-1}(\Delta)$ induces an 
epimorphism $H^n(\omega^{-1}(\Delta);G)\to H^n(\omega^{-1}(\partial\Delta);G)$.
\end{enumerate}
Here we regard a contractible space as $K(\bigoplus G,n)$ with zero 
number of summands $G$.
We recall that $\Z_{(\sL)}$ denotes the localization of integers at set 
of primes $\sL\subset\sP$.

We say that an abelian group $G$ is {\it $\sL$-local modulo torsion} if
$G/\Tor(G)=G\otimes\Z_{(\sL)}$.
\end{definition}

\begin{lemma}
For any of the groups $\Z$, $\Z_{(\sL)}$, $\Z_p$ for any $n\in\N$ and for 
any simplicial complex $L$ over a simplex $\Delta^m$ there is an 
Edwards-Walsh resolution $\omega\colon EW(L,G,n)\to |L|$ with the 
additional property for $n>1$:
\begin{itemize}
\item[(4-$\Z$)]
the $(n+1)$-skeleton of $EW(L,\Z,n)$ is isomorphic to $L^{(n)}$,
\item[(4-$\Z_p$)]
the $(n+1)$-skeleton of $EW(L,\Z_p,n)$ is obtained from $L^{(n)}$ by
attaching $(n+1)$-cells by a map of degree $p$ to the boundary 
$\partial\Delta^{n+1}$ for every $(n+1)$-dimensional simplex 
$\Delta^{n+1}$.
\item[(4-$\Z_{(\sL)}$)]
for every subcomplex $N\subset L$ the homomorphism 
$j_*\colon H_n(N^{(n)};\Z_{(\sL)})\to H_n(\omega^{-1}(N);\Z_{(\sL)})$ 
generated by the inclusion of the $n$-skeleton of $N$ in~$\omega^{-1}(N)$ 
is an isomorphism and the kernel of the homomorphism
$\omega_*\colon H_n(\omega^{-1}(N))\to H_n(N)$ is $\sL$-local modulo 
torsions.
\end{itemize}
\end{lemma}

\begin{proof}
First we consider the case when $n>1$. We consider three different cases.

($\Z$). 
Induction on $m$. If $m\leq n$, we define $EW(L,\Z,n) = |L|$ and
$\omega=\id_L$. Assume that there is a resolution with the properties 
(1)--(4) for $m$-dimensional complex $L$. Consider a simplex 
$\Delta^{m+1}$ of the dimension $m+1$. 
The barycentric subdivision of its boundary 
$K=\beta^1\partial\Delta^{m+1}$ is a complex over $\Delta^m$ and, hence, 
we can apply the induction assumption. The $n$-dimensional homotopy group
$\pi_n(EW(K,\Z,n))$ is equal to $\pi_n(K^{(n)})$ by the property 
(4-$\Z$). Since $K^{(n)}$ is homotopy equivalent to the wedge of 
$n$-spheres, the $n$-th homotopy group equals $\bigoplus \Z$. Therefore 
there exists a complex $\bar K\supset EW(K,\Z,n)$ of the type 
$K(\bigoplus\Z,n)$ such that its $n{+}1$-skeleton coincides with the 
$n{+}1$-skeleton of $EW(K,\Z,n)$ and, hence, coincides with $K^{(n)}$. We 
define a map $\bar\omega\colon \bar K\to\Delta^{m+1}$ such that
$\bar\omega$ is an extension of $\omega\colon EW(K,\Z,n)\to K$ and 
$\omega^{-1}(t)=\bar\omega^{-1}(t)$ for every $t\in K$. We note that 
$\bar\omega$ has all properties 1-4 and, hence is a resolution of
$m{+}1$-simplex. Then we apply the Williams construction to obtain a 
resolution of an arbitrary complex over $\Delta^{m+1}$. All properties 
(1)--(4) are easy to verify.

($\Z_p$).
The same, induction on $m$. Now we apply the property (4-$\Z_p$) to 
compute the $n$-dimensional homotopy group $\pi_n(EW(K,\Z_p,n))$. 
The result is $\bigoplus\Z_p$. Similarly we construct $\bar K$ from 
$EW(K,\Z_p,n)$ by attaching cells in dimensions $n+2$ and higher. 
Then we apply the Williams construction.

($\Z_{(\sL)}$).
We apply induction on $m$. 
If $m\leq n$, we define $EW(L,\Z_{(\sL)},n)= |L|$ and 
$\omega=\id_L$. 
Let $m\ge n+1$ and let $K$ be as in ($\Z$). 
If $m=n+1$ we attach to the $n$-sphere $K$ a complex $K(\Z_{(\sL)},n)$ 
having a $(\sP\setminus\sL)$-telescope as the $(n{+}1)$-skeleton to obtain 
a complex $\bar K$.
If $m>n+1$, then $K$ is $n$-connected and hence the condition 
(4-$\Z_{(\sL)}$) and the induction assumption imply that the group 
$H_n(\omega^{-1}(K))$ is $\sL$-local modulo torsion. 
Hence there is a short exact sequence
$$0\to G\to H_n(\omega^{-1}(K))\to H_n(\omega^{-1}(K))\otimes\Z_{(\sL)}\to 0$$
where $G$ is $(\sP\setminus\sL)$-torsion group. 
We note that 
$$H_n(\omega^{-1}(K))\otimes\Z_{(\sL)}
= H_n(\omega^{-1}(K);\Z_{(\sL)})
= H_n(K^{(n)};\Z_{(\sL)})
= \bigoplus\Z_{\sL}$$ 
by the induction assumption.
Since $\omega^{-1}(K)$ is $(n{-}1)$-connected, then by the Hurewicz 
theorem 
we have a short exact sequence 
$0\to G\to\pi_n(\omega^{-1}(K))\to\bigoplus\Z_{(\sL)}\to 0$.
We attach $(n{+}1)$-cells to $\omega^{-1}(K)$ along generators of the 
group 
$G$ and then we attach cells of higher dimension to obtain a complex 
$\bar K$ of the type $K(\bigoplus\Z_{(\sL)},n)$. 

As above we define an extension $\bar\omega\colon \bar K\to\Delta^{m+1}$ 
of $\omega\colon \omega^{-1}(K)\to K$ such that 
$\bar\omega(\bar K\setminus\omega^{-1}(K))\subset \Int\Delta^{m+1}$.
Then we apply Williams' construction to obtain a resolution of an 
arbitrary complex over $\Delta^{m+1}$. The conditions (1)--(2) of an 
EW-resolution hold automatically.

To verify (3) we show that every map 
$f\colon \omega^{-1}(K)\to K(\Z_{(\sL)},n)$ admits an
extension $\bar f\colon \bar K\to K(\Z_{(\sL)},n)$. 
It holds true when $m=n+1$. 
For $m>n+1$ we note that $f_*(G)=0$ where the homomorphism 
$f_*\colon \pi_n(\omega^{-1}(K))\to\pi_n(K(\Z_{(\sL)},n))$ 
is induced by $f$.
It means that $f$ can be extended to the $(n{+}1)$-dimensional skeleton of 
$\bar K$.
Since there is no obstruction for extending this map over higher 
dimensional cells, the required extension exists.

Now we check the property (4-$\Z_{(\sL)}$) by induction on the number of 
$(m{+}1)$-simplices in $N$. If that number is zero, the condition holds 
by the induction assumption. Let $N=N_1\cup\Delta$ where $\Delta$ is an 
$(m{+}1)$-simplex with $\Delta\cap N_1=\partial\Delta$. We consider the 
diagram generated by the inclusion $j\colon N^{(n)}\subset\omega^{-1}(N)$ 
and the homology Mayer-Vietoris sequence with $\Z_{(\sL)}$-coefficients 
for the triples $(N^{(n)},N^{(n)}_1,\Delta^{(n)})$ and 
$(\omega^{-1}(N^{(n)}), \omega^{-1}(N^{(n)}_1),\omega^{-1}(\Delta^{(n)}))$.
We note that the spaces $\partial\Delta^{(n)}$ and 
$\omega^{-1}(\partial\Delta^{(n)})$ are $(n{-}1)$-connected.
Then the induction assumption and the five lemma imply that $j_*$ is an 
isomorphism.
Since the homomorphism $H_n(N^{(n)})\to H_n(N)$ is an epimorphism which 
is factored through the homomorphism 
$\omega_*\colon H_n(\omega^{-1}(N))\to H_n(N)$, 
the latter is also an epimorphism. 
We apply the homomorphism generated by tensoring with $\Z_{\sL}$ to the 
short exact sequence $0\to K_N\to H_n(\omega^{-1}(N))\to H_n(N)\to 0$,
where $K_N$ is the corresponding kernel. We expand that diagram by taking 
the Mayer-Vietoris sequence for the triad $(N,N_1,\Delta)$ and its 
preimage $\omega^{-1}(N,N_1,\Delta)$.
For the kernels we obtain the diagram:
$$
\begin{CD}
K_{N_1}\oplus K_{\Delta} 
@>>>
K_N 
@>>> 
0\\
@VVV @V\phi VV\\
K_{N_1}\otimes\Z_{\sL}\oplus K_{\Delta}\otimes\Z_{\sL} 
@>>> 
K_N\otimes\Z_{\sL} 
@>>> 
0\\
\end{CD}
$$
Then the induction assumption implies that the homomorphism $\phi$ is an 
epimorphism. Since $H_n(N^{(n)})$ is torsion free, by already proven part 
of the property (4-$\Z_{(\sL)}$) we have that all torsions of the group 
$H_n(\omega^{-1}(N))$ are $(\sP\setminus\sL)$-torsions. Therefore all 
torsions of $K_N$ are of that type. Then we can conclude that the group 
$K_N$ is $\sL$-local modulo torsions. 

By an abelinization of a finite complex $L$ we understand a finite 
complex $\ab(L)$ obtained from $L$ by attaching 2-dimensional cells 
killing all nontrivial commutators of a finite set of generators of the 
fundamental group $\pi_1(L)$. 
If $L^{(1)}$ is the 1-dimensional skeleton of a simplicial complex $L$, 
then by $\ab_L(L^{(1)})$ we denote an inductively constructed complex
$\ab_L(L^{(1)})
= L_{\dim L-1} \supset \dots \supset L_3 \supset L_2 \supset L_1
= L^{(1)}$. 
Where $L_2$ is the union of abelinizations $\ab(\sigma^{(1)})$ of 
1-skeletons of 3-simpexes $\sigma\in L$. To construct $L_3$ we consider 
a 4-dimensional simplex $\delta\in L$ and consider 
$\ab_{\delta^{(3)}}\delta^{(1)}\subset L_2$ and take its
abelinization. Do it for all 4-simplexes, then $L_3$ will be the union of 
all those abelinizations and so on.
\end{proof}

If $n=1$ the property (4) for the group $\Z$ takes the following form:
$$EW(L,\Z,1)^{[2]}=\ab_L(L^{(1)}).$$
Here by $Y^{[k]}$ we denote the k-dimensional skeleton of CW-complex $Y$.
For the group $\Z_p$ the property (4) becomes the following: 
$$EW(L,\Z_p,1)^{[2]}
= \ab_L(L^{(1)}) \cup_p \{B^2 \mid \sigma\in L, \dim\sigma=2\}.$$
For the group $\Z_{(\sL)}$ the property (4) remains the same. 
Then the argument is basically the same as in the case $n>1$.

\begin{lemma}
Assume that a compact $X$ has the cohomological dimension $\dim_GX\leq n$. 
Then for every Edwards-Walsh resolution $\omega\colon EW(L,G,n)\to L$ and 
for every map $f\colon X\to L$ there is a map $f'\colon X\to EW(L,G,n)$ 
such that $\omega f'(x)$ lies in the same simplex of $L$ as $f(x)$ for 
every point $x\in X$.
\end{lemma}

\begin{proof}
The result follows from the property (2) of Edwards-Walsh resolution and 
the fact that $K(\bigoplus G,n)\in AE(X)$.
\end{proof}

Suppose that $\{X_i,p^{i+1}_i\}$ is an inverse sequence of pointed spaces 
and base point preserving bonding maps. Then for every $m$ there is a 
natural embedding of the product $X_1\times\dots\times X_m$ into the 
infinite product $\prod_{i=1}^{\infty}X_i$.
The sequence 
$$\begin{CD}
X_1	@<p^2_1<<	X_2	@<p^3_2<<	\cdots	@<P^m_{m-1}<<	X_m
\end{CD}$$ 
defines an imbedding of $X_m$ into the product 
$\prod_{i=1}^mX_i\subset\prod_{i=1}^{\infty}X_i$. The inverse sequence
$\{X_i,p^{i+1}_i\}$ defines an embedding of the limit space $X$ in 
$\prod_{i=1}^{\infty}X_i$. The projection in the inverse sequence 
$p^{\infty}_m\colon X\to X_m$ coincides with the restriction on $X$ of 
the projection onto the factor 
$\prod_{i=1}^{\infty}X_i\to\prod_{i=1}^mX_i$.
This system of imbeddings in $\prod_{i=1}^{\infty}X_i$ we call a {\it 
realization of the inverse sequence} $\{X_i,p^{i+1}_i\}$ in the product 
$\prod_{i=1}^{\infty}X_i$.

Let $\rho_i$ be a metric on $X_i$ and let $\delta_i$ denote the diameter 
of $X_i$.
We assume that $\sum_{i=1}^{\infty}\delta_i <\infty$. 
Then the formula 
$\rho(x,y)=\sum_{i=1}^{\infty}\rho_i(p^{\infty}_i(x),p^{\infty}_i(y))$ 
defines a metric $\rho$ on the product $\prod_{i=1}^{\infty}X_i$.

Let $\sM$ be a (finite) cover of a compact space $X$ with a given metric 
$\rho$.
By $d(\sM)$ we denote the diameter of $\sM$ = $\max\{\diam M \mid 
M\in\sM\}$ 
and by $\lambda(\sM)$ we denote the
Lebesgue number of $\sM$:
$$\lambda(\sM)
= \max\{r \mid \text{for any r-ball $O_r(x)$ there is $M\in\sM$, 
$O_r(x)\subset M$}\}.
$$
Here $O_r(x)$ is the ball in $X$ of a radius $r$ with respect to $\rho$, 
centered at $x\in X$. Let $M_x$ denote an arbitrary $M\in\sM$ with the 
property $x\in O_{\lambda(\sM}(x)
\subset \Cl(M)$.

\begin{lemma}
Let $X=\limleftarrow \{K_i,f^{i+1}_i\}$ and 
$Z=\limleftarrow \{L_i,g^{i+1}_i\}$
be limit spaces of inverse systems of compacta. Suppose the first 
sequence is realized in $\prod_{i=1}^{\infty}K_i$ and for every $i$ a 
finite cover $\sM^i$ of $K_i$, with the diameter $d_i$ and the Lebesgue 
number $\lambda_i$, and a mapping $\alpha_i\colon L_i\to K_i$ are defined 
such that
\begin{enumerate}
\item $\alpha_i(L_i)\cap M\neq\emptyset$ for every $M\in\sM^i$,
\item $d_i<\lambda_{i-1}/4$,
\item the diagram
$$
\begin{CD}
L_{i+1}	@>\alpha_{i+1}>>K_{i+1}\\
@Vg^{i+1}_iVV 		@Vf^{i+1}_iVV\\
L_i	@>\alpha_i>>	K_i\\
\end{CD}
$$
is $\lambda_i/4$-commutative.
\end{enumerate}
Then there exists a continuous map $\alpha\colon Z\to X$ onto $X$ such 
that the point preimage $\alpha^{-1}(x)$ is the limit space of 
$\limleftarrow \{\alpha^{-1}_i(M_{x_i}),q^{i+1}_i\}$ where 
$x_i=f^{\infty}_i(x)$ and $q^{i+1}_i$ is the restriction of $g^{i+1}_i$ 
on $\alpha_i^{-1}(M_{x_i})$.
\end{lemma}

\begin{proof}
(A). 
First we show that for any $i$ and $k$ the diagram
$$
\begin{CD}
L_{i+k}	@>\alpha_{i+k}>>K_{i+k}\\
@Vg^{i+k}_iVV		@Vf^{i+k}_iVV\\
L_i	@>\alpha_i>>	K_i\\
\end{CD}
$$
is $\lambda_i/2$-commutative, i.e.\ 
$\rho(\alpha_i(g^{i+k}_i(z)),f^{i+k}_i(\alpha_{i+k}(z)))
< \lambda_i/2$ for all $z\in L_{i+k}$.

We apply induction on $k$. For $k=1$ it follows by the condition (3) of 
the lemma. For $k>1$ we apply the triangle inequality to obtain 
\begin{multline*}
\rho(\alpha_i(g^{i+k}_i(z)),f^{i+k}_i(\alpha_{i+k}(z)))\\
\leq \rho(\alpha_i(g^{i+k}_i(z)), f^{i+1}_i\alpha_{i+1}g^{i+k}_{i+1}(z))
+ \rho(f^{i+1}_i\alpha_{i+1}g^{i+k}_{i+1}(z), f^{i+k}_i(\alpha_{i+k}(z))).
\end{multline*}
By the induction assumption we have 
$$
\rho(\alpha_i g^{i+1}_i(g^{i+k}_{i+1}(z)), f^{i+1}_i\alpha_{i+1}(g^{i+k}_{i+1}(z)))
< \lambda_{i+1}/4
$$ 
and
$$
\rho(\alpha_{i+1}g^{i+k}_{i+1}(z),f^{i+k}_{i+1} \alpha_{i+k}(z))
< \lambda_{i+1}/4.
$$ 
By the definition of the metric $\rho$ the map $f^{i+1}_i$ is a 
contraction, hence
$$
\rho(f^{i+1}_i\alpha_{i+1}g^{i+k}_{i+1}(z),f^{i+k}_i(\alpha_{i+k}(z)))
< \lambda_{i+1}/4.
$$
Therefore by the condition (2) of the Lemma we have the desired inequality
$$
\rho(\alpha_i(g^{i+k}_i(z)),f^{i+k}_i(\alpha_{i+k}(z)))
< \lambda_i/2
$$ 
for all $z\in L_{i+k}$.

(B).
Then we prove that the sequence of maps 
$\alpha_ig^{\infty}_i\colon Z\to\prod_{i=1}^{\infty}K_i$
has a limit.

Denote by $s_k$ the sum $\sum_{i=k}^{\infty}\delta_i$ where $\delta_i$ 
is the diameter of $K_i$. Then for any point $z\in Z$ the triangle 
inequality 
\begin{multline*}
\rho(\alpha_ig^{\infty}_i(z), \alpha_{i+k}g^{\infty}_{i+k}(z))\\
\leq \rho(\alpha_ig^{\infty}_i(z), f^{i+k}_i\alpha_{i+k}(g^{\infty}_{i+k}(z)))
+ \rho(f^{i+k}_i\alpha_{i+k}(g^{\infty}_{i+k}(z)), \alpha_{i+k}g^{\infty}_{i+k}(z))
\end{multline*}
and the property (A) imply that 
$$
\rho(\alpha_ig^{\infty}_i(z), \alpha_{i+k}g^{\infty}_{i+k}(z))
\leq \lambda_i/2+s_i.
$$ 
Then the proof follows from the Cauchy Criterion.

Denote the limit map by $\alpha$.

(C).
We show that $\alpha(Z)\subset X$. Indeed, for every $z\in Z$ the 
distance from $\alpha_ig^{\infty}_i(z)$ to the preimage 
$(f^{\infty}_i)^{-1}(\alpha_ig^{\infty}_i(z))$
does not exceed $s_i$. 
Hence $\lim_{i\to\infty}\rho(\alpha_ig^{\infty}_i(z),X)=0$.

(D).
Then we show that the inverse sequence 
$\limleftarrow \{\alpha^{-1}_i(M_{x_i}),q^{i+1}_i\}$ is well defined 
for any $x\in X$, i.e.\ we show that 
$g^{i+1}_i(\alpha_{i+1}^{-1}(M_{x_{i+1}}))
\subset\alpha_i^{-1}(M_{x_i})$.

Take an arbitrary point $y\in \alpha_{i+1}^{-1}(M_{x_{i+1}})$ and show 
that $\alpha_i(y)\in M_{x_i}$. By the triangle inequality we have
\begin{align*}
\rho(\alpha_ig^{i+1}_i(y),x_i)&
\leq \rho(\alpha_ig^{i+1}_i(y), f^{i+1}_i\alpha_{i+1}(y)) 
+ \rho(f^{i+1}_i\alpha_{i+1}(y), f^{i+1}_i(x_{i+1}))\\
&
\leq \lambda_i/4+d_{i+1}.
\end{align*}
By the condition (2) it does not exceed $\lambda_i/2$.
Hence $\alpha_ig^{i+1}_i(y)\in O_{\lambda_i/2}(x_i)\subset M_{x_i}$.

(E).
Show that 
$\alpha^{-1}(x) \supset\limleftarrow \{\alpha^{-1}_i(M_{x_i}\}$.
Let $z\in\limleftarrow \{\alpha^{-1}_i(M_{x_i})\}$. 
Since $\rho(\alpha_ig^{\infty}_i(z),x_i) \leq d_i$, then 
$\rho(\alpha_ig^{\infty}_i(z),x) \leq d_i+s_i\to 0$. 
Hence $\alpha(z)=x$ and $z\in\alpha^{-1}(x)$.

(F).
Then we show that 
$\alpha^{-1}(x) \subset \limleftarrow \{\alpha^{-1}_i(M_{x_i})\}$.

Let $z\in\alpha^{-1}(x)$ and suppose that $z$ does not belong to 
$\limleftarrow \{\alpha^{-1}_i(M_{x_i})\}$. 
Then there exists a number $i$ such that $g^{\infty}_i(z)$ does not belong 
to $M_{x_i}$. 
Therefore $\rho(\alpha_ig^{\infty}_i(z),x_i) > \lambda_i$. 
The property of the metric $\rho$ and the triangle inequality imply that 
\begin{align*}
\rho(\alpha_{i+k}g^{\infty}_{i+k}(z),x)&
\geq \rho(f^{i+k}_i\alpha_{i+k} g^{\infty}_{i+k}(z),x_i)\\
&
\geq \rho(\alpha_ig^{\infty}_i(z),x_i) 
- \rho(f^{i+k}_i\alpha_{i+k} g^{\infty}_{i+k}(z),\alpha_ig^{\infty}_i(z))\\
&
\geq \lambda_i - \lambda_i/2
\end{align*}
by (A).
Hence, $\rho(\alpha(z),x)\geq \lambda_i/2$. Contradiction.

Thus, (E) and (F) imply that 
$\alpha^{-1}(x) = \limleftarrow \{\alpha^{-1}_i(M_{x_i})\}$.
\end{proof}

We recall that a compact space $X$ is called {\it cell-like} if every map 
$f\colon X\to K$ of $X$ to a CW-complex is null homotopic. In that case 
$X$ can be imbedded in the Hilbert cube as the intersection of a nested 
sequence of sets homeomorphic to the Hilbert cube. If $X$ is finite 
dimensional, then it can be imbedded in the Euclidean space $\R^n$ as the 
intersection of a nested sequence of topological $n$-dimensional cells. 
This property of $F$ explains the name `cell-like'.

\begin{proposition}
If a compactum $X$ is the limit space of an inverse sequence of compact 
spaces with homotopy trivial bonding maps, then $X$ is cell-like.
\end{proposition}

\begin{proof}
Let $X=\limleftarrow \{X_i,p^{i+1}_i\}$. We assume that spaces $X_i$ 
are pointed and the bonding maps are point preserving. Then the system is 
realized in $\prod_{i=1}^{\infty}X_i$. Then $X=\bigcap_{k=1}^{\infty} 
(X_k\times\prod_{i=k+1}^{\infty}X_i)$.
Given a map $f\colon X\to K$ there is an extension $\bar f$ over an open 
neighborhood $O$ of $X$ in $\prod_{i=1}^{\infty}X_i$. Because of 
compactness there is a number $k$ such that 
$X_k\times\prod_{i=k+1}^{\infty}X_i\subset O$. For large enough $k$ the 
diameter of the set $\bar f(x\times\prod_{i=k+1}^{\infty}X_i)$ less then a 
given $\epsilon$ for all $x\in X_k$. For a CW-complex $K$ there is an 
$\epsilon>0$ such that every $\epsilon$-close to $f$ map $g\colon X\to K$ 
is homotopic to $f$. Take $k$ chosen for this $\epsilon$. Then the two 
maps $f$ and $\bar f\circ\pi_k$ are homotopic. Here 
$\pi_k\colon \prod_{i=1}^{\infty}X_i\to\prod_{i=1}^kX_i$ is the 
projection of the product onto the factor. Note that 
$\bar f\circ\pi_k
= \bar f\circ f^{\infty}_k
= \bar f \mymid_{X_k}\circ f^{k+1}_k\circ f^{\infty}_{k+1}$. 
Since the map $f^{k+1}_k$ is homotopically trivial, the map 
$\bar f \mymid_{X_k}\circ f^{k+1}_k\circ f^{\infty}_{k+1}$ is null 
homotopic. Hence the map $f$ is null homotopic.
\end{proof}

A map between spaces $F\colon X\to Y$ is called cell-like if $f^{-1}(y)$ 
is a cell-like set for every $y\in Y$. Since the empty set is not 
cell-like, a cell-like map is always a map onto.

\begin{theorem}[Edwards Resolution Theorem]
Let $X$ be a compactum of cohomological dimension $\dim_{\Z}X=n$. 
Then there is a compactum $Z$ of dimension $\dim Z\leq n$ and a 
cell-like map $\alpha\colon Z\to X$.
\end{theorem}

\begin{proof}
Let $X=\limleftarrow \{P_i,p^{i+1}_i\}$ be a limit space of an inverse 
sequence of compact polyhedra. We construct inverse systems 
$\{K_i,f^{i+1}_i\}$ and $\{L_i,g^{i+1}_i\}$ as in Lemma 8.3 with 
$X=\limleftarrow \{K_i,f^{i+1}_i\}$. In order to obtain a cell-like map 
$\alpha$ in the view of Proposition 8.4 we add one more condition on the 
sequences:
\begin{itemize}
\item[(4)]
a map 
$q^{i+1}_i\colon \alpha^{-1}_{i+1}(M_{x_{i+1}})\to\alpha^{-1}_i(M_{x_i})$ 
is null homotopic.
\end{itemize}

We construct sequences by induction on $i$. Let $K_1=P_1$ and let 
$\tau_1$ be a triangulation on $K_1$. We define $L_1$ as an 
$n$-dimensional skeleton $L_1^{(n)}$ of $K_1$ with respect to
triangulation $\tau_1$ and let $\alpha_1\colon L_1\to K_1$ be the 
inclusion. We define a cover $\sM^1$ of $K_1$ by closed subsets as the 
union of stars of all vertices $\{\Star(v) \mid v\in\tau^{(0)}_1\}$. Also 
we fix a metric $\rho_1$ on $K_1$.

Now assume that we have constructed sequences $\{K_i,f^i_{i-1}\}$, 
$\{L_i,g^i_{i-1}\}$, $\alpha_i\colon L_i\to K_i$ together with metrics 
$\rho_i$, triangulations $\tau_i$ on $K_i$ and covers $\sM^i$ for all 
$i\leq m$, satisfying the properties (1)-(4) of Lemma 8.3 and additionally 
$K_i=P_{r_i}$ for some $r_i$, a complex $L_i$ is an $n$-dimensional skeleton
of $K_i$ with respect to a subdivision $\tau'_m$ of the triangulation 
$\tau_i$ with the mesh $<\lambda_i/8$ and $\alpha_i$ is the inclusion map 
for all $i$. 
Also assume that a cover $\sM^i$ is defined as $\{\Star(v) \mid 
v\in\tau^{(0)}_i\}$. 
Moreover we assume that all spaces $K_i$ are pointed and, hence, 
naturally imbedded in the product $\prod_{i=1}^mK_i$ and we assume that a 
metric $\rho_i$ on each $K_i$ is the induced metric from a metric $\rho^m$ 
on the product.

We consider the Edwards-Walsh resolution $\omega\colon 
EW(\tau_m',\Z,n)\to K_m$ and apply Lemma 8.2 to the map 
$f=p^{\infty}_{r_m}\colon X\to P_{r_m}=K_m$ to obtain a lift 
$f'\colon X\to EW(\tau'_m,\Z,n)$.
Since an Edwards-Walsh space is an ANR, there is a number $k>r_m$ and a 
map $\tilde f\colon P_k\to EW(\tau_m,\Z,n)$ such that 
$\rho^m(\omega\tilde f(x),p^k_{r_m}(x))<\lambda_m/4$ for all $x\in P_k$. 
We define $K_{m+1}=P_k$, $f^{m+1}_m=p^k_{r_m}$. 
We define a metric $\rho^{m+1}$ on the product $\prod_{i=1}^{m+1}K_i$ as 
the sum of metrics $\rho^m$ on $\prod_{i=1}^mK_i$ and a metric 
$\rho_{m+1}$ bounded from above by $\frac{1}{2^m}$ on $K_{m+1}$. Then we 
imbed $K_{m+1}$ into the product $\prod_{i=1}^{m+1}K_i$ by maps 
$f^{m+1}_1,f^{m+1}_2,\dots, f^{m+1}_m, \id_{K_{m+1}}$.
Fix a base point in $K_{m+1}$ to get a canonical imbedding of 
$\prod_{i=1}^mK_i$ in $\prod_{i=1}^{m+1}K_i$.
Consider a triangulation $\tau_{m+1}$ on $K_{m+1}$ such that 
$d_{m+1}
= d(\sM^{m+1})
= d(\{\Star(v) \mid v\in\tau_{m+1}^{(0)}\})
< \lambda_m/4$ 
with respect to the metric $\rho^{m+1}$.
Then the condition (2) of Lemma 8.3 is satisfied. We define $L_{m+1}$ as 
an $n$-dimensional skeleton of a subdivision $\tau'_{m+1}$ of
$\tau_{m+1}$ with the mesh $<\lambda_{m+1}/8$ and $\alpha_{m+1}$ as the 
inclusion. 
Then the condition (1) holds.
We define $g^{m+1}_m=\omega\circ\bar f \mymid_{L_{m+1}}$ where $\bar f$ is 
a cellular approximation of $\tilde f$. Then $\omega\circ\tilde f(x)$ and 
$g^{m+1}_m(x)$ lie in one simplex of $\tau'_m$ for any $x\in L_{m+1}$. By 
the triangle inequality we have
\begin{align*}
\rho^m(g^{m+1}_m(x),f^{m+1}_m(x))&
\leq \rho^m(g^{m+1}_m(x),\omega\tilde f(x)) 
+ \rho^m(\omega\tilde f(x), p^k_{r_m}(x))\\
&
\leq \mesh\tau'_m 
+ \lambda_m/8\\
&
\leq \lambda_m/4.
\end{align*}
Hence (3) also holds for $i=m$. 
We note that by the construction $X=\limleftarrow \{K_i,f^{i+1}_i\}$. 
It means that according to Lemma 8.3 there is a map
$\alpha\colon Z\to X$ where $Z=\limleftarrow \{L_i,g^{i+1}_i\}$ with 
$\alpha^{-1}(x) 
= \limleftarrow \{M^{(n)}_{x_i},g^{i+1}_i \mymid_{\cdots}\}$. 
Note that $Z$ is at most $n$-dimensional as a limit space of 
$n$-dimensional complexes. If additionally we will have the property (4), 
then by Proposition 8.4 the map $\alpha$ will be cell-like.

Show that the condition (4) holds. For that we prove the inclusion
$$g^{m+1}_m(M^{(n+1)}_{x_{m+1}})\subset M^{(n)}_{x_m}.$$ 

Let $\Delta$ be an $n{+}1$-dimensional simplex from $M_{x_{m+1}}$. Then 
the image of $\Delta$ under the cellular map $\bar f$ lies in 
$n{+}1$-dimensional skeleton with respect to the CW-structure on the
Edwards-Walsh complex: $\bar f(\Delta)\subset EW(\tau'_m,\Z,n)^{[n+1]}$.
By the property 4-$\Z$ of the Edwards-Walsh resolution the 
$n{+}1$-skeleton $EW(\tau'_m\Z,n)^{[n+1]}$ is equal to 
$|(\tau'_m)^{(n)}|$. 
From the construction of the Edwards-Walsh complex follows 
$\bar f(\Delta)\subset\sigma^{(n)}$ for some simplex $\sigma\in\tau'_m$ 
containing $\omega\tilde f(\Delta)$.
In the proof of Lemma 8.3 part (D) it was shown that
$\alpha_ig^{i+1}_i(\alpha^{-1}_{i+1}(M_{x_{i+1}}))
\subset O_{\lambda_i/2}(x_i)$.

In our case it means that 
$g^{m+1}_m(M^{(n)}_{x_{m+1}})\subset O_{\lambda_i/2}(x_i)$. 
Hence 
$g^{m+1}_m(\partial\Delta) \subset O_{\lambda_i/2}(x_i)$. 
Hence, 
$\sigma\cap O_{\lambda_i/2}(x_i) \neq \emptyset$. 
Since $\diam\sigma<\lambda_m/8$, we have 
$\sigma \subset O_{\lambda_i}(x_i)$. 
Therefore 
$g^{m+1}_m(\Delta) \subset O_{\lambda_i}(x_i)\subset M_{x_i}$.
Since $g^{m+1}_m(\Delta)\subset |(\tau'_m)^{(n)}|$, we have the 
desired inclusion $g^{m+1}_m(\Delta)\subset M_{x_m}^{(n)}$. Since $M_x$ is 
contractible, the inclusion $M_x^{(n)}\subset M_x^{(n+1)}$ is homotopy 
trivial for any $x$. Hence, the map 
$g^{m+1}_m \mymid_{\cdots}\colon M_{x_{m+1}}^{(n)}\to M_{x_m}^{(n)}$ 
is null homotopic.
The condition (4) is checked.
\end{proof}

The following is a relative version of the Edwards Resolution theorem.

\begin{theorem}
Let $(X,A)$ be a compact pair with $\dim_{\Z}(X\setminus A)\leq n$.
Then there exists a pair $(Z,A)$ with $\dim(Z\setminus A)\leq n$ and a 
cell-like map $\alpha\colon (Z,A)\to (X,A)$ which is the identity on 
$A$.
\end{theorem}

\begin{proof}
The proof is exactly the same as in Theorem 8.5 with the only difference,
that we present $(X,A)$ as the limit space of relative polyhedra 
$(P_i,A)$ with triangulations on $P_i\setminus A$ having simplices with 
sizes tending to zero when one approaches the subset
$A$.
\end{proof}

A map between compacta $f\colon Y\to X$ is called $UV^n$-map if every 
fiber $f^{-1}(y)$ is approximately $n$-connected. We call a compactum $Z$ 
approximately $n$-connected if it has the $UV^n$-property, i.e.\ for any 
imbedding of $Z$ to ANR for every neighborhood $U\supset Z$ there is a 
smaller neighborhood $V\supset Z$ such that the inclusion $V\subset U$ 
induces a zero homomorphism for $k$-dimensional homotopy groups
$\pi_k(V)\to\pi_k(U)$ for $k\leq n$.

\begin{theorem}
Let $X$ be a compactum of the cohomological dimension $\dim_{\Z_p}X=n$. 
Then there is a compactum $Z$ of dimension $\dim Z\leq n$ and a 
$\Z_p$-acyclic $UV^{n-1}$-map $\alpha\colon Z\to X$ onto $X$.
\end{theorem}

\begin{proof}
As in the proof of Theorem 8.5 we start from an inverse system of 
polyhedra $\{P_i,p^{i+1}_i\}$ with the limit space $X$ and construct two 
inverse sequences $\sS_1=\{K_i,f^{i+1}_i\}$ and $\sS_2=\{L_i,g^{i+1}_i\}$ 
with limits $X=\limleftarrow \sS_1$ and $Z=\limleftarrow \sS_2$, 
satisfying the conditions (1)-(3) of Lemma 8.3. 
In order to get the above properties of the limit map $\alpha$ we add the 
following condition:
\begin{itemize}
\item[(4)]
the map $q^{i+1}_i\colon M^{(n)}_{x_{i+1}}\to M^{(n)}_{x_i}$ 
induces zero homomorphism in cohomologies with $\Z_p$-coefficients.
\end{itemize}

The construction of $\sS_1$ and $\sS_2$ is the same as in the proof of 
Theorem 8.5 with the only difference that we consider the Edwards-Walsh 
resolutions for the group $\Z_p$ instead of $\Z$.
We recall that $M_x$ is the star of some vertex. Hence $M_x^{(n)}$ is 
$n{-}1$-connected.
Hence the limit map $\alpha$ is approximately $n{-}1$-connected, i.e.\ 
$UV^{n-1}$.
All we have to show is that 
$(q^{i+1}_i)^*\colon H^n(M^{(n)}_{x_i};\Z_p) \to 
H^n(M^{(n)}_{x_{i+1}};\Z_p)$ 
is zero homomorphism. 
By the Universal Coefficient Theorem it suffices to show that $q^{i+1}_i$ 
induces zero homomorphism for $\Z_p$-homologies.
By the argument of Theorem 8.5 we know that 
$q^{i+1}_i = \omega \circ \bar f\mymid_{M^{(n)}_{x_{i+1}}}$. 
Denote by $h$ the restriction on $M^{(n+1)}_{x_{i+1}}$ of the map $\bar 
f\colon K_{i+1}\to EW(M_{x_i},\Z_p,n)$
defined in the proof of Theorem 8.5.
We recall $\bar f$ is a cellular map and 
$\bar f\circ i
= j\circ\omega\circ\bar f
= j\circ q^{i+1}_i$ 
where 
$i\colon M^{(n)}_{x_{i+1}}\to M^{(n+1)}_{X_{i+1}}$ 
and $j\colon M^{(n)}_{x_i}\to EW(M_{x_i},\Z_p,n)$ are inclusions. 
Hence $\bar f(M^{(n+1)}_{x_{i+1}})\subset EW(M_{x_i},\Z_p,n)^{[n+1]}$.
The property (4-$\Z_p$) of Edwrads-Walsh resolution implies that the 
inclusion $M_x^{(n)}\subset EW(M_X,\Z_p,n)^{[n+1]}$ induces a monomorphism 
of homologies with $\Z_p$-coefficients. 
Then the commutative diagram
$$
\begin{CD}
H_n(M^{(n)}_{x_{i+1}};\Z_p)	@>i_*>>	H_n(M^{(n+1)}_{x_{i+1}};\Z_p)\\
@V(q^{i+1}_i)_*VV			@Vh_*VV\\
H_n(M^{(n)}_{x_i};\Z_p)		@>j_*>>	H_n(EW(M_{x_i},\Z_p,n);\Z_p)\\
\end{CD}
$$
implies that $j_*\circ (q^{i+1}_i)_*$ is zero homomorphism. Hence 
$(q^{i+1}_i)_*$ is zero homomorphism.
\end{proof}

\begin{remark}
Let $\sL\subset\sP$ be a family of prime numbers and let 
$\dim_{\Z_p}X\leq n$ for all $p\in\sL$. 
Then there exists $\Z_p$-acyclic, $p\in\sL$, $UV^{n-1}$-map 
$f\colon Z\to X$ of $n$-dimensional compactum $Z$ onto $X$.
\end{remark}

\begin{proof}
In the construction of the inverse sequences $\sS_1$ and $\sS_2$ we apply 
the Edwards-Walsh resolutions with different $p\in\sL$ and with every $p$ 
infinitely many times. Then the result follows.
\end{proof}

The following is a relative version of Theorem 8.7.

\begin{theorem}
Let $(X,A)$ be a compact pair with $\dim_{\Z_p}(X\setminus A)\leq n$ for 
prime $p\in\sL$. Then there exists a compact pair $(Z,A)$ with 
$\dim(Z\setminus A)\leq n$ and a $Z_p$-acyclic, $p\in\sL$, $UV^{n-1}$-map 
$\alpha\colon (Z,A)\to (X,A)$ which is the identity on $A$.
\end{theorem}

\section{Resolutions preserving cohomological dimensions}

\begin{notation*}
Let $g\colon X\to K$ be a map onto a simplicial complex $K$ with a 
triangulation $\tau$.
By $\dim_G(g,\tau)\leq n$ we denote the following property of $g$:
\begin{quote}
For every subcomplex $L\subset K$ with respect to $\tau$ an every 
extension problem $\phi\colon L\to K(G,n)$ is resolved by $g$. 
\end{quote}
\end{notation*}

We note that $\dim_G(g,\tau)$ is not a number. 
We consider the inequality $\dim_G(g,\tau)\leq n$ as one symbol.

\begin{proposition}
An Edwards-Walsh resolution $\omega\colon EW(\tau,G,n)\to K$ of a finite 
complex $K$ with a triangulation $\tau$ has the property 
$\dim_G(\omega,\tau)\leq n$.
\end{proposition}

\begin{proof}
Consider a map $\phi\colon L\to K(G,n)$. It can be extended without 
problems over $n$-dimensional skeleton $K^{(n)}$ of $K$. Then by 
induction we can show that the map 
$w_m 
= \phi\omega \mymid_{\omega^{-1}(K^{(m)})\cup L}\colon 
\omega^{-1}(K^{(m)}\cup L)
\to K(G,n)$
has an extension $w_{m+1}$ over $\omega^{-1}(K^{(m+1)}\cup L)$. 
This follows from the property (3) of the Edwards-Walsh resolution. 
The union $w$ of maps $w_m$ will be a solution of the extension problem 
$\phi\omega \mymid_{\omega^{-1}(L)}$.
\end{proof}

A map $f\colon K\to L$ between two simplicial complexes is called {\it 
combinatorial} if the preimage $f^{-1}(M)$ of every subcomplex 
$M \subset L$ is a subcomplex of $K$.

\begin{lemma}
Let $X$ be a limit space of the inverse sequence of polyhedra 
$\{K_i,q^{i+1}_i\}$ with fixed metrics $\rho_i$ and fixed triangulations 
$\tau_i$ on $K_i$ such that
$$\lim_{i\to\infty}\mesh(q^{k+i}_k(\tau_{k+i}))=0$$
for all $k$. 
Assume that all bonding maps $q^{i+1}_i$ are combinatorial with respect 
to $\tau_{i+1}$ and $\tau_i$ and $\dim_G(q^{i+1}_i,\tau_i)\leq n$ for 
infinitely many $i$. 
Then $\dim_GX\leq n$.
\end{lemma}

\begin{proof}
First we note that extension problems on $X$ of the type 
$(q^{\infty}_i)^{-1}(L,\phi)$ form a basis of extension problems for the 
mappings to $K(G,n)$, where $L$ is a subcomplex of $K_i$ with respect to 
the triangulation $\tau_i$. 
By the assumption of Lemma, every such problem $(L,\phi)$ is resolved by a 
map $q^j_i$ for some $j$. 
Indeed, take $j>i$ with $\dim_G(q^j_{j-1},\tau_{j-1})\leq n$. 
Since the map $q^{j-1}_i$ is combinatorial, the problem 
$(q^{j-1}_i)^{-1}(L,\phi)$ is resolved by $q^j_{j-1}$. 
Hence, the original problem $(L,\phi)$ is resolved by $q^j_i$. 
Therefore every basic problem $(q^{\infty}_i)^{-1}(L,\phi)$ is solvable 
and, hence, by Proposition 5.4, all extension problems on $X$ to $K(G,n)$ 
have a solution. 
Hence, by Theorem 1.1, $\dim_GX\leq n$.
\end{proof}

The Proof of Lemma 8.1 can be applied to prove the following.

\begin{proposition}
Let $G$ be one of the following groups $\Z$, $\Z_p$ or $\Z_{(\sL)}$ where 
$\sL\subset\sP$ is a set of primes. 
Let $K\subset L$ be a subcomplex of a simplicial complex $L$ and let 
$\omega\colon EW(K,G,n)\to K$ be an Edwards-Walsh resolution. 
Then there exists an Edwards-Walsh resolution 
$\bar\omega\colon EW(L,G,n)\to L$ with 
$\bar\omega \mymid_{\bar\omega^{-1}(K)}=\omega$.
\end{proposition}

\begin{theorem}
Let $\sL$ be a set of prime numbers and let $n\geq 2$. Then for every 
compactum $X$ of the cohomological dimension $\dim_{\Z_{(\sL)}}X\leq n$ 
there exists a compactum $Y$ having dimensions $\dim Y\leq n+1$ and 
$\dim_{\Z_{(\sL)}}Y\leq n$ and a $\Z_{(\sL)}$-acyclic map 
$\alpha\colon Y\to X$ of $Y$ onto $X$.
\end{theorem}

\begin{proof}
We construct inverse sequences $\{K_i,f^{i+1}_i\}$, $\{L_i,g^{i+1}_i\}$ 
and a sequence of maps $\{\alpha\colon L_i\to K_i\}$ having the 
properties (1)--(3) of Lemma 8.3 with 
$X = \limleftarrow \{K_i,f^{i+1}_i\}$ and $\dim L_i= n+1$ for all $i$. 
Then a compactum $Z$ of Lemma 8.3 will be at most $n{+}1$-dimensional. 
In order to obtain the acyclicity of the map $\alpha$, we require the 
following
\begin{itemize}
\item[(4)]
A homomorphism 
$(q^{i+1}_i)^*\colon \tilde H^*(\alpha^{-1}_i(M_{x_i});\Z_{(\sL)})
\to \tilde H^*(\alpha^{-1}_{i+1}(M_{x_{i+1}});\Z_{(\sL)})$ 
is trivial.
\end{itemize}

To obtain the inequality $\dim_{\Z_{(\sL)}}Z\leq n$ we want to apply 
Lemma 9.2 and, hence, we require the existence of metrics $\hat\rho_i$ and 
triangulations $\kappa_i$ on $L_i$ such that a map $g^{i+1}_i$ is 
combinatorial with $\dim_{\Z_{(\sL)}}(g^{i+1}_i,\kappa_i)\leq n$ and with 
$$
\lim_{k\to\infty}\mesh (g^{i+k}_i(\kappa_{i+k}))=0
\quad \text{for any $i$.}
$$

We construct that by induction. 
First we fix an inverse sequence of compact polyhedra $\{P_r,p^{r+1}_r\}$ 
with limit space $X$.
By the induction on $m$ we construct the following diagram:
$$
\begin{CD}
L_1 @<g^2_1<< L_2 @<<< \cdots @<g^m_{m-1}<< L_m\\
@V\alpha_1VV @V\alpha_2VV @. @V\alpha_mVV\\
K_1 @<f^2_1<< K_2 @<<< \cdots @<f^m_{m-1}<< K_m\\
\end{CD}
$$
On each $K_i$ we define a cover $\sM^i$ with the diameter $d_i$ and the 
Lebesgue number $\lambda_i$, a triangulation $\tau_i$, a metric $\rho_i$ 
and a base point $x^*_i$. 
On each $L_i$ we define a triangulation $\kappa_i$, a metric $\hat\rho_i$ 
having the following properties
\begin{enumerate}
\renewcommand{\labelenumi}{{\normalfont (\roman{enumi})}}
\item (1) of Lemma 8.3, i.e.\ $\alpha_i(L_i)\cap M\neq\emptyset$ for 
every $M\in\sM^i$,
\item (2) of Lemma 8.3, i.e.\ $d_i<\lambda_{i-1}/4$,
\item (3) of Lemma 8.3, i.e.\ the diagram
$$
\begin{CD}
L_{i+1}	@>\alpha_{i+1}>>K_{i+1}\\
@Vg^{i+1}_iVV		@Vf^{i+1}_iVV\\
L_i	@>\alpha_i>>	K_i\\
\end{CD}
$$
is $\lambda_i/4$-commutative.
\item A homomorphism 
$(q^{i+1}_i)^*\colon \tilde H^*(\alpha^{-1}_i(M_{x_i});\Z_{(\sL)})
\to \tilde H^*(\alpha^{-1}_{i+1}(M_{x_{i+1}});\Z_{(\sL)})$ 
is trivial, where $M_{x_i}\in\sM^i$ and 
$x_i\in O_{\lambda_i}(x)\subset \Cl(M_{x_i})$
\item All spaces $K_i$ are imbedded into the product $\prod_{j=1}^mK_j$ 
by the mapping
$(f^i_1,f^i_2,\dots,f^i_{i-1},\id_{K_i},x^*_{i+1},\dots,x^*_m)$ 
and the metric $\rho_i$ is the induced from a brick metric 
$\rho^1+\cdots+\rho^m$ on the product. 
Also we assumed that $\diam_{\rho^i}K_i\leq 1/2^i$,
\item $\mesh_{\rho_i}(\tau_i)<\lambda_i/16$,
\item For every $M\in\sM^i$, $M$ is a contractible subcomplex of $K_i$ 
with respect to $\tau_i$,
\item For every $i$ there is $r(i)$ such that $K_i=P_{r(i)}$ and
$f^{i+1}_i=p^{r(i+1)}_r(i)$,
\item A complex $L_i$ has the following CW-complex structure: Take 
$n{+}1$-skeleton $K^{(n+1)}_i$ of $\tau_i$ subdivide some of its 
$n{+}1$-cells into a finite union of $n{+}1$-cells and replace some of the 
smaller $n{+}1$-cells by $n+1$-cells attached to the same boundary by 
maps of degree having all prime factors in $\sP\setminus\sL$. 
Then $\alpha_i$ is the natural projection of $L_i$ onto $K^{(n+1)}_i$ 
taking new $n{+}1$-cells to original,
\item The cellular structure on $L_i$ agrees with the triangulation 
$\kappa_i$, i.e.\ every CW-subcomplex is a simplicial complex with respect 
to $\kappa_i$,
\item Every complex $L_i$ is supplied with a metric $\hat\rho_i$ and 
$\mesh_{\hat\rho_j}(g^i_j(\kappa_i))\leq 1/2^i$ for every $j\leq i$,
\item $g^{i+1}_i$ is combinatorial and 
$g^{i+1}_i=\omega_i\circ\tilde f_i$ where 
$\omega_i\colon EW(\kappa_i, \Z_{(\sL)},n)\to L_i$ is an 
Edwards-Walsh 
resolution.
\end{enumerate}

The beginning of the induction: let $K_1=P_1$, let $\tau_1'$ be a 
triangulation on $K_1$ and let $\rho^1$ be a metric on $K_1$ of the 
diameter $\leq 1/2$. Let $\sM^1$ be a cover of $K_1$ by stars 
$M=\{\Star(v) \mid v\in(\tau_1')^{(0)}\}$. We define $\rho_1=\rho^1$ and
consider a subdivision $\tau_1$ of $\tau_1'$ with 
$\mesh_{\rho_1}(\tau_1)\leq\lambda_1/8$ where $\lambda_1$ is the Lebesgue 
number of $\sM^1$ with respect to the metric $\rho_1$.
Define $L_1$ to be the $n{+}1$-skeleton of $K_1$ with respect to 
triangulation $\tau_1$ and define $\alpha_1\colon L_1\to K_1$ as the 
inclusion. Take any metric $\hat\rho_1$ on $L_1$ and fix a triangulation 
$\kappa_1$ on $L_1$ with $\mesh_{\hat\rho_1}(\kappa_1)<1/2$. Fix a point
$x^*_1\in K_1$. All conditions (i)--(xii) are satisfied.

Now we assume that the diagram
$$
\begin{CD}
L_1 @<g^2_1<< L_2 @<<< \dots @<g^m_{m-1}<< L_m\\
@V\alpha_1VV @V\alpha_2VV @. @V\alpha_mVV\\
K_1 @<f^2_1<< K_2 @<<< \dots @<f^m_{m-1}<< K_m\\
\end{CD}
$$
is constructed satisfying the properties (i)--(xii). 
We consider the map
$$
\alpha_m\circ\omega_m\colon 
EW(\kappa_m,\Z_{(\sL)},n) \to |\tau_m^{(n+1)}|.
$$
According to the condition (ix) the homology group 
$H_n(\alpha^{-1}(\sigma))$ is a finite $(\sP\setminus\sL)$-torsion group 
for every $n{+}1$-dimensional simplex $\sigma\in\tau_m$.
The same holds true for every $n$-connected subcomplex $N\subset K_m$.
By the property (4-$\Z_{(\sL)}$) for every simplex $\sigma$ of dimension 
$\ge n+1$ there is a short exact sequence 
$$
0
\to K
\to \pi_n(\omega_m^{-1}(\alpha^{-1}(\sigma^{(n+1)})))
\to H_n(\alpha^{-1}(\sigma^{(n+1)}))
\to 0
$$
where $K$ is $\sL$-local modulo torsions. 
Hence $K/\Tor(K)=\bigoplus\Z_{(\sL)}$ and $\Tor(K)$ consists of 
$(\sP\setminus\sL)$-torsions. 
We consider an exact sequence
$$
0 
\to K/\Tor(K)
\to 
\pi_n(\omega_m^{-1}(\alpha^{-1}(\sigma^{(n+1)})))/\Tor(K)
\to H_n(\alpha^{-1}(\sigma^{(n+1)}))
\to 0.
$$
Since $Ext(G,\bigoplus\Z_{(\sL)})=0$ for any finite 
$(\sP\setminus\sL)$-torsion group $G$, 
we can present the group
$\pi_n(\omega_m^{-1}(\alpha^{-1}(\sigma^{(n+1)})))/\Tor(K)$ as the direct 
sum of $\bigoplus\Z_{(\sL)}$ and some $(\sP\setminus\sL)$-torsion group 
$G_{\sigma}$.
Thus, we have an epimorphism 
$$\pi_n(\omega_m^{-1}(\alpha^{-1}(\sigma^{(n+1)})))\to\bigoplus\Z_{(\sL)}$$ 
with a 
$(\sP\setminus\sL)$-torsion kernel $U_{\sigma}$.
Now for every $\sigma\in\tau_m$ of dimension $\ge n+1$ we kill the 
elements of the above group $U_{\sigma}$ by attaching $n{+}1$-cells. 
Then by attaching cells of higher dimensions we turn the space 
$EW(\kappa_m,\Z_{(\sL)},n)$ into a EW-resolution 
$w_m\colon EW(\tau_m,\Z_{(\sL)},n)\to K_m$ of $\tau_m$. 
Here the projection $w_m$ takes new open cells to the interior of 
corresponding simplices $\sigma$.
Since $\dim_{\Z_{(\sL)}}X\leq n$, by Lemma 8.2 there is a combinatorial 
lift $p'_m\colon X\to EW(\tau_m,\Z_{(\sL)},n)$ of 
$p^{\infty}_{r(m)}\colon X\to P_{r(m)}=K_m$ (see (8)). 
Since $EW(\tau_m,\Z_{(\sL)},n)$ is an absolute neighborhood extensor, 
there is a number $k$ and a map 
$f'_m\colon P_k\to EW(\tau_m,\Z_{(\sL)},n)$ such that
\begin{equation*}
\rho_m(w_mf'_m,p^k_{r(m)})<\lambda_m/16.\tag{$\ast$}
\end{equation*}
We define $r(m+1)=k$, $K_{m+1}=P_k$ and $f^{m+1}_m=p^k_{r(m)}$. 
Take a metric $\rho^{m+1}$ on $K_{m+1}$ of the diameter $\leq 1/2^{m+1}$ 
and define a metric $\rho_{m+1}$ on the product $\prod_{i=1}^{m+1}K_i$ as
the sum $\rho_m+\rho^{m+1}$. 
Fix a point $x^*_{m+1}\in(f^{m+1}_m)^{-1}(x^*_m)$.
The properties (v) and (viii) are satisfied.

Consider a triangulation $\tau_{m+1}'$ of $K_{m+1}$ with 
$$d_{m+1} = d(\{\Star(v) \mid v\in(\tau_{m+1}')^{(0)}\})<\lambda_m/4$$ 
and define 
$$\sM^{m+1} = \{\Star(v) \mid v\in(\tau_{m+1}')^{(0)}\}.$$ 
Then (ii) and (vii) are satisfied.

Let $\tau_{m+1}$ be a subdivision of $\tau_{m+1}'$ with 
$\mesh_{\rho_{m+1}}(\tau_{m+1})<\lambda_{m+1}/16$, where $\lambda_{m+1}$ 
is the Lebesgue number of $\sM^{m+1}$ with respect to $\rho_{m+1}$. 
Then (vi) holds.

Let $\bar f_m\colon K_{m+1}\to EW(\tau_m,\Z_{(\sL)},n)$ be a cellular 
approximation of $f'_m$ with respect to $\tau_{m+1}$ and the standard 
CW-structure on $ EW(\tau_m,\Z_{(\sL)},n)$.
By the construction the $n{+}1$-dimensional skeleton of 
$EW(\tau_m,\Z_{(\sL)},n)$ admits the following description:
$$EW(\tau_m,\Z_{(\sL)},n)^{[n+1]}
= EW(\kappa_m,\Z_{(\sL)},n)^{[n+1]}\cup_{\beta_i}B_i^{n+1},$$
where 
$\beta_i\colon \partial B^{n+1}_i
\to EW(\kappa_m,\Z_{(\sL)},n)^{[n+1]}$ 
defines a $(\sP\setminus\sL)-$torsion element $(\beta_i)_*$ 
in the homotopy group $\pi_n(EW(\kappa_m,\Z_{(\sL)},n)$.

Now we construct a finite CW-complex $L_{m+1}$ as follows. 
Consider $n{+}1$-skeleton $K^{(n+1)}_{m+1} = |\tau^{(n+1)}_{m+1}|$ 
and the restriction of $\bar f_m$ on it.
We may assume that for every $n{+}1$-simplex $\Delta$ in $\tau_{m+1}$ 
there is a partition of $\Delta$ into finitely many PL cells 
$\Delta=D^{n+1}_1 \cup \dots \cup D^{n+1}_s$ such that the image 
$\bar f_m(D^{n+1}_i)$ is an $n{+}1$-cell in 
$EW(\tau,\Z_{(\sL)},n)^{[n+1]}$. 
If $\bar f_m(D^{n+1}_i)=B_j^{n+1}$ for some $j$, we delete the interior of 
$D_i^{n+1}$ and attach an $n{+}1$-cell $\bar D_i^{n+1}$ by means of a map 
$\partial\bar D_i^{n+1}\to\partial D_i^{n+1}$ of the degree equal to the
order of the element $(\beta_j)_*$. 
We define 
$\alpha_{m+1}\colon L_{m+1}
\to |\tau^{(n+1)}_{m+1}| subset K_{m+1}$ 
by taking every cell $\bar D_i^{n+1}$ to $D_i^{n+1}$. 
Then the properties (i), (ix) hold.

Denote $N_m=EW(\kappa_m,\Z_{(\sL)},n)^{[n+1]}$. 
Now the map 
$\bar f_m \mymid_{\bar f^{-1}(N_m)}\colon \bar f^{-1}(N_m)\to N_m$ 
has an extension $\tilde f_m\colon L_{m+1}\to N_m$. 
We define $g^{m+1}_m=\omega_m\circ\tilde f_m$.
Then (xii) holds. 

Fix a metric $\hat\rho_{m+1}$ on $L_{m+1}$. 
We may assume that $L_{m+1}$ is a polyhedron and we take a triangulation 
$\kappa_{m+1}$ on it with 
$\mesh_{\hat\rho_j}(g^{m+1}_j(\kappa_{m+1})<1/2^{m+1}$ for all 
$j\leq m+1$. 
Then (x) and (xi) hold.

In order to verify (iii) we have to show that
$$\rho_m(\alpha_mg^{m+1}_m(x),f^{m+1}_m\alpha_{m+1}(x)) < \lambda_m/4.$$ 
Indeed,
\begin{align*}
\rho_m(\alpha_mg^{m+1}_m(x),f^{m+1}_m\alpha_{m+1}(x))&
\leq \rho_m(\alpha_mg^{m+1}_m(x), w_m\bar f_m\alpha_{m+1}(x))\\
&
\qquad
+ \rho_m(w_m\bar f_m(\alpha_{m+1}(x)),f^{m+1}_m(\alpha_{m+1}(x)))\\
\text{by (6) and ($\ast$)\qquad}&
< \rho_m(\alpha_m\omega_m\tilde f_m(x), w_m\bar f\alpha_{m+1}(x))
+ \lambda_m/8\\
&
= \rho_m(w_m\tilde f_m(x), w_m\bar f\alpha_{m+1}(x)) + \lambda_m/8\\
&
< \lambda_m/8 + \lambda_m/8\\
&
= \lambda/4.
\end{align*}

Now we check (iv). 
Since a complex $M_{x_{m+1}}$ is contractible, its $n$-skeleton
$M^{(n+1)}_{x_{m+1}}$ is $n$-connected. 
Hence by the construction the preimage 
$\alpha^{-1}_{m+1}(M_{x_{m+1}})
= \alpha^{-1}_{m+1}(M^{(n+1)}_{x_{m+1}})$ 
is $n{-}1$-connected. 
Note that 
$$H^n(\alpha^{-1}_{m+1}(M^{(n+1)}_{x_{m+1}});\Z_{(\sL)})
= H^n(M^{(n+1)}_{x_{m+1}};\Z_{(\sL)})=0.$$ Since $\alpha^{-1}_{m+1}(M^{(n+1)}_{x_{m+1}})$
is $n{+}1$-dimensional, it suffices to check (iv) in the dimension $n+1$. 
Note that
$H^{n+1}(EW(L,\Z_{(\sL)},n);\Z_{(\sL)})=0$.
Then by (12) $(q^{m+1}_m)^*$ is a zero homomorphism in the dimension 
$n+1$.
\end{proof}

We note that if $\sL=\emptyset$, then $\Z_{(\sL)}=\Q$.

There is a relative version of Theorem 9.4.

\begin{theorem}
Let $\sL$ be a set of primes and let $n\geq 2$. Let $(X,A)$ be a compact 
pair with $\dim_{\Z_{(\sL)}}(X\setminus A)\leq n$. Then there exists a 
compact pair $(Z,A)$ with $\dim(Z\setminus A)\leq n+1$ and 
$\dim_{\Z_{(\sL)}}(Z\setminus A)\leq n$ and a $\Z_{(\sL)}$-acyclic map 
$\alpha\colon (Z,A)\to (X,A)$ onto $X$ which is identity on $A$.
\end{theorem}

\begin{proposition}
For every finite simplicial complex $L$ there is the equality 
$$\pi_i(EW(L,\Z_p,k))=0\quad \text{for $k<i<2k-1$.}$$
\end{proposition}

\begin{proof}
Induction on $\dim L$. If $\dim L=0$, then the proposition holds.

Assume that it holds for all $m$-dimensional complexes and let $L$ be 
$m{+}1$-dimen\-sional.
We apply induction on the number of simplices in $L$. If $L$ consists of 
one simplex, then $EW(L,\Z_p,k)=K(\bigoplus\Z_p,k)$ and hence the 
proposition holds. 
Let $L=K\cup\Delta$ where $\Delta$ is a simplex of the dimension $m+1$. 
Since $i<2k-1$ and $\omega^{-1}(K')$ is $k{-}1$-connected for any 
subcomplex $K'\subset L$ where $\omega\colon EW(L,\Z_p,k)\to K'$ is
the Edwards-Walsh resolution, the Mayer-Vietoris sequence holds for 
homotopy groups:
$$\pi_i(\tilde\Delta)\oplus\pi_i(\tilde K)
\to \pi_i(L)
\to \pi_{i-1}(\tilde C)
\to \pi_i(\tilde\Delta)\oplus\pi_i(\tilde K).$$
Here by $\tilde A$ we denote the preimage $\omega^{-1}(A)$ for a 
subcomplex $A\subset K$.
Note that $\tilde A$ is an Edwards-Walsh resolution of $A$.
The induction assumption implies that $\pi_i(L)=0$ for $k+1<i<2k-1$.
Note that $\pi_k(\tilde C)\to\pi_k(\tilde\sigma)$ is a monomorphism. 
Then the exactness of the Mayer-Vietoris sequence implies that 
$\pi_{k+1}(L)=0$.
\end{proof}

\begin{proposition}
For any $p$, $k$ and any simplicial complex $L$ there exists an 
Edwards-Walsh resolution $\omega\colon EW(L,\Z_p,k)\to L$ such that 
$\omega(EW(L,\Z_p,k)^{[n]})\subset L^{(k+1)}$ for all $n<2k-1$. 
Moreover, any such resolution $\omega\colon EW(L,\Z_p,k)\to L$ given over 
a subcomplex $L\subset K$ can be extended to 
$\bar\omega\colon EW(K,\Z_p,k)\to K$ with the same property 
$\omega(EW(K,\Z_p,k)^{[n]})\subset K^{(k+1)}$ for any $n<2k-1$.
\end{proposition}

\begin{proof}
Induction on $m=\dim K$. 
If $m=k+1$ the statement is correct.
Let $\partial\Delta$ be a boundary of an $m$-simplex with some 
triangulation $\tau$. 
Assume that $\omega\colon EW(\tau,\Z_p,k)\to\partial\Delta$ be an
Edwards-Walsh resolution with the above property. 
By Proposition 9.6 $\pi_i(EW(\tau,\Z_p,k))=0$ for $k<i<2k-1$. 
Note that $\pi_i(EW(\tau,\Z_p,k))=0$ for $i<k$. 
By the property of Edwards-Walsh resolutions,
$\pi_k(EW(\tau,\Z_p,k))=\bigoplus\Z_p$. 
Therefore we can construct $EW(\Delta,\Z_p,k)$ out of $EW(\tau,\Z_p,k)$ by 
attaching cells in the dimensions $2k-1$ and higher.
Hence
$\omega(EW(\Delta,\Z_p,k)^{[n]})
= \omega(EW(\tau,\Z_p,k))
\subset (\partial\Delta)^{(k+1)}
\subset \Delta^{(k+1)}$. 
Then by induction on the number of simplices in a complex $K$ we can
construct the required Edwards-Walsh resolution.
\end{proof}

\begin{theorem}
For any set of primes $\sL\subset\sP$ and for every compactum $X$ with 
$\dim_{\Z}X\leq n$ and $\dim_{\Z_p}X\leq k$ for $p\in\sL$ with $n<2k-1$ 
there exists a compactum $Y$ with $\dim Y\leq n$ and $\dim_{\Z_p}Y\leq k$ 
for all $p\in\sL$ and a cell-like map $f\colon Y\to X$.
\end{theorem}

\begin{proof}
Define a sequence $\{p(i)\}$ of primes such that each $p\in\sL$ enters 
the sequence infinitely many times.
We construct inverse sequences of polyhedra
$$
\begin{CD}
L_1	@<g^2_1<<	L_2		@<<<	\cdots	@<g^m_{m-1}<<	L_m		@<<<	\cdots\\
@V\alpha_1VV		@V\alpha_2VV		@.			@V\alpha_mVV		@.\\
K_1	@<f^2_1<<	K_2		@<<<	\cdots	@<f^m_{m-1}<<	K_m		@<<<	\cdots\\
\end{CD}
$$
as in the proof of Theorem 9.4 with the properties
\begin{enumerate}
\renewcommand{\labelenumi}{{\normalfont (\alph{enumi})}}
\item (1) of Lemma 8.3,
\item (2) of Lemma 8.3,
\item (3) of Lemma 8.3,
\item $L_i$ is an $n$-skeleton of $K_i$ with respect to $\tau_i$ and 
$\alpha_i$ is the inclusion $K_i^{(n)}\subset K_i$,
\item A map $q^{i+1}_i\colon M^{(n)}_{x_{i+1}}\to M^{(n)}_{x_i}$ is 
null-homotopic for odd $i$,
\item (v) of Theorem 9.4,
\item (vi) of Theorem 9.4,
\item (vii) of Theorem 9.4,
\item (vii) of Theorem 9.4,
\item (x) of Theorem 9.4,
\item (xi) of Theorem 9.4,
\item $g^{i+1}_i$ is combinatorial and $g^{i+1}_i=\omega_i\circ\tilde 
f_i$, where $\omega_i\colon EW(\kappa_i,\Z_{p(\frac{i+1}{2})},k)$ for 
every odd $i$.
\end{enumerate}
Then it yields a cell-like map $\alpha\colon Z\to X$. 
Since $\dim L_i=n$, a compact $Z$ is at most $n$-dimensional. 
Propositions 9.1 and 9.2 imply that $\dim_{\Z_p}Z\leq k$ for all
$p\in\sL$.

We construct the sequences above by induction. If $m$ is even, we 
construct $\alpha_{m+1}\colon L_{m+1}\to K_{m+1}$, $g^{m+1}_m$ and 
$f^{m+1}_m$ as in the proof of Theorem 8.5.

If $m$ is odd, we consider an Edwards-Walsh resolution
$$\omega_m\colon EW(\kappa_m,\Z_{p(\frac{m+1}{2})},k)
\to L_m = K^{(n)}_m$$ 
as in Proposition 9.7. 
Again, by Proposition 9.7 there exists an extension 
$$w_m\colon EW(\tau_m,\Z_{p(\frac{m+1}{2})},k)\to K_m.$$
We construct $K_{m+1}$, $f^{m+1}_m$ and $L_{m+1}$ together with a 
cellular map $\bar f_m\colon K_{m+1}\to EW(K_m,\Z_{p(\frac{m+1}{2})},k)$ 
as in Theorem 9.4.
Then by the property of this Edwards-Walsh resolution, stated in 
Proposition 9.7, 
$w_m\circ\bar f_m(L_{m+1})
\subset K_m^{(k+1)}
\subset K_m^{(n)}
= L_m$. 
Define $g^{m+1}_m=w_m\circ\bar f_m$.
\end{proof}

\section{Imbedding and approximation}

According to the classical theorem every $n$-dimensional compactum can be 
imbedded in $\R^{2n+1}$. In this section we study the following question: 
Given cd-type $F$ find the least possible number $m$ such that $F$ has a 
representative $X$ embeddable in $R^m$. This question makes sense for 
cd-types with bounded norm $\| F\|<\infty$. The main result in this 
section is the following

\begin{theorem}
For every cd-type $F$ of $\| F\|=n$ there is a compactum 
$X\subset\R^{n+2}$ having cd-type $F$.
\end{theorem}

The proof of this theorem gives an independent proof of the Realization 
Theorem.

We recall that $M(G,n)$ denotes a Moore space, i.e.\ a CW-complex with 
trivial homology groups in dimensions $i\neq n$ and with $H_n(M(G,n))=G$.

\begin{proposition}
Suppose that the join product $L\ast M(G,1)$ is $(n{+}1)$-connected for some 
countable complex $L$ and some abelian group $G$. Then there exists an 
$n$-dimensional compactum $Y\subset\R^{n+2}$ with nontrivial Steenrod 
homology group $H_n(Y;G)\neq 0$ and with $L\in AE(Y)$.
\end{proposition}

\begin{proof}
Let $A=S^1\subset S^{n+2}$ be a circle in $(n+2)$-dimensional sphere and 
let $g\colon A\to M(G,1)$ induce a nontrivial element of the fundamental 
group $\pi_1(M(G,1))$.
Since $L\ast M(G,1)$ is $(n{+}1)$-connected, we have $L\ast M(G,1)\in 
AE(S^{n+2})$.
By Generalized Eilenberg-Borsuk theorem \cite{Dr3} there exists a 
compactum $Y\subset S^{n+2}$ with $L\in AE(Y)$ and an extension 
$\bar g\colon S^{n+2}\setminus Y\to M(G,1)$.
Since the natural inclusion $i\colon M(G,1)\to K(G,1)$ induces an 
isomorphism of the fundamental groups, the composition $i\circ g$ is a 
homotopically nontrivial map. 
Therefore $i\circ \bar g$ is a homotopically nontrivial map. 
The map $i\circ\bar g$ represents some nontrivial element 
$\alpha \in \check H^1(S^{n+2}-Y;G)$. 
By the Sitnikov duality there is a dual nontrivial element 
$\beta \in H_n(Y;G)$. 
This implies that $\dim Y\ge n$. It is possible to show that $\dim Y=n$ 
but probably the easiest way to complete the proof is by taking an
$n$-dimensional subset of $Y$.
\end{proof}

\begin{proposition}
The suspension over the smash product of two CW-complexes is homotopy 
equivalent to their joint product $\Sigma (K\wedge L)\sim K\ast L$.
\end{proposition}

\begin{proof}
Since by the definition 
$\Sigma (\limto M_{\alpha})= \limto \Sigma M_{\alpha}$,
\begin{align*}
K\wedge L&
= \limto \{L_{\alpha}\wedge K_{\beta} \mid 
\text{$L_{\alpha} \subset L$}, 
\text{$K_{\beta}\subset K$ are finite subcomplexes}\}
\quad\text{and}\\
K\ast L&
= \limto \{L_{\alpha}\ast  K_{\beta} \mid 
\text{$L_{\alpha}\subset L$, 
$K_{\beta}\subset K$ are finite subcomplexes}\},
\end{align*} 
it suffices to show that $\Sigma (K\wedge L)\sim K\ast L$ for compact 
CW-complexes.
For any pair of compact based spaces $(X,x_0)$ and $(Y,y_0)$ there is a 
closed contractible set $C=X\ast \{y_0\}\cup\{x_0\}\ast Y$ lying in $X\ast Y$ such 
that the quotient space $X\ast Y/C$ is homeomorphic to the reduced 
suspension over the smash product $X\wedge Y$. 
We note that the quotient map is a homotopy equivalence.
\end{proof}

\begin{lemma}
Suppose that two countable abelian groups have the properties $H\otimes 
G=0$ and $\Tor(H,G)=0$ ($\Tor$ means the torsion product). Then for every
$n$ there exists an $n$-dimensional compactum $Y\subset\R^{n+2}$ with
$\dim_HY\le 1$ and with nontrivial the Steenrod homology group 
$H_n(Y;G)\ne 0$.
\end{lemma}

\begin{proof}
By virtue of Proposition 10.3, we may compute homology groups 
$H_i(M(H,1)\ast M(G,1))$ via homology groups of the smash product. 
The homology group of the smash product $X\wedge Y$ is equal to the 
homology group of the pair $(X\times Y,X \vee Y)$. 
Now the homology exact sequence of the pair
$(M(H,1)\times M(G,1),M(H,1)\vee M(G,1))$ and the Kunneth formula imply 
that $H_i(M(H,1)\ast M(G,1))=0$ for all $i>0$. 
Since $\pi_1(M(H,1)\ast M(G,1))=0$, the space $M(H,1)\ast M(G,1)$ is 
$n$-connected for all $n$ by the Hurewicz theorem. 
Proposition 10.2 yields an $n$-dimensional compact $Y\subset S^{n+2}$ 
with $M(H,1)\in AE(Y)$.
By Theorem 6 of \cite{Dr4} (see also Theorem 11.4) the property 
$M\in AE(Y)$ implies the property $ SP^\infty M\in AE(Y)$ where 
$SP^\infty$ is the infinite symmetric power.
According to the Dold-Thom theorem \cite{D-T} $SP^\infty M(H,1)=K(H,1)$.

So, we have the property $K(H,1)\in AE(X)$ and hence $\dim_HY\le 1$.
\end{proof}

\begin{theorem}
For every $n$ and every $G\in\sigma$ there is a fundamental compactum
$X$ of the type $F(G,n)$ lying in $\R^{n+2}$.
\end{theorem}

\begin{proof}
We have four series of fundamental compacta. 
So, let us consider four cases.

(1) $F(\Q,n)$. 
We define $H=\bigoplus_{\text{all $p$}}\Z_p$ and $G=\Q$. 
Then the properties $G\otimes H=\Tor(G,H)=0$ hold. 
Apply Lemma 10.4 to obtain an $n$-dimensional compactum $
Y\subset \R^{n+2}$ with $\dim_HY\le 1$. 
Then it follows that $\dim_{\Z_p}Y\le 1$ for all primes $p$. 
The Bockstein inequality BI2 implies that $\dim_{\Z_{p^\infty}}Y\le 1$.
The Bockstein inequality BI5 implies $\dim_{\Z_{(p)}}Y \le \dim_{\Q} Y$ 
provided $\dim_{\Q} Y\ge 2$.
According to Lemma 10.4 $H_n(Y;\Q)\ne 0$. 
That implies $\check H^n(Y;\Q)\neq 0$ and hence $\dim_{\Q} Y\ge n\ge 2$.

Since $Y$ is $n$-dimensional, $\dim_{\Q} Y\le n$ and hence $\dim_{\Q} Y=n$.
The Bockstein inequality BI4: $\dim_{\Q} \le \dim_{\Z_{(q)}}$ completes 
the proof in the first case.

(2) $F(\Z_{(p)},n)$. 
Define $H=\bigoplus_{q\ne p}\Z_q$ and $G=\Z_{(p)}$.
Then we obtain $n$-dimensional $Y\subset \R^{n+2}$ which is 
one-dimensional with respect to $\Z_{q^\infty}$ and $\Z_q$. 
By virtue of the Bockstein inequality BI6 
$\dim_{\Z_{p^\infty}}Y
\le \max \lbrace \dim_{\Q} Y, \dim_{\Z_{(p)}}Y-1\rbrace$ 
it is sufficient to show that $\dim_{\Z_{p^\infty}}Y=n$.

Lemma 10.4 implies that $H_n(Y;G)\ne 0$. 
Therefore $\Hom(\check H^n(Y),G)\ne 0$. 
Hence the group $\check H^n(Y)$ can not be divisible by $p$. 
This means that $\check H^n(Y)\otimes \Z_{p^\infty}\ne 0$ and hence
$\check H^n(Y;\Z_{p^\infty})\ne 0$.

(3) $F(\Z_p,n)$. 
Define $H=\Z\lbrack \frac{1}{p}\rbrack$ and $G=\Z_p$. 
By Lemma 10.4 we obtain an $n$-dimensional compactum $Y\subset \R^{n+2}$ 
which is one-dimensional with respect to the groups $\Q$, $\Z_{(q)}$, 
$\Z_q$, $\Z_{q^\infty}$ ($q\ne p$) and $H_n(Y,\Z_p)\ne 0$.
Since $\Hom(\check H^n(Y),\Z_p)$ is nontrivial, the product 
$\check H^n(Y)\otimes \Z_p$ is nontrivial and hence c-$\dim_{\Z_p}Y=n$. 
The equality c-$\dim_{\Z_{(p)}}Y=n$ follows by the Bockstein theorem 
which claims that for a finite dimensional compact space $Y$ there is a 
prime $p$ such that $\dim Y = \dim_{\Z_{(p)}}Y$, and the equality 
$\dim_{\Z_{p^\infty}}Y=n-1$ follows from the Bockstein inequalities.

(4) $F(\Z_{p^\infty},n)$. 
Consider $L=M(\Z\lbrack \frac{1}{p}\rbrack,1)\vee M(\Z_p,n-1)$.

First we show that $L\ast M(\Z_{p^\infty},1)$ is an $n{+}1$-connected 
space. 
We have 
\begin{align*}
H_i(L\ast M(\Z_{p^\infty},1))&
= H_{i-1}(L\wedge M(\Z_{p^\infty},1))\\
&
= H_{i-1}(M(\Z\lbrack \tfrac{1}{p}\rbrack,1)\wedge M(\Z_{p^\infty},1))\\
&
\qquad
\oplus H_{i-1}(M(\Z_p,n-1)\wedge M(\Z_{p^\infty},1)).
\end{align*}
Since $\Z\lbrack \frac{1}{p}\rbrack \otimes \Z_{p^\infty}=0$ and
$\Tor(\Z\lbrack \frac{1}{p}\rbrack, \Z_{p^\infty})=0$, it follows that
the space $M(\Z\lbrack \frac{1}{p}\rbrack,1)\wedge M(\Z_{p^\infty},1)$ is 
contractible.
Notice that $H_{i-1}(M(\Z_p,n-1)\wedge M(\Z_{p^\infty},1))=0$ for
$i-1\le n$. 
Then the Hurewicz theorem implies that $L\ast M(\Z_{p^\infty},1)$ is
$n{+}1$-connected.

Proposition 10.2 implies that there exist an $n$-dimensional compactum 
$Y\subset\R^{n+2}$ with the property $L\in AE(Y)$.
Hence we have $M(\Z\lbrack \frac{1}{p}\rbrack,1)\in AE(Y)$ and 
$M(\Z_p,n-1)\in AE(Y)$.
Therefore $\dim_{\Z\lbrack \frac{1}{p}\rbrack}Y\le 1$ and
$\dim_{\Z_p}Y\le n-1$. 
These inequalities completely define the space $F(\Z_{p^\infty},n)$.
\end{proof}

\begin{proof}[Proof of Theorem 10.1]
By Theorem 4.11 we have $F=\bigvee\Phi(G,k_G)$. 
Since $k_G\le n$, by Theorem 10.5 every fundamental type can be realized 
by a compactum $X_G\subset\R^n$. 
The countable wedge $X=\bigvee X_G$ can be imbedded in $\R^{n+2}$. 
Note that $D_X=F$.
\end{proof}

\begin{proposition}
Let $X\subset\R^n$ be an arbitrary compactum. 
Then every map $f\colon X\to\R^n$ can be approximated by maps which do 
not change the cd-type.
\end{proposition}

\begin{proof}
Let $f\colon X\to\R^n$ be given. 
Take a compact, $n$-dimensional polyhedron $P\subset\R^n$ such that 
$X\subset P$ and extend $f$ over $P$, i.e.\ get a map 
$\bar f\colon P\to\R^n$ such that $\bar f \mymid_X=f$. 
Approximate $\bar f$ by a simplicial, general position map 
$g\colon P\to\R^n$. 
Then $g \mymid_{\Delta}\colon \Delta\to\R^n$ is an embedding for every 
simplex $\Delta$ in $P$. 
Consider $f'=g \mymid_X$. 
Since $X = \bigcup\{X\cap\Delta \mid \Delta\subset P\}$, it follows that 
$f'(X) = \bigcup\{f'(X\cap\Delta) \mid \Delta\subset P\}$. 
Then 
$D_{f'(X)} = \bigvee D_{f'(X\cap\Delta)} = \bigvee D_{X\cap\Delta} = D_X$.
\end{proof}

\begin{theorem}
For every compactum $X$ of $\dim X<n-2$ every map $f\colon X\to\R^n$ 
can be approximated by maps $f'$ with 
$D_X \preceq D_{f'(X)} \preceq D_X \vee 2$. 
\end{theorem}

\begin{proof}
Let $Z\subset\R^n$ be a realization of the cd-type of $X$ in $\R^n$ given 
by Theorem 10.1. 
Denote by $X'=Z\vee I^2\subset\R^n$. 
Let $C\subset C(X',\R^n)$ be a dense countable subset such that 
$D_{g(X')}=D_{X'}$ for all $g\in C$. 
The existence of such $C$ follows from Proposition 10.6.
Denote by $N$ the union of images $\bigcup_{g\in C}g(X')$. 
By the Completion theorem \cite{Ol} there is a $G_{\delta}$ set 
$W\supset N$ such that $\dim_GW=\dim_GN$ for all $G\in\sigma$. 
Then every compactum $Z'\subset W$ has a cd-type $\preceq$ the cd-type of 
$X'$. 
Then the complement of $W$ in $\R^n$ is a countable union of compacta 
$\bigcup Y_i$. Note that every map $q\colon X'\to\R^n$ can be approximated 
by maps avoiding $Y_i$ for every $i$. 
Then by the main Theorem of \cite{Dr5} $\dim(Y_i\times X')<n$. 
It implies that $\dim(Y_i\times X)<n$ and $\dim(Y_i\times I^2)<n$ for all 
$i$. 
The last inequality means that $\dim Y_i<n-2$. 
Then by the Disjoining Theorem \cite{D-R-S1} $f$ can be approximated by 
maps $f'$ having the empty intersection with $\bigcup Y_i$. 
Since $f'(X)\subset W$, $D_{f'(X)}\preceq D_{X'}$.
We may assume that $f'$ is a light map, then $D_X\preceq D_{f'(X)}$.
\end{proof}

\begin{corollary}
For every compactum $X$ with dimensions $\dim_GX\geq 2$ and $\dim X<n-2$ 
every map $f\colon X\to\R^n$ can be approximated with mappings $f'$ with 
$\dim_Gf'(X)=\dim_GX$.
\end{corollary}

\begin{remark*}
If $X$ is not dimensionally full-valued compactum of $\dim X=2$, say, 
$X$ is a Pontryagin surface, then according to Theorem 1.10 a map 
$f\colon X\to\R^3$ can not be approximated by maps $f'$ preserving the 
cd-type.
\end{remark*}

\section{Classifying spaces for cohomological dimension}

\begin{proposition}
The following conditions for an abelian group $G$ are equivalent:
\begin{enumerate}
\item $G$ is $p$-divisible,
\item $\Ext(\Z_{p^{\infty}}, G)=0$,
\item $\Ext(\Z_{p^{\infty}}, G)$ is $p$-divisible.
\end{enumerate}
\end{proposition}

\begin{proof}
This is a direct consequence of the short exact sequence
$$
0
\to {\lim}^1\Hom(\Z_{p^n},G)
\to \Ext(\Z_{p^{\infty}}, G)
\to \hat{G}_p
\to 0,
$$
where $\hat{G}_p=\lim G/p^nG$ is the $p$-adic completion of $G$.
\end{proof}

We recall that a space $M$ is called {\it simple} if the action of
the fundamental group $\pi_1(M)$ on all homotopy groups is trivial.
In particular this implies that $\pi_1(M)$ is abelian.

\begin{lemma}
Let $M$ be a simple CW-complex and let $X$ be a compactum. If
$\dim_{H_i(M)}X\leq i$ for all $i$, then
$\dim_{\pi_i(M)}X\leq i$ for all $i$.
\end{lemma}

\begin{proof}
Let $\pi_n=\pi_n(M)$ and $H_n=H_n(M)$. 
We prove $\dim_{\pi_n}X\leq n$ by induction on $n$. 

Since $H_1(M)=\pi_1(M)$, the claim holds for $n=1$. 

Let $\dim_{\pi_i(M)}X\leq i$ hold for all $i<n$.
For the group $\pi_n$ there is a short exact sequence
$
0
\to \left(\bigoplus_{\text{$p$ prime}} G^n_p \right)
\to \pi_n
\to F(\pi_n)
\to 0,
$
where $G^n_p$ is the Sylow $p$-subgroup of $\pi_n$ and $F(\pi_n)$ is
torsion-free.
By Lemma 2.2 it suffices to show $\dim_{F(\pi_n)}X\leq n$ and 
$\dim_{G^n_p}X\leq n$.

Let us first show that $F(\pi_n)\neq 0$ implies $\dim_{\Q}X\leq n$.
If $\pi_i$, $i<n$, are torsion groups, the Hurewicz theorem modulo the
generalized Serre class of torsion groups implies $F(H_n)\neq 0$ and 
hence $\dim_{\Q}X\leq n$. If, however, at least one of the groups $\pi_i$ 
is not a torsion group, then by the same Hurewicz theorem we obtain 
$F(H_j)\ne 0$ for some $j<n$. Therefore, $\Q\in\sigma(F(H_j))$ and 
${\dim}_{\Q}X \le \dim_{H_j}X \le j<n$.

Let $p$ be a prime number. We consider the case when $F(\pi_n)$ is not 
$p$-divisible. In that case $\Z_{(p)}\in\sigma(F(\pi_n))$. We show that 
$\dim_{\Z_{(p)}}X\leq n$. 

We may assume that all groups $H_i$, $i<n$, are $p$-divisible without
$p$-torsions. Otherwise, $\Z_p\in\sigma(H_i)$ or 
$\Z_{p^{\infty}}\in\sigma(H_i)$ and we have 
${\dim}_{\Z_p}X \le {\dim}_{H_i}X \le i < n$ or 
${\dim}_{\Z_{p^{\infty}}}X \le {\dim}_{H_i}X \le i < n$. 
In view of the inequality BI2, in both cases we have 
${\dim}_{\Z_{p^{\infty}}}X+1 \le n$. 
Then the inequality $\dim_{\Q}X\leq n$ and the Bockstein Alternative 
(Theorem 2.7) imply that $\dim_{\Z_{(p)}}X\leq n$.

Because of induction assumption, similarly we may assume that all groups
$\pi_i$, $i<n$ are $p$-divisible and without $p$-torsions.

Since $M$ is a simple CW-complex its $p$-completion $\hat{M}_p$ exists 
\cite{Bo-Ka}.
Our assumptions, Proposition 11.1 and the exact sequence 
$$
0
\to \Ext(\Z_{p^{\infty}}, \pi_i)
\to \pi_i(\hat{M}_p)
\to \Hom(\Z_{p^{\infty}}, \pi_{i-1})
\to 0
$$
imply $\pi_i(\hat{M}_p)=0$ for $i<n$. 

From the Hurewicz theorem we obtain $\pi_n(\hat{M}_p)=H_n(\hat{M}_p)$. 
This group is $\pi_n(\hat{M}_p)=\Ext(\Z_{p^{\infty}},\pi_n)$ and its 
$p$-divisibility would imply that it is the trivial group.
Since $F(\pi_n)$ is not $p$-divisible, by Proposition 11.1 
$\Ext(\Z_{p^{\infty}},F(\pi_n))$ is not $p$-divisible. Note that the 
$p$-adic completion of a torsion free group 
$\hat{F(\pi_n)}=\Ext(\Z_{p^{\infty}},F(\pi_n))$ is without torsion. The 
exactness of the sequence
$$
\Ext(\Z_{p^{\infty}},\pi_n)
\to \Ext(\Z_{p^{\infty}},F(\pi_n))
\to 0
$$
implies that 
$\Ext(\Z_{p^{\infty}}, \pi_n)
= \pi_n(\hat{M}_p)
= H_n(\hat{M}_p)$ 
is not a $p$-torsion group and is not $p$-divisible.
Therefore $H_n(\hat{M}_p)\otimes \Z_{p^{\infty}}\neq 0$ and by the 
universal coefficient theorem $H_n(\hat{M}_p;\Z_{p^{\infty}})\neq 0$.

One of the main properties of the $p$-completion $M\mapsto\hat{M}_p$ is
that it induces an isomorphism of homology with coefficients in $\Z_p$ 
\cite{Bo-Ka}. 
With exact sequences
$$
0 
\to \Z_{p^k}
\to \Z_{p^{k+1}}
\to \Z_p
\to 0
$$ 
and induction we can prove that the $p$-completion induces an isomorphism 
in homology with coefficients in $\Z_{p^n}$ for arbitrary $n$. 
Since the tensor product and homology commute with the direct limit 
the $p$-completion induces also an isomorphism in homology with 
coefficients in $\Z_{p^{\infty}}$. 

Therefore $H_n(M;\Z_{p^{\infty}})\neq 0$.
Since $H_{n-1}$ has no $p$-torsion this implies 
$H_n\otimes\Z_{p^{\infty}}\neq 0$. 
Thus and $\dim_{\Z_{(p)}}X\leq n$.

Thus, we proved the inequality $\dim_{\Z_{(p)}}X\leq n$ for all $p$ for 
which $F(\pi_n)$ is $p$-divisible. 
Since the Bockstein family $\sigma(F(\pi_n))$ consists of all such $p$'s, 
we proved the inequality $\dim_{F(\pi_n)}X\le n$.

To perform the induction step we still have to prove the inequalities 
$\dim_{G^n_p}X\le n$ for all $p$. When $F(\pi_n)$ is not $p$-divisible we 
have shown 
$\dim_{G^n_p}X
\leq \dim_{\Z_p}X
\leq \dim_{\Z_{(p)}}X
\leq n$.

Assume now $F(\pi_n)$ is $p$-divisible.
We consider two cases:

(1) $G^n_p$ is not $p$-divisible.
In this case $\sigma(G^n_p)=\{\Z_p\}$ and we have to show the inequality 
$\dim_{\Z_p}X\leq n$. 
Like above we can assume that all groups $\pi_i$, $H_i$, $i\leq n-1$,
have no $p$-torsion and are $p$-divisible. 
From the exact sequence
$$
0
\to \Ext(\Z_{p^{\infty}}, \pi_i)
\to \pi_i(\hat{M}_p)
\to \Hom(\Z_{p^{\infty}}, \pi_{i-1})
\to 0
$$
and Proposition 11.1 we obtain $\pi_i(\hat{M}_p)=0$ for $0<i<n$.
Since $G^n_p$ is not $p$-divisible, Proposition 11.1 and the exactness of 
the sequence
$$
0
= \Hom(\Z_{p^{\infty}},F(\pi_n))
\to \Ext(\Z_{p^{\infty}},G^n_p)
\to \Ext(\Z_{p^{\infty}},\pi_n)
\to \Ext(\Z_{p^{\infty}},F(\pi_n))
= 0
$$
imply that the group
$\pi_n(\hat{M}_p)
= \Ext(\Z_{p^{\infty}}, \pi_n)
= \Ext(\Z_{p^{\infty}}, G^n_p)$
is not trivial and is not $p$-divisible.

Thus the Hurewicz theorem implies $H_i(\hat{M}_p)=0$ for $0<i<n$ and the 
group $H_n(\hat{M}_p)$ is not $p$-divisible. Therefore 
$H_n(\hat{M}_p)\otimes\Z_p\neq 0$ and $H_n(\hat{M}_p;\Z_p)\neq 0$. 
From the main properties of the $p$-completion we obtain
$H_n(M;\Z_p)\neq 0$ and since $H_{n-1}$ is without $p$-torsion, 
$H_n\otimes\Z_p\neq 0$. Therefore $\Z_p\in\sigma(H_n)$ or 
$\Z_{(p)}\in\sigma(H_n)$. By virtue of Bockstein's inequality BI3 in 
both cases we have $\dim_{\Z_p}X\leq n$ and $\dim_{G^n_p}X\leq n$.

(2) $G^n_p\neq 0$ is $p$-divisible.
Then the group $\pi_n$ is $p$-divisible.

Since $\sigma(G^n_p)=\{\Z_{p^{\infty}}\}$, we have to show that
$\dim_{\Z_{p^{\infty}}}X\leq n$. We obtain this directly if $H_n$ has 
$p$-torsion elements, so assume $H_n$ has no $p$-torsion.
Again we can assume also that all the groups $\pi_i$, $H_i$, 
$1\leq i\leq n-1$, are without $p$-torsion. Therefore the exact sequence
$$
0
\to \Ext(\Z_{p^{\infty}}, \pi_i)
\to \pi_i(\hat{M}_p)
\to \Hom(\Z_{p^{\infty}}, \pi_{i-1})
\to 0
$$
implies $\pi_n(\hat{M}_p)=0$ and the group $\pi_{n+1}(\hat{M}_p)$ maps 
epimorphically onto $\Hom(\Z_{p^{\infty}}, \pi_n)$. The latter group 
includes the $p$-adic integers $\hat{\Z}_p=\limleftarrow \Z_{p^n}$ 
since $\Hom(\Z_{p^{\infty}},\Z_{p^{\infty}})\cong\hat{\Z}_p$. 
Therefore $\Hom(\Z_{p^{\infty}}, \pi_n)$ is not a $p$-torsion~group and 
since $\Z_{p^{\infty}}$ is divisible, the group 
$\Hom(\Z_{p^{\infty}}, \pi_n)$ contains 
$\Hom(\Z_{p^{\infty}},\Z_{p^{\infty}})$ which is not $p$-divisible, as a 
direct summand.
Thus the group $\pi_{n+1}(\hat{M}_p)=H_{n+1}(\hat{M}_p)$ is neither a 
$p$-torsion group nor $p$-divisible. 
Therefore $H_{n+1}(\hat{M}_p)\otimes\Z_{p^{\infty}}\neq 0$ and 
$H_{n+1}(\hat{M}_p;\Z_{p^{\infty}})\neq 0$. 
This implies $H_{n+1}(M;\Z_{p^{\infty}})\neq 0$ and since by assumption 
$H_n$ has no $p$-torsion elements the universal coefficient theorem gives 
$H_{n+1}\otimes\Z_{p^{\infty}}\neq 0$ which in turn implies 
$\dim_{\Z_{(p)}}X\leq n+1$. 

If all the groups $\pi_i$, $1\leq i\leq n-1$, are torsion groups, the 
Hurewicz theorem modulo the generalized Serre class of torsion groups 
without $p$-torsion implies that $H_n$ has $p$-torsion and thus
$\dim_{\Z_{p^{\infty}}}X\leq n$. If, however, $F(\pi_i)\neq 0$ for 
some $i$, $1\leq i\leq n-1$, we obtain $\dim_{\Q}X\leq i\leq n-1$.
Bockstein's inequality BI6 then implies $\dim_{\Z_{p^{\infty}}}X\leq n$.
\end{proof}

We recall that the $n$-th symmetric power $SP^nX$ of a space $X$ is the 
orbit space $X^n/S_n$ of the action of the symmetric group $S_n$ of degree 
$n$ by permutations of coordinates on the $n$-th power $X^n$. For a 
pointed space $X$ the inclusion $X^n\times\{x_0\}\subset X^{n+1}$ induces 
an embedding $SP^nX\to SP^{n+1}X$. The infinite symmetric power 
$SP^{\infty}X$ is the direct limit $\lim_{to}SP^nX$.

\begin{lemma}
The infinite symmetric power $SP^{\infty}M$ of a CW-complex $M$ is 
homotopy equivalent to the direct limit 
$\limto \{\prod_{i=1}^nK(H_i(M),i) \mid n\in\N\}$.
\end{lemma}

\begin{proof}
By Dold-Thom's theorem we have $\pi_i(SP^{\infty}(M))=H_i(M)$. Therefore 
there is a map of a Moore space $f_i\colon M(H_i(M),i)\to SP^{\infty}M$ 
which induces an isomorphism of $i$-dimensional homotopy groups. Note 
that the natural inclusion 
$\xi_i\colon M(H_i(M),i)\to SP^{\infty}M(H_i(M),i)$ induces an 
isomorphism of $i$-dimensional homotopy group. 
Consider a map $g_i\colon SP^{\infty}(M(H_i(M),i))\to SP^{\infty}M$ 
generated by $f_i$: regard $SP^{\infty}Y$ as the free abelian monoid over 
a space $Y$, then 
$$g_i(n_1x_1+n_2x_2+\cdots+n_mx_2)
= n_1f_i(x_1)+n_2f_i(x_2)+\cdots+n_mf_i(x_m)$$ 
where $x_j$ are points in $M(H_i(M),i)$ and $n_j\in\N$. 
Then $f_i=g_i\circ\xi_i$. Therefore $g_i$ induces an isomorphism of 
$i$-dimensional homotopy groups. We define a map 
$\mu_n\colon \prod_{i=1}^nSP^{\infty}M(H_i(M),i)\to SP^{\infty}M$ by the 
formula 
$$\mu_n((w_1,\dots,w_n)) = g_1(w_1)+\cdots+g_n(w_n).$$ 
Note that the base point in $SP^{\infty}M(H_{n+1}(M),n+1)$ defines the 
natural 
imbedding 
$\prod_{i=1}^nSP^{\infty}M(H_i(M),i)
\subset \prod_{i=1}^{n+1}SP^{\infty}M(H_i(M),i)$.
Then $\mu_{n+1}$ restricted to $\prod_{i=1}^nSP^{\infty}M(H_i(M),i)$ 
coincides with $\mu_n$. Note that for $n\ge i$ the map $g_i$ can be 
naturally factored through $\mu_n$, $mu_n\circ\gamma=g_i$. 
It implies that $\mu_n$ induces an isomorphism of homotopy groups in
dimensions $i\le n$. Hence 
$$\mu 
= \bigcup\mu_n \colon 
\limto \{\prod_{i=1}^n SP^{\infty}M(H_i(M),i) \mid n\in\N\}
\to SP^{\infty}M$$ 
is a weak homotopy equivalence. 
Since both spaces are CW-complexes, $\mu$ is a homotopy equivalence.
\end{proof}

\begin{theorem}
Let $M$ be a simple CW complex and let $X$ be a finite dimensional compactum.
Then the following are equivalent:
\begin{enumerate}
\item $M\in AE(X)$,
\item $SP^{\infty}M\in AE(X)$,
\item $\dim_{H_k(M)}X\le k$ for all $k$,
\item $\dim_{\pi_k(M)}X\le k$ for all $k$.
\end{enumerate}
\end{theorem}

\begin{proof}
(1) $\Rightarrow$ (2). 
Since $X$ is compact, it suffices to show that $SP^nM\in AE(X)$ for all 
$n$. We recall that the support ${\support}(\mu)$ of an element 
$\mu\in SP^nY\subset Y$ is the unordered set of coordinates of $\mu$. We 
may assume that $M$ is a subcomplex of a contractible complex $C$.
Then there is a natural embedding $SP^nM\subset SP^nC$ and $SP^nC$ is an
absolute extensor for compact metric spaces.
Let $\phi\colon A\to SP^nM$ be a continuous map of a closed subset 
$A\subset X$.
Then there exists an extension $\psi\colon X\to SP^nC$. 
Let 
$\Gamma_{\psi}
= \{(x,y)\in X\times C \mid y\in{\support}(\psi(x))\} 
\subset X\times C$ 
and let $F=\Gamma_{\psi}\cap(X\times M)$.
Assume that we can prove the property $M\in AE(\Gamma_{\psi})$. 
Then the map $\pi\colon F\to M$, defined by the projection $\pi(x,c)=c$, 
admits an extension $\xi\colon \Gamma_{\psi}\to M$. Consider the map
$$
\bar\phi
= SP^n(\xi)\circ SP^n(j)^{-1}\circ i\circ(\psi,\id_X) \colon X
\to SP^nM,
$$
where $j\colon \Gamma_{\psi}\to X\times C$ and 
$i\colon X\times SP^n(C)\to SP^n(X\times C)$ are the natural embeddings. 
It is easily seen that $\bar\phi$ is an extension of $\phi$ over $X$.

Now we prove the property $M\in AE(\Gamma_{\psi})$. We consider the 
following filtration on $X$: 
$X_1\subset X_2\subset\dots\subset X_n$, where
$X_k=\{x\in X \mid |{\support(\psi(x))}|\le k\}$. 
Observe that the sets $X_k$ are closed for all $k$. 
Let $p\colon \Gamma_{\psi}\to X$ be the restriction to $\Gamma_{\psi}$ of 
the projection $X\times C\to X$. 
Put $\Gamma_k=p^{-1}(X_k)$. 
In view of the Finite Union Theorem (see \cite{Dr4}), it suffices to show 
that $M\in AE(\Gamma_k)$ for all $k$. 
Since $\Gamma_1=X_1$, the condition (1) implies $M\in AE(\Gamma_1)$.
Assume that $M\in AE(\Gamma_k)$. 
The space $\Gamma_{k+1}\setminus\Gamma_k$ has a locally trivial fibration 
over the space $X_{k+1}\setminus X_k$ with $k{+}1$-point fiber. 
This implies that $M\in AE(\Gamma_{k+1}\setminus\Gamma_k)$. 
Therefore, $M\in AE(\Gamma_{k+1})$ \cite{Dr4}.

(2) $\Rightarrow$ (3).
By Lemma 11.3 we may conclude that 
$\limto \prod_{i=1}^nK(H_i(M),i)\in AE(X)$. 
Since $X$ is finite dimensional, we have 
$\prod_{i=1}^nK(H_i(M),i)\in AE(X)$ where $n=\dim X$.
Hence $K(H_k(M),k)\in AE(X)$ for all $k\le n$. 
Since $X$ is $n$-dimensional, this property holds for all $k$. 
Theorem 1.1 implies (3).

(3) $\Rightarrow$ (4). Apply Lemma 11.2.

(3) $\Rightarrow$ (1). 
By Theorem 1.1 we have $\check H^{k+1}(X,A;\pi_k(M))=0$ for every closed
subset $A\subset X$. It follows that all obstructions to an extension of 
a map $f\colon A\to M$ are vanishing. Since $X$ is finite dimensional, 
there is an extension $\bar f\colon X\to M$. Hence $M\in AE(X)$.
\end{proof}

\begin{corollary}
For finite dimensional compacta and for $k>1$ the following conditions 
are equivalent
\begin{enumerate}
\item $\dim_GX\leq k$,
\item $M(G,k)\in AE(X)$.
\end{enumerate}
\end{corollary}

This Corollary is a generalization of Alexandroff Theorem (Theorem 1.4.) 
for all abelian groups. Thus, for finite dimensional compacta Moore spaces 
are classifying spaces for the cohomological dimension as well as 
Eilenberg-MacLane spaces. 
The only possible exception is in the dimension one.

\begin{problem*}
Does the property $RP^2\in AE(X)$ hold for finite dimensional compactum 
$X$ with $\dim_{\Z_2}X=1$?
\end{problem*}

\begin{theorem}
For any compactum $X$ of dimension $\dim X= n$ and any abelian group $G$ 
such that $\dim_GX\le k$ and $k\ge 2$ there exists a closed subset 
$Y\subset X$ with $\dim Y= n-1$ and $\dim_GY\le k-1$.
\end{theorem}

\begin{proof}
By virtue of the Bockstein theorem it suffices to proof that for 
$G\in\sigma$. Since $k\ge 2$, the join product $M(G,k-1) \ast S^0$ is a Moore 
space $M(G,k)$. By Corollary 11.5, $M(G,k)\in AE(X)$. There exist two 
closed subsets $ Z^+,Z^-\subset X$ such that every separator $C\subset X$ 
has dimension $\ge n-1$. Let $f\colon \{Z^+,Z^-\}\to S^0$ be the 
separating map.
By the Generalized Eilenberg-MacLane theorem there is a compactum 
$Y\subset X$ with $M(G,k-1)\in AE(Y)$ such that $f$ is extendible to 
$X\setminus Y$. Hence $Y$ is a separator and hence $\dim Y\ge n-1$. 
By Corollary 11.5 $\dim_GY\le k-1$. 
We always may assume that $\dim Y=k-1$.
\end{proof}

\begin{theorem}
For any ring $R$, any $k\le n$ for finite dimensional compactum $X$ the 
following conditions are equivalent:
\begin{enumerate}
\item 
$\dim_R X \le n$,
\item 
every map $f\colon A \to K(R,k)$ given on a closed subset $A \subset X$ 
can be extended over to the complement $X\setminus Y$ of a compact set $Y$ 
of $\dim_R Y \le n-k-1$.
\end{enumerate}
\end{theorem}

\begin{proof}
It is sufficient to prove this theorem for rings $R\in\sigma$.

(1) $\Rightarrow$ (2). 
Let $M=M(R,n-k-1) \ast K(R,k)$ be the join product. 
It is easy to verify that $\dim_{H_k(M)}X\le\dim_RX$. Then Theorem 11.4 
yields the property $M\in AE(X)$. Then by the Generalized Eilenberg-Borsuk 
theorem \cite{Dr3} every partial map $f\colon A\to K(R,k)$ can be 
extended over the complement of compactum $Y$ with $M(R,n-k-1)\in AE(Y)$. 
By Corollary 11.5, $\dim_R Y\le n-k-1$.

(1) $\Leftarrow$ (2). 
Let $\{f_i\colon A_i\to K(R,k)\}$ be a countable basis of extension 
problems. The condition (2) gives us a compactum $Y_i$ of 
$\dim_R Y_i \le n-k-1$ and an extension 
$\bar f_i\colon X\setminus Y_i\to K(R,k)$. 
By the Countable Union theorem $\dim_R\cup Y_i\le n-k-1$.
By the Completion Theorem, there is a $G_{\delta}$ set $Z\supset\cup Y_i$ 
of $\dim_R\le n-k-1$. 
Note that every compactum $C\subset X\setminus Z$ has the property 
$K(R,k)\in AE(C)$. 
Hence, by Theorem 1.1 and the Countable Union theorem, we have 
$\dim_R(X\setminus Z)\le k$.
The Uryhson-Menger formula for the cohomological dimension \cite{Dy2}
implies that
$\dim_R X \le \dim_R Z + \dim_R (X\setminus Z) + 1 \le (n-k-1)+k+1 = n$.
\end{proof}

We note when $k=n$ the above theorem is contained in Theorem 1.1.

\section{Cohomological dimension of ANR compacta}

Absolute neighborhood retracts are locally contractible. This conditions 
gives a strict restriction on cohomological dimension. Surprisingly enough 
that locally contractible compacta can be dimensionally non-full-valued.

\begin{lemma}
Let $X$ be an ANR-compactum and $Y\in AE(X)$, then $K\in AE(X)$ for any 
CW-complex $K$ homotopy equivalent to $Y$.
\end{lemma}

\begin{proof}
Let $h\colon K\to Y$ be a weak homotopy equivalence. The important 
property of $h$ is that $h_*\colon [Z,K]\to[Z,y]$ is a bijection for all 
spaces $Z$ which are homotopy equivalent to CW-complexes. 
Let $f\colon A\to K$ is a map of a closed subset $A\subset X$. Extend $f$ 
to $f'\colon V\to K$, where $ V$ is a closed neighborhood of $A$ in $X$. 
Let $\bar f\colon X\to Y$ be an extension of $h\circ f'\colon V\to Y$. 
Take a homotopy lift $\tilde f$ of $\bar f$.
Since $h\circ\tilde f\mymid_{\Int V}$ is homotopic to 
$\bar f\mymid_{\Int V} = h\circ f'\mymid_{\Int V}$
and $\Int V$ is homotopy equivalent to a CW-complex, it follows that 
$\tilde f\mymid_{\Int V}$ is homotopic to $f'\mymid_{\Int V}$. 
Hence, $\tilde f\mymid_A\sim f'\mymid_A=f$. 
Thus, $f$ extends over $X$ up to homotopy, so it extends over $X$ by the 
Homotopy Extension Theorem.
\end{proof}

\begin{theorem}
Let $G=\prod_{s\in S}G_s$ be the direct product of abelian groups. 
Then
$\dim_GX = \max\{\dim_{G_s}X \mid s\in S\}$
for any compactum $X$.
\end{theorem}

\begin{proof}
Since each $G_s$ is a direct summand of $G$, Corollary 1.7 implies that
$\dim_{G_s}X\le\dim_GX$. 
Hence, $\max\{\dim_{G_s}X \mid s\in S\}\le\dim_GX$.
Suppose that $\max\{\dim_{G_s}X \mid s\in S\}=n$.
Note that $Y=\prod_{s\in S}K(G_s,n)\in AE(X)$. 
Note that $Y$ is weakly homotopy equivalent to $K(G,n)$. 
By Lemma 12.1, $K(G,n)\in AE(X)$, hence, by Theorem 1.1 
$\dim_G \le n = \max\{\dim_{G_s}X \mid s\in S\}$.
\end{proof}

\begin{theorem}
Let $X$ be an ANR compactum, then
\begin{enumerate}
\item $\dim_{\Z_{(p)}}X=\dim_{\Z_p}X$ for all prime $p$,
\item $\dim_GX\ge\dim_{\Q}X$ for any abelian group $G\neq 0$.
\end{enumerate}
\end{theorem}

\begin{proof}
(1). In view of the Bockstein inequality BI3 it suffices to show that 
$\dim_{\Z_{(p)}}X\le \dim_{\Z_p}X$. Consider $G=\prod_k\Z_{p^k}$. Then by 
Theorem 12.2 and Proposition 2.3, 
$\dim_GX=\max\{\dim_{\Z_{p^k}}X\}=\dim_{\Z_p}X$. Since $G$ contains an 
element of infinite order and not divisible by $p$, we have 
$\Z_{(p)}\in\sigma(G)$. By the Bockstein theorem (Theorem 2.1) 
$\dim_{\Z_{(p)}}X\le\dim_GX=\dim_{\Z_p}X$.

(2). If $\Q\in\sigma(G)$, then the inequality follows from Theorem 2.1.

If $\Z_{(p)}\in\sigma(G)$, then the inequality BI4 implies the required 
inequality.

If $\Z_p\in\sigma(G)$, then we apply BI4, (1) and Theorem 2.1 
to obtain 
$\dim_{\Q}X
\le \dim_{\Z_{(p)}}X
= \dim_{\Z_p}X
\le \dim_GX$.

If $\Z_{p^{\infty}}\in\sigma(G)$, then we consider the group 
$A = \Z_{p^{\infty}}\times\Z_{p^{\infty}}\times \cdots$. 
By Theorem 12.2 $\dim_AX=\dim_{\Z_{p^{\infty}}}X\le\dim_GX$. 
Since $A$ is not a torsion group and it is divisible by all $p$, by the 
definition $\Q\in\sigma(A)$. 
By Theorem 2.1 
$\dim_{\Q} X \le \dim_A X \le \dim_G X$.
\end{proof}

\begin{corollary}
Every ANR-compactum $X$ is of the basic type, i.e.\ the formula 
$\dim(X\times X) = 2 \dim X$ holds for $X$.
\end{corollary}

\begin{proof}
We consider a finite dimensional ANR compactum $X$. By Theorem 1.4 and 
Theorem 2.1, $\dim X=\dim_{\Z}X=\dim_{\Z_{(p)}}X$ for some prime $p$. 
By Theorem 12.3, $\dim X=\dim_{\Z_p}X$. 
Then Criterion 3.17 completes the proof.
\end{proof}

\begin{theorem}
Every 2-dimensional ANR compactum $X$ is dimensionally full-valued.
\end{theorem}

\begin{proof}
By the Universal Coefficient Theorem the simplicial 1-dimensional 
cohomology is a free abelian group. Therefore the \v{C}ech 1-dimensional 
cohomology is a torsion free group. Hence by the Universal Coefficient 
Formula $\check H^1(A;\Q)\ne 0$ for any $A$. Take a closed neighborhood 
$U\subset X$ which is contractible in $X$ and with $\dim U=2$. Then there 
is a compact subset $A\subset U$ with $\check H^2(U,A)\ne 0$. The 
homomorphism $\gamma\colon \check H^2(U,A)\to \check H^2(U)$ is trivial, 
since $\phi$ is trivial in the following diagram and $\alpha$ is 
surjective because of 2-dimensionality of $X$.
$$
\begin{CD}
\check H^2(U) 
@<\gamma<< 
\check H^2(U,A) 
@<<< 
\check H^1(A)\\
@A{\phi}AA 
@A{\alpha}AA 
@.\\
\check H^2(X) 
@<\nu<< 
\check H^2(X,A)\\
\end{CD}
$$
Therefore $\check H^1(A)\ne 0$ and hence $\check H^1(A;\Q)\ne 0$.
Since the inclusion $A\subset X$ is homotopically trivial, the induced 
homomorphism in rational cohomologies is trivial. Hence $\check 
H^2(X,A;\Q)\ne 0$. Now by Theorem 12.3 $X$ is dimensionally full-valued.
\end{proof}

\begin{theorem}
For any prime $p$ there exists an AR compactum $M_p$, having dimensions:
$\dim M_p = \dim_{\Z_{(p)}}M_p = \dim_{\Z_p}M_p = 4$ and
$\dim_{\Q}M_p = \dim_{\Z_{p^{\infty}}}M_p = \dim_{\Z_q}M_p = 3$ where 
$q\ne p$ is prime.
\end{theorem}

For a map $f\colon A\to B$ we denote by 
$S_f=\{x\in A \mid f^{-1}f(x)\ne x\}\subset A$
the singularity set of $f$.
We use the following theorem which generalizes Borsuk's ANR pasting 
theorem.

\begin{theorem}
Let $A,B,X$ be ANR compacta and let $\alpha\colon A\to X$ and 
$f\colon A\to B$ have the following property: $\alpha$ restricted to the 
singularity set $S_f$ is one-to one.
Then the pushout $Y$ of the diagram
$$
\begin{CD}
A @>f>> B\\
@V{\alpha}VV @.\\
X\\
\end{CD}
$$
is an ANR compactum provided it is finite dimensional.
\end{theorem}

\begin{proof}
Consider the diagram:
$$
\begin{CD}
A @>f>> B \\
@V{\alpha}VV @V{\beta}VV\\
X @>\phi>> Y\\
\end{CD}
$$
Since $\alpha$ is injective on $S_f$, the map $\phi$ is defined by the 
decomposition 
$\sF = \{\alpha f^{-1}(y) \mid \text{$y\in f(S_f)$ and singletons}\}$. 
It is clear that this decomposition is upper semicontinuous. 
Hence $Y$ is a compact metric space.
There is the natural map $q\colon DM_{\alpha,f}\to Y$ of the double 
mapping cylinder onto the pushout. 
By Borsuk's ANR pasting theorem $DM_{\alpha,f}$ is ANR.
We show that map $q$ is cell-like, then the result follows.

We consider three cases.

(1) $y\in Y\setminus\phi\alpha A$. 
In that case $q^{-1}$ is a singleton, i.e.\ is cell-like.

(2) $y\in \phi\alpha A\setminus\phi\alpha S_f$. 
In this case $q^{-1}$ is homeomorphic to the cone over $\alpha^{-1}(x)$ 
where $\phi(x)=y$. Hence it is cell-like.

(3) $y\in\phi\alpha S_f$. 
In this case the restriction of $\alpha$ on $\alpha^{-1}\phi^{-1}(y)$ is a 
retraction $r$ onto $\phi^{-1}(y)$. 
Let $S=\alpha^{-1}\phi^{-1}(y)\cap S_f$. 
Then $q^{-1}$ is homeomorphic to the union of the mapping cylinder of $r$ 
and the cone over $S\subset\alpha^{-1}\phi^{-1}(y)$. We can define a 
contraction of this union to a point as follows. First we can deform the 
mapping cylinder $M_r$ to the image space $\phi^{-1}(y)\cong S$. This 
deformation can be extended to a homotopy of the whole $q^{-1}(y)$. As the 
result we have a deformation of the space to the union of the mapping 
cylinder of $\alpha$ restricted over $S$ and the cone over $S$. Since this 
is homeomorphic to the cone over $S$, we can contract that to a point.
\end{proof}

\begin{lemma}
There is an imbedding of an infinite tree $T=\bigcup T_i$ in a 
four-dimensional cube $I^4$ such that there is a sequence of regular 
neighborhoods $N_1\subset N_2\subset\cdots$ of the finite trees 
$T_1\subset T_2\subset\cdots$ with the properties:
\begin{enumerate}
\item The union $\bigcup N_i=N$ is dense in $I^4$,
\item For every $i$ there is an $\epsilon_i$-retraction
$h_i\colon N_{i+1}\setminus \Int(N_i)
\to \Cl(\partial N_{i+1}\setminus\partial N_i)$,
\item $\sum \epsilon_i < \infty$,
\item The restriction 
$h_i\mymid_{\Cl(\partial N_i\setminus\partial N_{i+1})}$
is an imbedding.
\end{enumerate}
\end{lemma}

\begin{proof}
We construct $T$ and $N$ by induction. Assume that $\diam(I^4)=1$ and 
choose a point $x_0\in\partial I^4$. We define $T_1$ as the segment from 
$x_0$ to the center $c$ of the cube $I^4$. Take $\epsilon_1=2$ and let 
$N_1$ be a regular neighborhood of $T_1$ in $I^4$. There is an 
$\epsilon_1$-retraction 
$h_1\colon N_1\to \Cl(\partial N_1\setminus\partial I^4)$.
Consider a finite 1/2-net in $\Int(I^4\setminus N_1)$. 
Then we join points of the net by smooth arcs in $I^4$ of the length 
$\le 1$ with $c$.
We may assume that all arcs are disjoint and transversal to $\partial N_1$.
The union of these arcs with $T_1$ gives $T_2$. Then we consider a 
regular neighborhood $N'_2$ of $T_2$ such that there is an 
$\epsilon_2$-retraction
$h_2\colon N'_2\setminus \Int(N_1)\to \Cl(\partial N'_2\setminus N_1)$ 
with $\epsilon_2=1$.
Define $N_2=N'_2\cup N_1$. Consider a 1/4-net in $\Int(I^4\setminus N_2)$ 
and join every point of the net with one of the closest point of the 
previous net by an arc of length $\le 1/2$ and transversal to 
$\partial N_2$ and so on.
\end{proof}

\begin{proof}[Proof of Theorem 12.6]
Let $N$ and $T$ be as above.
Let $A=B=N\cap\partial I^4=D$ be a 3-dimensional disk. 
Define $X=I^4\setminus \Int(N)$.
Since $\sum \epsilon_i < \infty$, the composition 
$\bar h = \dots \circ h_2 \circ h_1$
is a retraction of $I^4$ onto $X$. 
Hence $X\in AR$. 
We define $\alpha=\bar h\mymid_D$.
We define $f\colon D\to D$ as follows. 
Denote 
$D_k=\alpha^{-1}(\bigcup_{i=1}^k\Cl(\partial N_i\setminus N_{i+1})$.
Then we define $f_0\colon \partial D\to\partial D$ as a map of degree $p$.
Since the second homotopy group is abelian, we can extend $f_0$ to 
$f_1\colon D_1\to D_1$ in such way that the restriction of $f_1$ on every 
component of the boundary $\partial D_1$ is a mapping the component to 
itself with the degree $p$. 
Then we can extend $f_1$ to $f_2\colon D_2\to D_2$ in the similar fashion 
and so on. 
Let 
$U=\bigcup_{i=1}^{\infty}\alpha^{-1}\Cl(\partial N_i\setminus N_{i+1})$.
Then $D\setminus U=C$ is a Cantor set. We define $f$ on $U$ as the union 
of $f_i$ and $f\mymid_C=\id_C$. We note that $\alpha\mymid_U$ is injective 
and $S_f\subset U$. 
Also it is easy to see that the pushout in this case is at most 
4-dimensional.
Then Theorem 12.7 defines an $AR$-space $M_p=Y$.

Note that $Z=\phi(\Cl(\partial I^4\setminus D))$ is homeomorphic to the
Moore space $M(\Z_p,2)$. 
Since $H^3(Z;\Z_p)\ne 0$ and $Z\subset Y\in AR$, the exact sequence of 
pair $(Y,Z)$ implies that $\dim_{\Z_p}Y\ge 4$.
Therefore the Bockstein Theorem and the Alexandroff Theorem together with 
BI3 imply that
$$\dim M_p
= \dim_{\Z_{(p)}}M_p
= \dim_{\Z_p}M_p
= 4.$$
We show that for every closed subset $F\subset Y$ the equality 
$\check H^3(F;\Z_q)=0$ holds for all prime $q\ne p$. 
Then Theorem 12.3 and the Bockstein Alternative imply that
$$\dim_{\Q}M_p
= \dim_{\Z_q}M_p
= \dim_{\Z_{p^{\infty}}}M_p
= 3.$$

Let $K=\beta^{-1}(F)$. There is a sequence of open 3-balls $\{B_i\}$ in $D$
such that
\begin{enumerate}
\item each ball is a component of a complement to $D_l$ for some $l$,
\item $C\setminus K\subset\bigcup_{i=1}^{\infty}B_i$,
\item $B_i\cap K=\emptyset$.
\end{enumerate}
Denote $D'=D\setminus\bigcup_{i=1}^{\infty}B_i$ and consider 
$F'=\beta(D')$.
We show that the inclusion $F\subset F'$ induces an epimorphism in 
3-dimensional cohomologies. Let $g\colon F\to K(G,3)$ be a map to 
Eilenberg-MacLane complex.
Since $\dim D'=3$, there is an extension $\nu\colon D'\to K(G,3)$ of a 
map $g\circ\beta\mymid_{\beta^{-1}(F)}$. 
We define $\bar g\colon F'\to K(G,3)$ by the formula: 
$\bar g(z)=\nu\beta^{-1}\mymid_{D'}(z)$ for $z\in F'$. 
Since $\bar g$ is an upper semi-continuous multi-valued map, it suffices 
to show that $\nu\beta^{-1}(z)$ consists of one point for all $z\in F'$. 
By the definition this holds for $z\in F$.
Let $z\in F'\setminus F$. 
Then by the definition of $D'$ we have that 
$\beta^{-1}(z)\cap D' \subset U = f(U) \subset S_f$. 
Since $\alpha\mymid_{S_f}$ is injective, 
$|f(f^{-1}\beta^{-1}(y))\cap S_f| \le 1$. 
This implies that $\beta^{-1}(z)\cap D'$ consists of one point. 
Next, we show that $\check H^3(F';\Z_q)=0$. 
We consider the map $\gamma\colon M_{\alpha}\to DM_{\alpha,f}$
generated by the map $f$. Let $\bar q\colon DM_{\alpha,f}\to Y$ be the 
quotient map of Theorem 12.7. 
Consider the diagram generated by the map $\gamma$ restricted to the 
pairs $(\bar q^{-1}(F'),D')$ and $(\gamma^{-1}\bar q^{-1}(F'),D')$ where 
$D'$ is considered here as the subset of $D=B$ and $D=A$ respectively;
$$
\begin{CD}
0	@<<<	\check H^3(\bar q^{-1}(F');\Z_q)		@<\gamma<<	\check H^3(\bar q^{-1}(F'),D';\Z_q)		@<<<	\check H^2(D';\Z_q)\\
@.		@V{\gamma^*}VV							@V{\phi_2}VV						@V{\phi_3}VV\\
0	@<<<	\check H^3(\gamma^{-1}\bar q^{-1}(F');\Z_q)	@<<<		\check H^3(\gamma^{-1}\bar q^{-1}(F'),D';\Z_q)	@<<<	\check H^2(D';\Z_q)\\
\end{CD}
$$
The homomorphism $\phi_2$ is generated by a relative homeomorphism and, 
hence, is an isomorphism. The homomorphism $\phi_3$ is generated by the 
restriction $f\mymid_{D'}$ which is a map of degree $p$ of an infinite 
wedge of spheres to itself. Hence it induces an isomorphism of 
cohomologies with coefficients in $\Z_q$ for $q$ relatively prime to $p$. 
The Five lemma implies that $\gamma^*$ is an isomorphism.

Let $\bar \alpha\colon M_{\alpha}\to X$ be the natural projection to the 
range.
The diagram
$$
\begin{CD}
DM_{\alpha,f} @>{\bar q}>> Y\\
@A{\bar\gamma}AA @A{\phi}AA\\
M_{\alpha} @>{\bar\alpha}>> X\\
\end{CD}
$$
restricted to $F'\subset Y$ produces isomorphisms diagram for
3-dimensional cohomology:
$$
\begin{CD}
\check H^3(\bar q^{-1}(F');\Z_q)
@<<< \check H^3(F';\Z_q)\\
@V{\gamma^*}VV @VVV\\
\check H^3(\gamma^{-1}\bar q^{-1}(F');\Z_q)
@<<< \check H^3(\phi^{-1}(F');\Z_q)\\
\end{CD}
$$
Since $X$ is 3-dimensional $AR$-space, $\check H^3(\phi^{-1}(F');\Z_q)=0$.
Hence $\check H^3(F';\Z_q)=0$ and $\check H^3(F;\Z_q)=0$.
\end{proof}

\begin{remark} 
For relatively prime $p$ and $q$ the dimension of the product does not 
comply to the logarithmic law: $\dim M_p\times\dim M_q=7$.
\end{remark}

\begin{proof}
By Alexandroff and Bockstein theorems we have
\begin{align*}
\dim(M_p\times M_q)&
= \dim_{\Z}(M_p\times M_q)\\
&
= \max\{\dim_{\Z_{(r)}}(M_p\times M_q)\}
\quad \text{by Theorem 12.3}\\
&
= \max\{\dim_{\Z_r}(M_p\times M_q)\}\\
&
= \max\{\dim_{\Z_r}M_p+\dim_{\Z_r}M_q\}\\
&
= 7.
\end{align*}
\end{proof}


\begin{thebibliography}{10}

\bibitem{A}
P.S. Alexandroff, \emph{Dimensiontheorie, ein {B}etrag zur {G}eometire der
  abgeschlossen {M}engen}, Math. Ann. \textbf{106} (1932), 161--238.

\bibitem{A2}
\bysame, \emph{Einige {P}roblemstellungen in der mengentheoretischen
  {T}opologie}, Math. Sbor. \textbf{43} (1936), 619--634.

\bibitem{A1}
\bysame, \emph{Introduction to homological dimension theory and general
  combinatorial topology}, Izdat. ``Nauka'', Moscow, 1975 (Russian). \MR{58
  \#24234a}

\bibitem{A-H}
D.W. Anderson and L.~Hodgkin, \emph{The {$K$}-theory of {E}ilenberg-{M}ac{L}ane
  complexes}, Topology \textbf{7} (1968), 317--329. \MR{37 \#6924}

\bibitem{Bo1}
M.F. Bockstein, \emph{Homology invariants of topological spaces}, Trudy Moskov.
  Mat. Ob\v{s}\v{c}. \textbf{5} (1956), 3--80 (Russian). \MR{18,813h}

\bibitem{Bo2}
\bysame, \emph{Homological invariants of topological spaces. {II}}, Trudy
  Moskov. Mat. Ob\v{s}\v{c}. \textbf{6} (1957), 3--133 (Russian). \MR{19,875a}

\bibitem{Bol2}
V.G. Boltyanskij, \emph{An example of a two-dimensional compactum whose
  topological square is three-dimensional}, Doklady Akad. Nauk SSSR (N.S.)
  \textbf{67} (1949), 597--599. \MR{11,45e}

\bibitem{Bol1}
\bysame, \emph{On the dimensional full-valuedness of compacta}, Doklady Akad.
  Nauk SSSR (N.S.) \textbf{67} (1949), 773--776. \MR{11,195j}

\bibitem{Bo-Ka}
A.K. Bousfield and D.M. Kan, \emph{Homotopy limits, completions and
  localizations}, Lecture Notes in Math., no. 304, Springer-Verlag, Berlin,
  1972. \MR{51 \#1825}

\bibitem{B-M}
V.M. Buh{\v{s}}taber and A.S. Mi{\v{s}}{\v{c}}enko, \emph{A {$K$}-theory on the
  category of infinite cell complexes}, Izv. Akad. Nauk SSSR Ser. Mat.
  \textbf{32} (1968), 560--604. \MR{40 \#6537}

\bibitem{Co}
H.~Cohen, \emph{A cohomological definition of dimension for locally compact
  {H}ausdorff spaces}, Duke Math. J. \textbf{21} (1954), 209--224. \MR{16,609b}

\bibitem{D-T}
A.~Dold and R.~Thom, \emph{Quasifaserungen und unendliche symmetrische
  {P}rodukte}, Ann. of Math. (2) \textbf{67} (1958), 239--281. \MR{20 \#3542}

\bibitem{Dr8}
A.N. Dranishnikov, \emph{The cohomological dimension is not preserved under the
  {S}tone-\v{C}ech compactification}, C. R. Acad. Bulgare Sci. \textbf{41}
  (1988), no.~12, 9--10. \MR{90e:55002}

\bibitem{Dr1}
\bysame, \emph{Homological dimension theory}, Uspekhi Mat. Nauk \textbf{43}
  (1988), no.~4(262), 11--55, 255. \MR{90e:55003}

\bibitem{Dr2}
\bysame, \emph{On a problem of {P}. {S}. {A}leksandrov}, Mat. Sb. (N.S.)
  \textbf{135(177)} (1988), no.~4, 551--557, 560. \MR{90e:55004}

\bibitem{Dr7}
\bysame, \emph{Generalized cohomological dimension of compact metric spaces},
  Tsukuba J. Math. \textbf{14} (1990), no.~2, 247--262. \MR{92d:55002}

\bibitem{Dr4}
\bysame, \emph{Extension of mappings into {CW}-complexes}, Mat. Sb.
  \textbf{182} (1991), no.~9, 1300--1310. \MR{93a:55002}

\bibitem{Dr5}
\bysame, \emph{On intersection of compacta in {E}uclidean space. {II}}, Proc.
  Amer. Math. Soc. \textbf{113} (1991), no.~4, 1149--1154. \MR{92c:54015}

\bibitem{Dr12}
\bysame, \emph{On intersections of compacta in {E}uclidean space}, Proc. Amer.
  Math. Soc. \textbf{112} (1991), no.~1, 267--275. \MR{91h:54024}

\bibitem{Dr13}
\bysame, \emph{{$K$}-theory of {E}ilenberg-{M}ac {L}ane spaces and cell-like
  mapping problem}, Trans. Amer. Math. Soc. \textbf{335} (1993), no.~1,
  91--103. \MR{93e:55003}

\bibitem{Dr6}
\bysame, \emph{The {E}ilenberg-{B}orsuk theorem for mappings in an arbitrary
  complex}, Mat. Sb. \textbf{185} (1994), no.~4, 81--90. \MR{95j:54028}

\bibitem{Dr3}
\bysame, \emph{On the mapping intersection problem}, Pacific J. Math.
  \textbf{173} (1996), no.~2, 403--412. \MR{97e:54030}

\bibitem{Dr9}
\bysame, \emph{On the virtual cohomological dimensions of {C}oxeter groups},
  Proc. Amer. Math. Soc. \textbf{125} (1997), no.~7, 1885--1891. \MR{98d:55001}

\bibitem{Dr11}
\bysame, \emph{Rational homology manifolds and rational resolutions}, Topology
  Appl. \textbf{94} (1999), no.~1-3, 75--86. \MR{2000j:55001}

\bibitem{Dr10}
\bysame, \emph{On the dimension of the product of two compacta and the
  dimension of their intersection in general position in {E}uclidean space},
  Trans. Amer. Math. Soc. \textbf{352} (2000), no.~12, 5599--5618.
  \MR{2001j:55002}

\bibitem{D-D2}
A.N Dranishnikov and J.~Dydak, \emph{Extension dimension and extension types},
  Tr. Mat. Inst. Steklova \textbf{212} (1996), no.~Otobrazh. i Razmer., 61--94.
  \MR{99h:54049}

\bibitem{D-D1}
A.N. Dranishnikov and J.~Dydak, \emph{Extension theory of separable metrizable
  spaces with applications to dimension theory}, Trans. Amer. Math. Soc.
  \textbf{353} (2001), no.~1, 133--156. \MR{2001f:55002}

\bibitem{D-D-W}
A.N. Dranishnikov, J.~Dydak, and J.J. Walsh, \emph{Cohomological dimension
  theory with applications}, 1992, Preprint.

\bibitem{D-R}
A.N. Dranishnikov and D.~Repov{\v{s}}, \emph{Cohomological dimension with
  respect to perfect groups}, Topology Appl. \textbf{74} (1996), no.~1-3,
  123--140. \MR{98e:55002}

\bibitem{D-R-S3}
A.N. Dranishnikov, D.~Repov{\v{s}}, and E.V. {\v{S}}{\v{c}}epin,
  \emph{Dimension of products with continua}, Topology Proc. \textbf{18}
  (1993), 57--73. \MR{96b:54054}

\bibitem{D-R-S1}
\bysame, \emph{On approximation and embedding problems for cohomological
  dimension}, Topology Appl. \textbf{55} (1994), no.~1, 67--86. \MR{94m:55001}

\bibitem{D-R-S2}
\bysame, \emph{On the failure of the {U}rysohn-{M}enger sum formula for
  cohomological dimension}, Proc. Amer. Math. Soc. \textbf{120} (1994), no.~4,
  1267--1270. \MR{94f:55001}

\bibitem{D-R-S4}
\bysame, \emph{Transversal intersection formula for compacta}, Topology Appl.
  \textbf{85} (1998), no.~1-3, 93--117. \MR{2000f:55001}

\bibitem{D-S}
A.N. Dranishnikov and E.V. Shchepin, \emph{Cell-like mappings. {T}he problem of
  the increase of dimension}, Uspekhi Mat. Nauk \textbf{41} (1986), no.~6(252),
  49--90, 230. \MR{88g:57021}

\bibitem{D-W}
A.N. Dranishnikov and J.E. West, \emph{Compact group actions that raise
  dimension to infinity}, Topology Appl. \textbf{80} (1997), no.~1-2, 101--114.
  \MR{98i:55006}

\bibitem{Dy1}
J.~Dydak, \emph{Cohomological dimension and metrizable spaces}, Trans. Amer.
  Math. Soc. \textbf{337} (1993), no.~1, 219--234. \MR{93g:55001}

\bibitem{Dy2}
\bysame, \emph{Cohomological dimension and metrizable spaces. {II}}, Trans.
  Amer. Math. Soc. \textbf{348} (1996), no.~4, 1647--1661. \MR{96h:55001}

\bibitem{Dy3}
\bysame, \emph{Cohomological dimension theory}, 1997, Preprint.

\bibitem{D-W1}
J.~Dydak and J.J. Walsh, \emph{Spaces without cohomological dimension
  preserving compactifications}, Proc. Amer. Math. Soc. \textbf{113} (1991),
  no.~4, 1155--1162. \MR{92c:54039}

\bibitem{D-W3}
\bysame, \emph{Complexes that arise in cohomological dimension theory: a
  unified approach}, J. London Math. Soc. (2) \textbf{48} (1993), no.~2,
  329--347. \MR{94d:55003}

\bibitem{D-W2}
\bysame, \emph{Infinite-dimensional compacta having cohomological dimension
  two: an application of the {S}ullivan conjecture}, Topology \textbf{32}
  (1993), no.~1, 93--104. \MR{94b:55002}

\bibitem{Dy}
E.~Dyer, \emph{On the dimension of products}, Fund. Math. \textbf{47} (1959),
  141--160. \MR{22 \#986}

\bibitem{Ed}
R.D. Edwards, \emph{A theorem and a question related to cohomological dimension
  and cell-like maps}, Notices Amer. Math. Soc. \textbf{25} (1978), A--259.

\bibitem{H-W}
W.~Hurewicz and H.~Wallman, \emph{Dimension {T}heory}, Princeton Mathematical
  Series, vol.~4, Princeton University Press, Princeton, N. J., 1941.
  \MR{3,312b}

\bibitem{Ko1}
Y.~Kodama, \emph{On a problem of {A}lexandroff concerning the dimension of
  product spaces. {I}}, J. Math. Soc. Japan \textbf{10} (1958), 380--404.
  \MR{21 \#5198}

\bibitem{Ko2}
\bysame, \emph{Test spaces for homological dimension}, Duke Math. J.
  \textbf{29} (1962), 41--50. \MR{25 \#4532}

\bibitem{Ku}
V.I. Kuz$'$minov, \emph{Homological dimension theory}, Uspehi Mat. Nauk
  \textbf{23} (1968), no.~5 (143), 3--49. \MR{39 \#2158}

\bibitem{Mi}
Haynes Miller, \emph{The {S}ullivan conjecture on maps from classifying
  spaces}, Ann. of Math. (2) \textbf{120} (1984), no.~1, 39--87. \MR{85i:55012}

\bibitem{Na}
K.~Nagami, \emph{Dimension theory}, Pure and Applied Mathematics, vol.~37,
  Academic Press, New York, 1970. \MR{42 \#6799}

\bibitem{Ol}
W.~Olszewski, \emph{Completion theorem for cohomological dimensions}, Proc.
  Amer. Math. Soc. \textbf{123} (1995), no.~7, 2261--2264. \MR{95k:54064}

\bibitem{P}
L.S. Pontryagin, \emph{Sur une hypothese fondamentale de la dimension}, C.R.
  Acad. Sci. \textbf{190} (1930), 1105--1107.

\bibitem{Sk}
E.G. Skljarenko, \emph{On the definition of cohomology dimension}, Dokl. Akad.
  Nauk SSSR \textbf{161} (1965), 538--539. \MR{31 \#739}

\bibitem{WA}
J.J. Walsh, \emph{Dimension, cohomological dimension, and cell-like mappings},
  Shape theory and geometric topology (Dubrovnik, 1981), Lecture Notes in
  Math., vol. 870, Springer, Berlin, 1981, pp.~105--118. \MR{83a:57021}

\bibitem{Wi}
R.F. Williams, \emph{A useful functor and three famous examples in topology},
  Trans. Amer. Math. Soc. \textbf{106} (1963), 319--329. \MR{26 \#4352}

\end{thebibliography}

\providecommand{\bysame}{\leavevmode\hbox to3em{\hrulefill}\thinspace}
\providecommand{\MR}{\relax\ifhmode\unskip\space\fi MR }
\providecommand{\MRhref}[2]{%
  \href{http://www.ams.org/mathscinet-getitem?mr=#1}{#2}
}
\providecommand{\href}[2]{#2}

\end{document}